\documentclass[10pt,a4paper,reqno]{amsart}
\usepackage{latexsym,amssymb,amsthm}

\theoremstyle{plain}
\newtheorem{theo+}           {Theorem}      [subsection]
\newtheorem{theo++}           {Theorem}      [section]
\newtheorem{prop+}  [theo+]  {Proposition}
\newtheorem{prop++}  [theo++]  {Proposition}
\newtheorem{coro+}  [theo+]  {Corollary}
\newtheorem{coro++}  [theo++]  {Corollary}
\newtheorem{lemm+}  [theo+]  {Lemma}
\newtheorem{lemm++}  [theo++]  {Lemma}
\newtheorem{problem}   {Problem}

\theoremstyle{definition}
\newtheorem{rema+}  [theo+]  {Remark}
\newtheorem{rema++}  [theo++]  {Remark}
\newtheorem{definition}  {Definition}

\newenvironment{theorem}{\begin{theo+}}{\end{theo+}}
\newenvironment{proposition}{\begin{prop+}}{\end{prop+}}
\newenvironment{corollary}{\begin{coro+}}{\end{coro+}}
\newenvironment{lemma}{\begin{lemm+}}{\end{lemm+}}
\newenvironment{remark}{\begin{rema+}}{\end{rema+}}

\numberwithin{equation}{section}

\begin{document}

\title[Dualities and Vertex Operator Algebras]{Dualities and
 Vertex Operator Algebras of Affine Type} 
\author[Julius Borcea]{Julius Borcea}
\address{Department of Mathematics, Stockholm University, SE-106 91, 
Stockholm, Sweden.} 
\email{julius@matematik.su.se}

\begin{abstract}
  We notice that for any positive integer $k$, the set of 
  $(1,2)$-specialized characters of level
  $k$ standard $A_{1}^{(1)}$-modules is the same as the set of
  rescaled graded dimensions of the subspaces of level $2k+1$
  standard $A_{2}^{(2)}$-modules that are vacuum spaces for the
  action of the principal Heisenberg subalgebra of $A_{2}^{(2)}$. We
  conjecture the existence of a semisimple category induced by the
  ``equal level'' representations of some algebraic structure which
  would naturally explain this duality-like property, and we study
  potential such structures in the context of generalized vertex operator 
  algebras. 
\end{abstract}

\maketitle

\section{Introduction}

In this paper we discuss a duality-like property for the rank two affine Lie
algebras $A_{1}^{(1)}$ and $A_{2}^{(2)}$, whose representation
theories are to a certain extent prototypical for the representation
theory of all the untwisted respectively twisted affine Lie
algebras. These two algebras have therefore been intensely studied
 during the past three decades (cf., e.g., \cite{Ca}, \cite{LM}, 
\cite{LW}, \cite{MP1}, \cite{MP2}).

By abuse of notations we let $\Lambda_{0}$ and $\Lambda_{1}$ denote the
fundamental weights for both algebra $A_{1}^{(1)}$ and algebra $A_{2}^{(2)}$. 
Let also $F_{(s_{1},s_{2})}$ be the $(s_{1},s_{2})$-specialization 
homomorphism (see \S 2.1 for notations). We start by noticing that for any 
positive integer $k$ the number of level $k$ standard $A_{1}^{(1)}$-modules 
is the same as the number of level $2k+1$ standard $A_{2}^{(2)}$-modules, 
namely $k+1$. The central theme of the paper is then the following 1-1
correspondence between the set of level $k$ standard $A_{1}^{(1)}$-modules 
$\mathcal{O}_{1}(k):=\big\{L_{k_{0},k}(A_{1}^{(1)})\mid 
k_{0}\in \{0,1,\ldots,k\} \big\}$ and the set 
$\big\{L_{k_{0},2k+1}(A_{2}^{(2)})\mid k_{0}\in \{0,1,\ldots,k\}
\big\}$ of level $2k+1$ standard $A_{2}^{(2)}$-modules:
\begin{align}
  F_{(1,2)}\big(e^{-\left(k_{0}\Lambda_{0}+(k-k_{0})\Lambda_{1}\right)}
\mbox{ch}\, & L_{k_{0},k}(A_{1}^{(1)})\big)\label{1.1} \\
  & =P(q)^{-1}
  F_{(1,1)}\big(e^{-(k_{0}\Lambda_{0}+(2k+1-2k_{0})\Lambda_{1})}\mbox{ch}\, 
L_{k_{0},2k+1}(A_{2}^{(2)})\big),\nonumber
\end{align}
where $k_{0}\in \{0,1,\ldots,k\}$ and $P(q)=\prod_{n\equiv \pm 1 \bmod 6}
(1-q^{n})^{-1}\in \mathbf{Z}[[q]]$.

The identities in \eqref{1.1} are obtained in \S 2.2 by means of the 
Lepowsky-Wakimoto product formulas for suitable specializations of
the Weyl-Kac character formula for standard $A_{1}^{(1)}$- and
$A_{2}^{(2)}$-modules respectively. These product formulas are 
consequences of the techniques developed in \cite{L1}, \cite{LM}, \cite{W}, 
which were based mainly on the fact that the generalized Cartan matrices 
(GCMs) of $A_{1}^{(1)}$ and $A_{2}^{(2)}$ are conjugate. However, 
there is reason to believe that the identities in \eqref{1.1} are 
in fact consequences of a  
representation-theoretical phenomenon that lies deeper than initially 
indicated by the above-mentioned conjugation property and product formulas. 
Indeed, the correspondences in \eqref{1.1} may actually hint at a new kind of 
duality between appropriately defined module categories for $A_{1}^{(1)}$ and
$A_{2}^{(2)}$, such as the category generated by the standard
$A_{1}^{(1)}$-modules of level $k$ and the category generated by the 
standard $A_{2}^{(2)}$-modules of level $2k+1$. These categories may
be equivalently defined as the categories of level $k$,
respectively level $2k+1$, integrable modules from the BGG category
$\mathcal{O}$ for $A_{1}^{(1)}$ and $A_{2}^{(2)}$ respectively 
(cf.~\cite{K1}). The objects in the category generated by the standard 
modules of level $l$ for an untwisted affine Lie
algebra $\hat{\mathfrak{g}}$ are finite direct sums of 
standard $\hat{\mathfrak{g}}$-modules of level $l$. This category is 
particularly important from the viewpoint of WZNW models in conformal field
theory and related mathematics (see, e.g., \cite{DL1}, \cite{EFK}, \cite{HL}).

As we note at the end at \S 2.2, purely Lie-algebraic considerations alone do 
not seem sufficient for finding satisfactory representation-theoretical 
explanations for the correspondences in \eqref{1.1}. We therefore propose a 
different yet natural approach, namely by using the representation
theory of (generalized) vertex operator algebras ((G)VOAs).

Let $\sigma$ be the principal automorphism of $\mathfrak{sl}(3,\mathbf{C})$ 
(see \cite{LW}) and denote by $\widehat{\mathfrak{sl}(3,\mathbf{C})}[\sigma]$ 
the principal realization of $A_{2}^{(2)}$ viewed as a
$\tfrac{1}{6}\mathbf{Z}$-graded affine Kac-Moody Lie algebra. Then the
twisted affinization $\hat{\mathfrak{s}}[\sigma]$ of the principal Cartan 
subalgebra (CSA) $\mathfrak{s}$ of $\mathfrak{sl}(3,\mathbf{C})$ is 
the principal Heisenberg subalgebra of
$\widehat{\mathfrak{sl}(3,\mathbf{C})}[\sigma]$, whose representation
theory explains in a natural way the presence of the factor $P(q)$ in
\eqref{1.1}. Indeed, the graded dimension of the twisted Fock space
representation at level $2k+1$ of $\hat{\mathfrak{s}}[\sigma]$ is
exactly $P(q^{1/6})$. Therefore, if $\Omega_{k_{0},2k+1}$
denotes the vacuum subspace of
$L_{k_{0},2k+1}\big(\widehat{\mathfrak{sl}(3,\mathbf{C})}[\sigma]\big)$
for the action of $\hat{\mathfrak{s}}[\sigma]$ and
$\mathcal{O}_{2}(k):=\big\{\Omega_{k_{0},2k+1} \mid 
k_{0}\in \{0,1,\ldots,k\} \big\}$, we get a bijection $\phi_{k}:\,
\mathcal{O}_{1}(k) \longrightarrow \mathcal{O}_{2}(k)$ such that
$\phi_{k}\big(L_{k_{0},k}(A_{1}^{(1)})\big)=\Omega_{k_{0},2k+1}$ and
\eqref{1.1} may be rewritten as
\begin{equation}\label{1.2}
F_{(1,2)}\big(e^{-(k_{0}\Lambda_{0}+(k-k_{0})\Lambda_{1})}\mbox{ch}\, 
L_{k_{0},k}(A_{1}^{(1)})\big)\Big|_{q\rightarrow  
  q^{\frac{1}{6}}}=\dim_{*}\phi_{k}\big(L_{k_{0},k}(A_{1}^{(1)})\big),\, 
k_{0}\in \{0,1,\ldots,k\},
\end{equation}
where $\dim_{*}$ denotes the graded dimension (cf.~Theorem 2.2.1). From the
viewpoint of representation theory, it is natural to ask whether
$\mathcal{O}_{1}(k)$ and $\mathcal{O}_{2}(k)$ are the same 
set of (inequivalent) simple graded objects in a semisimple category 
induced by the ``equal level'' representations of some algebraic structure 
(Problem 1). A weaker version of this question would be to ask 
whether there exists an algebraic structure for which the 
spaces $L_{k_{0},k}(A_{1}^{(1)})$ and $\Omega_{k_{0},2k+1}$ are isomorphic 
irreducible graded modules for all $k_{0}\in \{0,1,\ldots,k\}$ (Problem
  2). We conjecture that the answer to both these questions is 
affirmative, and we study various possibilities for such a structure.

Let $\nu$ be the $(1,2)$-specialization automorphism of 
$\mathfrak{sl}(2,\mathbf{C})$. It is known that
$V_{k}(A_{1}^{(1)}):=L(k\Lambda_{0};A_{1}^{(1)})$ is a $\nu$-rational VOA and 
that $L_{k_{0},k}\big(\widehat{\mathfrak{sl}(2,\mathbf{C})}[\nu]\big)$, 
$k_{0}\in \{0,1,\ldots,k\}$, are all its simple $\nu$-twisted modules 
(cf.~Theorem 3.1.5). Although there are several 
ways of interpreting appropriate modifications of the left-hand side of 
\eqref{1.1}, the most natural one appears to be that of the $q$-trace
$\chi_{_{k_{0},k}}^{\nu}(q)$ of the simple $\nu$-twisted 
$V_{k}(A_{1}^{(1)})$-module $L_{k_{0},k}\big(\widehat{\mathfrak{sl}(2,
\mathbf{C})}[\nu]\big)$ (Theorem 3.2.1).

In order to get similar interpretations of the right-hand side of \eqref{1.1},
we first construct a certain VOA $\Omega_{2k+1}^{0}\subset
L\big((2k+1)\Lambda_{0};A_{2}^{(1)}\big)$. This is done by means of
the coset (or commutant) construction associated to the irreducible
Fock space representation $M(2k+1)\subset V_{2k+1}(A_{2}^{(1)}):=
L\big((2k+1)\Lambda_{0};A_{2}^{(1)}\big)$ of the Heisenberg algebra
$\tilde{\mathfrak{s}}'\subset A_{2}^{(1)}$. The VOA
$\Omega_{2k+1}^{0}$ is actually a subspace of the vacuum space
$\Omega_{2k+1}\subset V_{2k+1}(A_{2}^{(1)})$ for
the action of $\tilde{\mathfrak{s}}'$, and the vacuum spaces
$\Omega_{k_{0},2k+1}$ defined above are $\sigma$-twisted
$\Omega_{2k+1}^{0}$-modules. When modified appropriately, the right-hand side 
of \eqref{1.1} becomes the $q$-trace $f_{k_{0},2k+1}(q)$ of the
$\sigma$-twisted $\Omega_{2k+1}^{0}$-module $\Omega_{k_{0},2k+1}$
(Theorem 3.2.5). Identity \eqref{1.2} may then be written as
$\chi_{_{k_{0},k}}^{\nu}(q^{1/2})=f_{k_{0},2k+1}(q)$, so that in fact the
equality is achieved only up to the transformation $q\rightarrow
q^{1/2}$. In addition to that, both the ranks of the VOAs 
$V_{k}(A_{1}^{(1)})$ and $\Omega_{2k+1}^{0}$ and the orders of the
automorphisms $\nu \in \text{Aut}\big(V_{k}(A_{1}^{(1)})\big)$ and
$\sigma \in \text{Aut}\big(\Omega_{2k+1}^{0}\big)$ differ by the same
factor 2. By using the recently developed permutation orbifold theory
(\cite{BDM}), we can arrange to remove these differences. We endow
$V_{k}(A_{1}^{(1)})^{\otimes 2}$ with a suitable VOA structure such
that $V_{k}(A_{1}^{(1)})^{\otimes 2}$ is $\tau$-rational and the level
$k$ standard $A_{1}^{(1)}$-modules are all its simple $\tau$-twisted
modules, where $\tau$ is a certain order 6 automorphism of
$V_{k}(A_{1}^{(1)})^{\otimes 2}$. Then the $q$-trace
$\chi_{_{k_{0},k}}^{\tau}(q)$ of the irreducible $\tau$-twisted
$V_{k}(A_{1}^{(1)})^{\otimes 2}$-module $L_{k_{0},k}(A_{1}^{(1)})$
coincides with $f_{k_{0},2k+1}(q)$ and
$\text{rank}\,V_{k}(A_{1}^{(1)})^{\otimes
  2}=\text{rank}\,\Omega_{2k+1}^{0}$ (Theo\-rem 3.2.7). Obviously, if 
$\Omega_{2k+1}^{0}$ were to fulfill the requirements of Problem 1 
then it would have to satisfy a rationality-like condition. 
The fact that $\Omega_{2k+1}^{0}$ is not actually rational 
(cf.~Remark 3.2.4 (i)) would therefore suggest that 
$\Omega_{2k+1}^{0}$ may be more appropriate for a solution to 
Problem 2 rather than Problem 1. However, this doesn't necessarily 
rule out $\Omega_{2k+1}^{0}$ as a potential solution to Problem 1 as    
a VOA $V$ may be $\sigma$-rational for some $\sigma \in \text{Aut}(V)$ 
without being rational itself (cf.~\cite{DN}). 

Several questions arise naturally at this stage. One of these is whether 
$\Omega_{k_{0},2k+1}$, $k_{0}=0,1,\ldots,k$, are irreducible 
$\sigma$-twisted $\Omega_{2k+1}^{0}$-modules. Although the results of 
\cite{DM} seem to indicate 
that this is actually true, for our purposes one would still have to deal 
with the question whether $\Omega_{2k+1}^{0}$ and 
$V_{k}(A_{1}^{(1)})^{\otimes 2}$ are related
by some VOA homomorphism.  As we point out in the Appendix 
(for the case $k=1$), these VOAs are definitely
not isomorphic. Besides, the fact that $V_{k}(A_{1}^{(1)})^{\otimes
  2}$ contains weight one vectors while $\Omega_{2k+1}^{0}$ does not
is also rather inconvenient in this context. We choose therefore to
embed $\Omega_{2k+1}^{0}$ into a larger structure, namely a simple 
GVOA in the sense of \cite{DL1} which we denote by
$\Omega_{2k+1}^{A}$ (Theorem 4.4.3). This is done in Section 4. 
In addition to the fact that 
it certainly contains weight one vectors, we can also prove that 
$\Omega_{2k+1}^{A}$ acts irreducibly on the spaces $\Omega_{k_{0},2k+1}$ 
without altering the $q$-traces $f_{k_{0},2k+1}(q)$ (Theorem 4.4.8). 
As we explain at the end of \S 4.4, the GVOA $\Omega_{2k+1}^{A}$ 
appears to be a more appropriate object than the VOA $\Omega_{2k+1}^{0}$ 
for further investigating Problem 1. Actually, the techniques that we 
use in Section 4 seem in fact flexible enough to allow the construction of 
structures even larger than $\Omega_{2k+1}^{A}$ that still satisfy the same
properties (cf.~Remark 4.4.9). Moreover, these techniques are easily adapted 
to a more general situation, and they may therefore be useful for an 
axiomatic study of the (yet-to-be-defined) notion of ``twisted module'' 
for a GVOA. The appropriate axiomatic setting once developed, the spaces 
$\Omega_{k_{0},2k+1}$ should provide natural examples of simple 
``twisted modules'' for the GVOA $\Omega_{2k+1}^{A}$.

The construction of $\Omega_{2k+1}^{A}$ and its action on
$\Omega_{k_{0},2k+1}$ are based on the theory of relative 
vertex operators and the $\mathcal{Z}$-algebra theory developed in 
\cite{DL1}, \cite{DL2}, \cite{LW}. It would be interesting 
to see whether structures and techniques such as simple current 
extensions of VOAs and intertwining operator algebras -- introduced in 
\cite{DL1}, \cite{DL2}, \cite{DLiM3}, \cite{HL}, \cite{Li3}, \cite{Li4} -- 
provide additional tools for a continuation of the present work. Finally, 
one should also investigate whether similar duality-like properties hold 
for affine Lie algebras of higher ranks.

This paper is organized as follows: in \S 2.1 we recall some basic facts 
about affine Lie algebras; in \S 2.2 we prove \eqref{1.2} and formulate 
the problems that we study in Sections 3 and 4. In \S 3.1 we review 
the results on GVOAs and their 
modules that we need for \S 3.2, where we construct the VOAs 
$V_{k}(A_{1}^{(1)})$, $\Omega_{2k+1}^{0}$, and 
$V_{k}(A_{1}^{(1)})^{\otimes 2}$ and reinterpret \eqref{1.2} in terms of 
the representation theory of these VOAs. In \S 4.1-\S 4.3 we change the 
setting by adapting the theories developed in  
\cite{DL1} and \cite{DL2} to our particular case. In \S 4.4 we use this new 
setting to embed $\Omega_{2k+1}^{0}$ into the 
GVOA $\Omega_{2k+1}^{A}$ and to prove that the latter acts irreducibly on the 
 spaces $\Omega_{k_{0},2k+1}$ without altering their $q$-traces.

Throughout this paper we shall use $\mathbf{N}$ 
for the nonnegative integers, $\mathbf{Z}_{+}$ for the positive integers, 
$\mathbf{Q}_{\geq 0}$ for the nonnegative rational numbers, 
and $\mathbf{C}^{\times}$ for the nonzero complex numbers.

\noindent
{\bf Acknowledgements.} I would like to thank Chongying Dong, Mirko Primc, and 
Arne Meurman for useful suggestions and discussions, and the referee for his 
patience and constructive criticism.

\numberwithin{equation}{subsection}

\section{A duality-like property for rank two affine Lie algebras}

We first introduce a few notations and review some standard material on 
affine Lie algebras and their 
highest weight modules (further details may be found in 
e.g.~\cite{K1}, \cite{L1}). We then concentrate on the rank two affine Lie 
algebras and formulate the problems that we study in Sections 3 and 4.

\subsection{Preliminaries and notations}

Let $\mathfrak{g}$ be a finite-dimensional simple complex Lie algebra. Fix a 
CSA $\mathfrak{t}$ of $\mathfrak{g}$ and let $\mu$ be an 
automorphism of $\mathfrak{g}$ of order $r$ ($=1,2,$ or 3) induced by
an automorphism of order $r$ of the Dynkin diagram of $\mathfrak{g}$
with respect to $\mathfrak{t}$. Let
$\varepsilon$ be a primitive $r$-th root of unity and denote by
$\mathfrak{g}_{[i]}$ the $\varepsilon^{i}$-eigenspace of $\mu$ in
$\mathfrak{g}$, $i\in \mathbf{Z}_{r}$. Then the fixed set 
$\mathfrak{g}_{[0]}$ is a simple subalgebra of
$\mathfrak{g}$, the space $\mathfrak{b}:=\mathfrak{g}_{[0]}\cap \mathfrak{t}$ 
is a CSA of $\mathfrak{g}_{[0]}$, and the $\mathfrak{g}_{[0]}$-modules
$\mathfrak{g}_{[1]}$ and $\mathfrak{g}_{[-1]}$ are irreducible and
contragredient. Set $l=\text{rank}\,\mathfrak{g}_{[0]}$, and let 
$\{\beta_{1},\ldots,\beta_{l}\}\subset \mathfrak{b}^{*}$ be a root basis of
$\mathfrak{g}_{[0]}$ and $\{E_{j}, F_{j}, H_{j}\mid j\in
\{1,\ldots,l\}\}$ a corresponding set of
canonical generators of $\mathfrak{g}_{[0]}$. Let $\beta_{0}\in
\mathfrak{b}^{*}$ be the lowest weight of the $\mathfrak{g}_{[0]}$-module
$\mathfrak{g}_{[1]}$, and let $E_{0}$ and $F_{0}$
be a lowest weight vector of the $\mathfrak{g}_{[0]}$-module
$\mathfrak{g}_{[1]}$ respectively a highest weight vector of the 
$\mathfrak{g}_{[0]}$-module $\mathfrak{g}_{[-1]}$. We assume that $E_{0}$ 
and $F_{0}$ are normalized so that $[H_{0},E_{0}]=2E_{0}$,
where $H_{0}=[E_{0},F_{0}]$. For $i,j\in \{0,1,\ldots,l\}$ define
$a_{ij}=\beta_{j}(H_{i})$ and denote by $A$ the matrix
$(a_{ij})_{i,j=0}^{l}$. Then $A$ is an affine GCM, and so there are positive
integers $a_{0},\ldots,a_{l}$ such that $(a_{0},\ldots,a_{l})A^{t}=0$. 
Equivalently, there exist positive integers
$a\check{}_{0},\ldots,a\check{}_{l}$ such that
$(a\check{}_{0},\ldots,a\check{}_{l})A=0$. Both these sets of
integers are assumed to be normalized so that
$\gcd(a_{0},\ldots,a_{l})=\gcd(a\check{}_{0},\ldots,a\check{}_{l})$ $=1$. The
integers $h:=\sum_{j=0}^{l}a_{j}$ and
$h\check{}:=\sum_{j=0}^{l}a\check{}_{j}$ are then the Coxeter
number respectively the dual Coxeter
number of the matrix $A$. If $\mathfrak{g}$ is of type $X_{N}$
($X=A,B,\ldots,G$ and $N\ge 1$), then $A$ will be denoted by
$X_{N}^{(r)}$. We shall use the Dynkin diagrams of the affine GCMs as
listed in \cite{KKLW}, that is, with the vertex corresponding to the 0-th
index always occurring at the left-end of the diagram. The $a_{j}$'s
are then the numerical labels next to the vertices of each diagram. In
particular, $a_{0}$ is always 1, while $a\check{}_{0}$ is 1 unless $A$
is of type $A_{2l}^{(2)}$ ($l\ge 1$), in which case
$a\check{}_{0}=2$.

Let $\mathbf{s}=(s_{0},s_{1},\ldots,s_{l})$ be a sequence of
nonnegative relatively prime integers, and set
$T=r\sum_{j=0}^{l}s_{j}a_{j}$. If $\eta$ is a primitive $T$-th root of
unity, the conditions
\begin{equation}\label{2.1.1}
  \nu(H_{j})=H_{j}, \; \nu(E_{j})=\eta^{s_{j}}E_{j}, \; j\in
  \{0,1,\ldots,l\},
\end{equation}
define a $T$-th order automorphism $\nu$ of
$\mathfrak{g}$, the so-called $\mathbf{s}$-automorphism. The 
$(1,0,\ldots,0)$-automor\-phism is just the 
original diagram automorphism $\mu$. For $j\in \mathbf{Z}_{T}$ 
let $\mathfrak{g}_{(j)}$ be the 
$\eta^{j}$-eigenspace of $\nu$ in $\mathfrak{g}$ ($j$ denotes here both an 
integer between 0 and $T-1$ and its residue class modulo $T$). The
$\mathbf{Z}_{T}$-gradation $\mathfrak{g}=\coprod_{j\in
  \mathbf{Z}_{T}}\mathfrak{g}_{(j)}$ is accordingly called the
$\mathbf{s}$-gradation. Let $\langle \cdot\, ,\cdot \rangle$ be a nondegenerate
symmetric $\mathfrak{g}$-invariant bilinear form on
$\mathfrak{g}$. Being a multiple of the Killing form,
$\langle \cdot\, ,\cdot \rangle$ is also $\nu$-invariant and 
remains nonsingular on
the CSA $\mathfrak{b}$ of $\mathfrak{g}_{[0]}$. We may
therefore identify $\mathfrak{b}$ with $\mathfrak{b}^{*}$ by means of the
restricted form. Furthermore, we may assume that 
$\langle \cdot\, ,\cdot \rangle$
is normalized so that $\langle \beta_{0},\beta_{0} \rangle =2a\check{}_{0}/r$,
which then implies that
$a\check{}_{j}=r\langle \beta_{j},\beta_{j} \rangle a_{j}/2$ 
for $j=0,1,\ldots,l$
(cf. ~\cite[Proposition 1.1]{KKLW}). This normalization of the form
$\langle \cdot\, ,\cdot \rangle$ amounts to the condition that
$\langle \alpha,\alpha \rangle =2$ whenever $\alpha\in \mathfrak{t}^{*}$ is a 
long root of $\mathfrak{g}$, in which case
the Killing form  equals $2h\check{}\langle \cdot\, ,\cdot \rangle$. 
Define the Lie algebras
\begin{equation}\label{2.1.2}
  \hat{\mathfrak{g}}[\nu]=\oplus_{j=0}^{T-1}\mathfrak{g}_{(j)}\otimes
  t^{\frac{j}{T}}\mathbf{C}[t,t^{-1}]\oplus \mathbf{C}c, \quad  
  \tilde{\mathfrak{g}}[\nu]= \hat{\mathfrak{g}}[\nu]\rtimes
  \mathbf{C}d,
\end{equation}
by the conditions
\begin{equation}\label{2.1.3}
\begin{split}
  & c \;\, \mbox{central},\; c\neq 0, \; [d,a\otimes t^{m}]=m\, a\otimes
  t^{m},\\
  & [a\otimes t^{m},b\otimes t^{n}]=[a,b]\otimes
  t^{m+n}+m\delta_{m+n,0}\langle a,b \rangle c,
\end{split}
\end{equation}
for $m,n\in \frac{1}{T}\mathbf{Z}$, $a\in \mathfrak{g}_{(mT \bmod T)}$,
$b\in \mathfrak{g}_{(nT \bmod T)}$, where $p \bmod T$ denotes the residue 
class of an integer $p$ modulo $T$. For $a\in \mathfrak{g}_{(j)}$, $n\in
\mathbf{Z}$, we shall frequently use $a(n+j/T)$ and
$\mathfrak{g}(n+j/T)$ to denote $a\otimes t^{n+j/T}$ and 
$\mathfrak{g}_{(j)}\otimes
  t^{n+j/T}$, respectively, and we often identify
  $\mathfrak{g}_{(0)}(0)$ with $\mathfrak{g}_{(0)}$. The space
  $\mathfrak{h}:=\mathfrak{b}\oplus \mathbf{C}c$ (respectively, 
  $\mathfrak{h}^{e}:=\mathfrak{h}\rtimes \mathbf{C}d$) is a CSA of 
$\hat{\mathfrak{g}}[\nu]$ (respectively, $\tilde{\mathfrak{g}}[\nu]$).
 Let $\delta\in \mathfrak{h}^{e^{*}}$ be such that
$\delta|_{\mathfrak{h}}=0$, $\delta(d)=1$, and define 
$\alpha_{j}\in \mathfrak{h}^{e^{*}}$ by
$\alpha_{j}|_{\mathfrak{b}}=\beta_{j}$, $\alpha_{j}(c)=0$,
$\alpha_{j}(d)=s_{j}T^{-1}$ if $j=1,\ldots,l$, and
$\alpha_{0}=r^{-1}\delta-\sum_{j=1}^{l}a_{j}\alpha_{j}$, so that in
particular $\alpha_{0}(d)=s_{0}T^{-1}$. For $j\in \{0,1,\ldots,l\}$ 
define also 
\begin{equation}\label{2.1.4}
  e_{j}=E_{j}\otimes t^{\frac{s_{j}}{T}},\, f_{j}=F_{j}\otimes
  t^{-\frac{s_{j}}{T}},\,
  h_{j}=H_{j}+\frac{2s_{j}}{T\langle \beta_{j},\beta_{j} \rangle}c.
\end{equation}
Then $\{e_{j},f_{j},h_{j},d\mid j\in \{0,1,\ldots,l\}\}$ is a system
of canonical generators of $\tilde{\mathfrak{g}}[\nu]$, viewed as the
$\nu$-twisted affine Kac-Moody Lie algebra of rank $l+1$ associated to 
the GCM $A$, in what is called the $\mathbf{s}$-realization of this algebra. 
Note that the canonical central element 
$c$ equals $\sum_{j=0}^{l}a\check{}_{j}h_{j}$, and let
$\tilde{\mathfrak{g}}[\nu]_{i}:=\{a\in \tilde{\mathfrak{g}}[\nu]\mid
[d,a]=i\,a\}$, $i\in \frac{1}{T}\mathbf{Z}$. The corresponding 
$\frac{1}{T}\mathbf{Z}$-gradation
$\tilde{\mathfrak{g}}[\nu]=\coprod_{i\in
  \frac{1}{T}\mathbf{Z}}\tilde{\mathfrak{g}}[\nu]_{i}$ is then called the
$\mathbf{s}$-gradation of $\tilde{\mathfrak{g}}[\nu]$.
\begin{remark}
  (i) Let us temporarily denote by
  $\{e^{\nu}_{j},f^{\nu}_{j},h^{\nu}_{j},d^{\nu}\mid j\in
  \{0,1,\ldots,l\}\}$ the canonical generators of
  $\tilde{\mathfrak{g}}[\nu]$ defined in \eqref{2.1.4}. One can show that
  there are uniquely determined scalars $x_{1},\ldots,x_{l}$
  such that the affine Lie algebras $\tilde{\mathfrak{g}}[\mu]$ and
  $\tilde{\mathfrak{g}}[\nu]$ are isomorphic (but of course not graded
  isomorphic) by the Lie algebra map defined by $c\mapsto c$, $d^{\mu}\mapsto
    d^{\nu}+\sum_{k=1}^{l}x_{k}H_{k}$, $e^{\mu}_{j}\mapsto
    e^{\nu}_{j}$, $f^{\mu}_{j}\mapsto
    f^{\nu}_{j}$, $h^{\mu}_{j}\mapsto h^{\nu}_{j}$, $j=0,1,\ldots,l$. \newline
(ii) There are two particularly important realizations of an affine Lie
algebra, namely the homogeneous and the principal ones,
given by the choices $\mathbf{s}=(1,0,\ldots,0)$ and
$\mathbf{s}=(1,1,\ldots,1)$ respectively.
\end{remark}

Recall that a $\hat{\mathfrak{g}}[\nu]$-module $V$ is said to be restricted 
if $\mathfrak{g}_{(j)}(n+j/T)\cdot v=0$ for any
$v\in V$ and $n\gg 0$. In particular, any highest weight module is restricted.
 Given a restricted $\hat{\mathfrak{g}}[\nu]$-module 
$V$ and $a\in
\mathfrak{g}_{(j)}$, we shall use generating functions of
operators on $V$ of the form 
\begin{equation}\label{2.1.5}
  a(\nu;z)=\sum_{n\in \mathbf{Z}}a(n+j/T)z^{-n-\frac{j}{T}-1}\in
  (\mbox{End}\,V)\big[\big[z^{1/T},z^{-1/T}\big]\big].
\end{equation}
We shall often write $a_{n+j/T}$ when we think of $a(n+j/T)$ as a
coefficient of $a(\nu;z)$, and $a(\mbox{id}_{\mathfrak{g}};z)$ will be
denoted by $a(z)$. We shall mostly consider 
highest weight modules with highest weight $\Lambda\in \mathfrak{h}^{e^{*}}$ 
satisfying $\Lambda(d)=0$. The canonical central element $c\in
\tilde{\mathfrak{g}}[\nu]$ acts on a highest weight module $V$ with
highest weight $\Lambda$ as multiplication by the scalar
$\Lambda(c)$, called the level of
$V$. The character of $V$ is defined as the formal infinite sum 
$\mbox{ch}V=\sum_{\lambda\in \mathfrak{h}^{e^{*}}}(\mbox{dim}\,
  V_{\lambda})e^{\lambda}$, where $V=\coprod_{\lambda\le \Lambda}V_{\lambda}$ 
is the weight decomposition. Note that $e^{-\Lambda}\mbox{ch}V\in 
\mathbf{Z}[[e^{-\alpha_{0}},\ldots,e^{-\alpha_{l}}]]$ and let
$q$ be an indeterminate. Provided that $s_{i}>0$, 
$0\le i\le l$, the sequence $\mathbf{s}=(s_{0},\ldots,s_{l})$
defines a homomorphism of power series rings
\begin{align*}
  F_{\mathbf{s}}:\mathbf{Z}[[e^{-\alpha_{0}},\ldots,e^{-\alpha_{l}}]]
  & \longrightarrow
\mathbf{Z}[[q]]\\
  e^{-\alpha_{i}} & \longmapsto
  F_{\mathbf{s}}(e^{-\alpha_{i}})=q^{s_{i}}, \quad i=0,\ldots,l,
\end{align*}
called the $q$-specialization of type $\mathbf{s}$. Let $j\in 
\mathbf{N}$ and set $V_{j}(\mathbf{s})=\coprod_{\lambda:\;
  \deg \lambda=j}V_{\lambda}$,  where 
$\deg\big(\Lambda-\sum_{i=0}^{l}k_{i}\alpha_{i}\big):=
\sum_{i=0}^{l}k_{i}s_{i}$, $(k_{0},\ldots,k_{l})\in 
\mathbf{N}^{l+1}$. Then 
$V=\coprod_{j\in \mathbf{N}}V_{j}(\mathbf{s})$ is the 
  $\mathbf{s}$-gradation of V and one has  
$F_{\mathbf{s}}(e^{-\Lambda}\mbox{ch}V)=\sum_{ j\in
  \mathbf{N}}(\mbox{dim}\,V_{j}(\mathbf{s}))q^{j}$.

For $i=0,1,\ldots,l$, let 
$\Lambda_{i}\in P_{+}$ be the fundamental weight determined by
$\Lambda_{i}(h_{j})=\delta_{ij}$, $\Lambda_{i}(d)=0$. The standard
$\tilde{\mathfrak{g}}[\nu]$-modules with highest weight $\Lambda$ such
that $\Lambda(d)=0$ are then parametrized up to equivalence by arbitrary
sequences $(k_{0},k_{1},\ldots,k_{l})\in \mathbf{N}^{l+1}$ 
satisfying $\Lambda=\sum_{i=0}^{l}k_{i}\Lambda_{i}$. We shall
interchangeably use the notations
$L_{k_{0},\ldots,k_{l-1},k}(\tilde{\mathfrak{g}}[\nu])$,
$L(\Lambda;\tilde{\mathfrak{g}}[\nu])$, 
and $L(\Lambda)$ (if no confusion is possible) in order to designate
the level $k$ standard $\tilde{\mathfrak{g}}[\nu]$-module with highest weight
$\Lambda=\sum_{i=0}^{l}k_{i}\Lambda_{i}$. The 
$\tilde{\mathfrak{g}}[\nu]$-modules $L(\Lambda_{i})$, $i=0,1,\ldots,l$, are 
called fundamental, and the level one standard
modules are called basic. Every standard
module occurs in a tensor product of fundamental modules, and in fact
every standard module of level $k$ occurs in the
tensor product of $k$ basic modules (cf. ~\cite{Ag}).

\subsection{Rank two affine Lie algebras}

We now apply the constructions described in \S 2.1 to 
the affine Lie algebras $A_{1}^{(1)}$ and $A_{2}^{(2)}$, and
we deduce a series of coincidences between specialized characters of 
certain standard
representations of these algebras (Theorem 2.2.1). We
then describe two ways of investigating the existence of conceptual 
explanations for these coincidences (Problems 1 and 2).

Starting with $\mathfrak{g}=\mathfrak{sl}(2,\mathbf{C})$,
$\mathbf{s}=(1,0)$ and $\nu=\mbox{id}_{\mathfrak{g}}$, one gets the
``homogeneous picture'' of the (untwisted) affine Lie algebra
$\widehat{\mathfrak{sl}(2,\mathbf{C})}$, also denoted $A_{1}^{(1)}$,
whose GCM is \(
\begin{pmatrix}
2 & -2 \\
-2 & 2
\end{pmatrix}
\). Taking instead $\mathbf{s}=(1,2)$ and the corresponding 
$(1,2)$-automorphism
$\nu$ of $\mathfrak{g}$ defined as in \eqref{2.1.1}, one gets
$\widehat{\mathfrak{sl}(2,\mathbf{C})}[\nu]$, that is, the $(1,2)$-realization
of $A_{1}^{(1)}$. We shall use the standard basis $\{e,f,h\}$ of
$\mathfrak{sl}(2,\mathbf{C})$, with $[e,f]=h$, $[h,e]=2e$, and 
$[h,f]=-2f$.

Let now $\mathfrak{g}=\mathfrak{sl}(3,\mathbf{C})$ and denote by $\mu$
the negative transpose
map of $\mathfrak{g}$. Then $\mu$ is a diagram automorphism of order 2
of $\mathfrak{g}$ ($\mu$ is the outer automorphism induced by the 
automorphism of the Dynkin diagram
of $\mathfrak{g}$ -- with respect to some appropriate CSA 
of $\mathfrak{g}$ -- which permutes the two simple
roots). Note that the fixed set of $\mu$ in $\mathfrak{g}$ is the
3-dimensional rank one subalgebra
$\mathfrak{g}_{[0]}=\mathfrak{so}(3,\mathbf{C})$ ($\cong
\mathfrak{sl}(2,\mathbf{C})$), while $\mathfrak{g}_{[1]}$ is the
5-dimensional subspace consisting of the symmetric traceless $3\times
3$ matrices. The discussion of \S 2.1 with $\mathbf{s}=(1,0)$
and $\mu$ as the $(1,0)$-automorphism leads then to the homogeneous
realization of the twisted affine Lie algebra
$\widehat{\mathfrak{sl}(3,\mathbf{C})}[\mu]$, also denoted
$A_{2}^{(2)}$, whose GCM is \(
\begin{pmatrix}
2 & -1 \\
-4 & 2
\end{pmatrix}
\). Using instead $\mathbf{s}=(1,1)$ and the principal automorphism
$\sigma$ of $\mathfrak{g}$ defined as in \eqref{2.1.1}, one gets the algebra
$\widehat{\mathfrak{sl}(3,\mathbf{C})}[\sigma]$, which is the
principal realization of $A_{2}^{(2)}$.

We shall consider later on different realizations of the rank three   
untwisted affine Lie algebra $A_{2}^{(1)}$, obtained as in the
previous section from $\mathfrak{g}=\mathfrak{sl}(3,\mathbf{C})$
and $\mu=\mbox{id}_{\mathfrak{g}}$.

Let $k\in \mathbf{Z}_{+}$. Recall that the standard
$A_{1}^{(1)}$-modules with a highest weight $\Lambda\in P_{+}$ such that 
$\Lambda(c)-k=\Lambda(d)=0$ are parametrized up to equivalence 
by pairs $(k_{0},k_{1})\in \mathbf{N}^{2}$ satisfying
\begin{equation}\label{2.2.1}
  \Lambda=k_{0}\Lambda_{0}+k_{1}\Lambda_{1},\quad k_{0}+k_{1}=k.
\end{equation}
Note that there are exactly $k+1$ such distinct pairs. Suppose that 
$\Lambda\in P_{+}$ is of the form \eqref{2.2.1}. If $k_{0}+k_{1}\neq 3k_{1}+1$,
then the following infinite product
expression for the $(1,2)$-specialized character of
$L_{k_{0},k}(A_{1}^{(1)})$ is a simple consequence of identity (11.1.12) in  
\cite{MP2}:
{\allowdisplaybreaks
\begin{equation}\label{2.2.2}
  F_{(1,2)}\big(e^{-\Lambda}\mbox{ch}\,
  L_{k_{0},k}(A_{1}^{(1)})\big)=\prod_{
\genfrac{}{}{0pt}{}{n\not\equiv 0, 2(k+2),\pm (k_{0}+1),}{\pm (2k-k_{0}+3),
\pm 2(k-k_{0}+1) \bmod 4(k+2)}}(1-q^{n})^{-1}.
\end{equation}
If $k\equiv 1 \bmod 3$ and $k_{0}+k_{1}=3k_{1}+1$ then identity  
(11.1.13) in \cite{MP2} gives instead:
\begin{equation}\label{2.2.3}
  F_{(1,2)}\big(e^{-\Lambda}\mbox{ch}\,
  L_{\frac{2k+1}{3},k}(A_{1}^{(1)})\big) 
  =\prod_{n\equiv
    \pm \frac{2(k+2)}{3} \bmod 4(k+2)}(1-q^{n})
\prod_{\genfrac{}{}{0pt}{}{n\not\equiv 0,2(k+2),\pm \frac{2(k+2)}{3},}{\pm 
\frac{4(k+2)}{3} \bmod 4(k+2)}}(1-q^{n})^{-1}.
\end{equation}
The identities in \cite{MP2} mentioned above were derived from the Weyl-Kac 
character and denominator formulas (\cite{K1}) and Wakimoto's generalization 
of Lepowsky's numerator formula (\cite{L1}, \cite{W}).}

We now look at standard $A_{2}^{(2)}$-modules with a highest weight
$\Lambda\in P_{+}$ such that $\Lambda(c)=2k+1$ and
$\Lambda(d)=0$. These are parametrized up to equivalence by 
pairs $(k_{0},k_{1})\in \mathbf{N}^{2}$ satisfying
\begin{equation}\label{2.2.4}
  \Lambda=k_{0}\Lambda_{0}+k_{1}\Lambda_{1},\quad 2k_{0}+k_{1}=2k+1.
\end{equation}  
Note that there are exactly $k+1$ such distinct pairs. Define
\begin{equation*}  
  P(q)=\prod_{n\equiv \pm 1 \bmod 6}(1-q^{n})^{-1}\in \mathbf{Z}[[q]]
\end{equation*}
and let 
$\Lambda\in P_{+}$ be as in \eqref{2.2.4}. The principally specialized
character of $L(\Lambda;A_{2}^{(2)})$ was expressed as an infinite
product in \cite[Theorem 6.8]{LM} by means of the character, denominator
and numerator formulas. One has:
{\allowdisplaybreaks
\begin{equation}\label{2.2.5}
  F_{(1,1)}\big(e^{-\Lambda}\mbox{ch}\,
  L_{k_{0},2k+1}(A_{2}^{(2)})\big)=P(q)\prod_{\genfrac{}{}{0pt}{}{n\not\equiv 
0,2(k+2),\pm (k_{0}+1),}{\pm (2k-k_{0}+3),
\pm 2(k-k_{0}+1) \bmod 4(k+2)}}(1-q^{n})^{-1}
\end{equation}    
if \( k_{0}\neq \frac{2k+1}{3} \), and
\begin{align}
  F_{(1,1)}&\big(e^{-\Lambda}\mbox{ch}\,
  L_{\frac{2k+1}{3},2k+1}(A_{2}^{(2)})\big)\label{2.2.6} \\
  &=P(q)\prod_{n\equiv
    \pm \frac{2(k+2)}{3} \bmod 4(k+2)}(1-q^{n})\prod_{
    \genfrac{}{}{0pt}{}{n\not\equiv 0,2(k+2), \pm \frac{2(k+2)}{3},}{\pm 
\frac{4(k+2)}{3} \bmod 4(k+2)}}(1-q^{n})^{-1}
\nonumber
\end{align}   
if \( k_{0}=\frac{2k+1}{3} \) (which can occur only for $k\equiv 1
\bmod 3$).

It is now obvious that except for the factor $P(q)$, the right-hand
sides of \eqref{2.2.2} and \eqref{2.2.3} coincide with those of \eqref{2.2.5} 
and \eqref{2.2.6}
respectively. Therefore, for each positive integer $k$ we have a 1-1 
correspondence
\begin{equation}\label{2.2.7}
  \big\{L_{k_{0},k}(A_{1}^{(1)})\mid k_{0}\in \{0,1,\ldots,k\} \big\}
  \longleftrightarrow \big\{L_{k_{0},2k+1}(A_{2}^{(2)})\mid k_{0}\in 
\{0,1,\ldots,k\} \big\}
\end{equation}
such that for every $k_{0}\in \{0,1,\ldots,k\}$
\begin{equation*}
  F_{(1,2)}\big(e^{-\left(k_{0}\Lambda_{0}+(k-k_{0})\Lambda_{1}\right)}
\mbox{ch}\, L_{k_{0},k}(A_{1}^{(1)})\big)=P(q)^{-1} 
  F_{(1,1)}\big(e^{-(k_{0}\Lambda_{0}+(2k+1-2k_{0})\Lambda_{1})}\mbox{ch}
\, L_{k_{0},2k+1}(A_{2}^{(2)})\big).
\end{equation*}

The presence of the factor $P(q)$ in 
\eqref{2.2.5}-\eqref{2.2.6} may be explained by the 
representation theory of the principal Heisenberg subalgebra of 
$\widehat{\mathfrak{sl}(3,\mathbf{C})}[\sigma]$, where $\sigma$ is as
before the principal automorphism of
$\mathfrak{sl}(3,\mathbf{C})$.} More precisely, $P(q)$ is the rescaled
graded dimension of the fermionic twisted Fock space representation
at level $2k+1$ of this Heisenberg algebra. Indeed, let $\tau :
\widehat{\mathfrak{sl}(3,\mathbf{C})}[\sigma] \rightarrow
\widehat{\mathfrak{sl}(3,\mathbf{C})}[\sigma]/\mathbf{C}c$ be the
canonical Lie algebra map, and note that the
principal gradation of $\widetilde{\mathfrak{sl}(3,\mathbf{C})}[\sigma]$ 
is actually a
$\frac{1}{6}\mathbf{Z}$-gradation. Consider the following subalgebras
of $\widetilde{\mathfrak{sl}(3,\mathbf{C})}[\sigma]$:
\begin{equation}\label{2.2.8}
  \begin{split}
  & \tilde{\mathfrak{s}}[\sigma]=\tau
  ^{-1}\Big(\mbox{Cent}_{\widehat{\mathfrak{sl}(3,\mathbf{C})}[\sigma]/
\mathbf{C}c}\tau (e_{0}+e_{1})\Big)\rtimes \mathbf{C}d,\quad 
\tilde{\mathfrak{s}}[\sigma]'=\big[\tilde{\mathfrak{s}}[\sigma],
\tilde{\mathfrak{s}}[\sigma]\big], \\
  & \tilde{\mathfrak{s}}[\sigma]'_{\pm}=\mbox{{\bf C}-span}\, \{x\in
  \tilde{\mathfrak{s}}[\sigma]' \mid \pm\deg x>0 \},
  \end{split}  
\end{equation}
where $e_{0}$ and $e_{1}$ are as in \eqref{2.1.4}. Then $
\tilde{\mathfrak{s}}[\sigma]'=  \tilde{\mathfrak{s}}[\sigma]'_{-}\oplus
\mathbf{C}c\oplus \tilde{\mathfrak{s}}[\sigma]'_{+}$, and one can show
that the commutator subalgebra of $\tilde{\mathfrak{s}}[\sigma]'$ is
one-dimensional and coincides with $\mathbf{C}c$ (see Remark 3.2.3). Hence
$\tilde{\mathfrak{s}}[\sigma]'$ is a Heisenberg Lie algebra, called the 
principal Heisenberg subalgebra of
$\widetilde{\mathfrak{sl}(3,\mathbf{C})}[\sigma]$ (cf. ~\cite{LW}). Define 
a $(\tilde{\mathfrak{s}}[\sigma]'_{+}\oplus \mathbf{C}c)$-module structure 
on $\mathbf{C}$ by $c\cdot 1=(2k+1)\cdot 1$, $\tilde{\mathfrak{s}}[\sigma]'_{+}
\cdot 1
=0$, $\deg 1=0$, and form the induced $\tilde{\mathfrak{s}}[\sigma]'$-module
\begin{equation}\label{2.2.9}
  M(\sigma;2k+1)=U(\tilde{\mathfrak{s}}[\sigma]')
\otimes_{U(\tilde{\mathfrak{s}}
[\sigma]'_{+}\oplus \mathbf{C}c)} \mathbf{C},
\end{equation}
which acquires a $\frac{1}{6}\mathbf{Z}$-gradation by letting $d$ act
as the degree operator. Then $M(\sigma;2k+1)$ is an irreducible
$\tilde{\mathfrak{s}}[\sigma]'$-module which is 
isomorphic to $S(\tilde{\mathfrak{s}}[\sigma]'_{-})$
as a $\frac{1}{6}\mathbf{Z}$-graded vector space. Recall that the vacuum 
space of an arbitrary $\tilde{\mathfrak{s}}[\sigma]'$-module $V$ is the 
($\frac{1}{6}\mathbf{Z}$-graded) subspace $\Omega_{V}=\{v\in V 
\mid \tilde{\mathfrak{s}}[\sigma]'_{+}\cdot v=0\}$. If $\Lambda\in P_{+}$ is 
as in \eqref{2.2.4}, then \cite[Theorem 1.7.3]{FLM} implies that when 
viewed as a $\tilde{\mathfrak{s}}[\sigma]'$-module,
$L_{k_{0},2k+1}\big(\widehat{\mathfrak{sl}(3,\mathbf{C})}[\sigma]\big)$ 
decomposes as 
\begin{equation}\label{2.2.10}
  L_{k_{0},2k+1}\big(\widehat{\mathfrak{sl}(3,\mathbf{C})}[\sigma]\big)=
M(\sigma;2k+1)\otimes \Omega_{k_{0},2k+1}=S(\tilde{\mathfrak{s}}[\sigma]'_{-})
\otimes \Omega_{k_{0},2k+1},
\end{equation}
where
$\Omega_{k_{0},2k+1}:=\Omega_{L_{k_{0},2k+1}(\widehat{\mathfrak{sl}(3,
\mathbf{C})}[\sigma])}$ is the vacuum subspace of
$L_{k_{0},2k+1}\big(\widehat{\mathfrak{sl}(3,\mathbf{C})}[\sigma]\big)$ for
the action of $\tilde{\mathfrak{s}}[\sigma]'$. It is then easy to see that 
\begin{align}
  F_{(1,1)}&\big(e^{-(k_{0}\Lambda_{0}+(2k+1-2k_{0})\Lambda_{1})}\mbox{ch}\,
L_{k_{0},2k+1}(A_{2}^{(2)})\big)\Big|_{q\rightarrow q^{\frac{1}{6}}}
\label{2.2.11}\\
  &=\dim_{*}
  L_{k_{0},2k+1}\big(\widehat{\mathfrak{sl}(3,\mathbf{C})}[\sigma]\big)=
\dim_{*}S(\tilde{\mathfrak{s}}[\sigma]'_{-})\dim_{*}\Omega_{k_{0},2k+1}
  =P(q^{\frac{1}{6}})\dim_{*}\Omega_{k_{0},2k+1},\nonumber
\end{align}
where $\dim_{*}$ denotes the graded dimension (cf. ~\cite[\S 1.10]{FLM}). Set 
\begin{equation}\label{O}
  \mathcal{O}_{1}(k)=\big\{L_{k_{0},k}(A_{1}^{(1)})\mid k_{0}\in
        \{0,1,\ldots,k\}\big\},\,
  \mathcal{O}_{2}(k)=\big\{\Omega_{k_{0},2k+1} \mid
k_{0}\in \{0,1,\ldots,k\} \big\}.
\end{equation}
We conclude that \eqref{2.2.7} can be restated as follows:

{\allowdisplaybreaks
\begin{theorem}
  For each
positive integer $k$ there is a bijection 
\begin{align}\label{2.2.12}
  \phi_{k}:\, & \mathcal{O}_{1}(k) \longrightarrow
    \mathcal{O}_{2}(k)\\
    & L_{k_{0},k}(A_{1}^{(1)}) \longmapsto
    \Omega_{k_{0},2k+1},\quad k_{0}\in \{0,1,\ldots,k\},\nonumber
\end{align}
such that 
$$ F_{(1,2)}\big(e^{-(k_{0}\Lambda_{0}+(k-k_{0})\Lambda_{1})}\mbox{{\em
      ch}}\, L_{k_{0},k}(A_{1}^{(1)})\big)\Big|_{q\rightarrow  
    q^{\frac{1}{6}}}=\dim_{*}\phi_{k}\big(L_{k_{0},k}(A_{1}^{(1)})\big)$$
for every $k_{0}\in \{0,1,\ldots,k\}$. \hfill $\Box$ 
\end{theorem}}

One is naturally led to assume that a deeper reason may lie behind
the bijections $\phi_{k}$, which in fact suggest a
new kind of duality between appropriately defined module categories
for $A_{1}^{(1)}$ and $A_{2}^{(2)}$. From the viewpoint of representation 
theory, natural approaches for finding conceptual explanations to 
Theorem 2.2.1 would be the following:
\begin{problem}
  Are $\mathcal{O}_{1}(k)$ and $\mathcal{O}_{2}(k)$ the same
  set of inequivalent simple graded objects in a semisimple category
  induced by the ``equal level'' representations of some algebraic structure?
\end{problem}
\begin{problem}
  Does there exist an algebraic structure for which the
  spaces $L_{k_{0},k}(A_{1}^{(1)})$ and 
  $\Omega_{k_{0},2k+1}$ are isomorphic simple graded modules for
  all $k_{0}\in 
  \{0,1,\ldots,k\}$?
\end{problem}
The ``equal level'' representations in the formulation of 
Problem 1 are of course to be understood in a broad sense, by
analogy with the (highest weight) modules of fixed level for an affine
Lie algebra, the (unitary) $Vir$-modules of fixed central charge, the modules
for a (rational) VOA with prescribed rank, etc. Problem 2 is 
obviously a weaker version of Problem 1 and we shall therefore 
concentrate mainly on the latter. We conjecture that 
$\mathcal{O}_{1}(k)$ and $\mathcal{O}_{2}(k)$ are in fact the same set of 
 simple twisted modules for a (rational) VOA or a 
related structure, and we produce some evidence in this direction in Sections 
3 and 4.
\begin{remark}
From a Lie-algebraic point of view, the natural candidates 
for the algebraic structures of Problems 1 and 2 would be
some (appropriately infinite-dimensional) Lie algebras. For instance,
one could try to construct suitable $A_{1}^{(1)}$-module structures
for the spaces in $\mathcal{O}_{2}(k)$. Although there are several
ways of embedding $A_{1}^{(1)}$ into $A_{2}^{(2)}$, these
procedures are not likely to serve our purposes mainly because they
fail to make the spaces $\Omega_{k_{0},2k+1}$ become 
$A_{1}^{(1)}$-modules (see also Remark 3.2.6 for a somewhat similar
point of view). Let us consider as an example the following elements
of $A_{2}^{(2)}$: $e_{0}'=\frac{1}{6}\big[e_{0},[[[e_{0},e_{1}],e_{1}],e_{1}]
\big]$, $f_{0}'=\frac{1}{6}\big[f_{0},[[[f_{0},f_{1}],f_{1}],f_{1}]\big]$, 
$e_{1}'=e_{1}$, $f_{1}'=f_{1}$, $h_{0}'=3h_{1}+8h_{0}$, $h_{1}'=h_{1}$, 
where $e_{i}$, $f_{i}$, $h_{i}$, $i=0,1$, are as in \eqref{2.1.4}. Then
$\{e_{i}', f_{i}', h_{i}'\mid i\in \{0,1\} \}$ generates a subalgebra
$\mathfrak{a}$ of $A_{2}^{(2)}$ which is isomorphic to the affine Lie
algebra $A_{1}^{(1)}$ with canonical central element
$c'=h_{0}'+h_{1}'=4c$ (where $c$ is the canonical central element of
$A_{2}^{(2)}$), so that one may identify $\mathfrak{a}$ with the full 
subalgebra of $A_{1}^{(1)}$ of depth 4. However, the action of $\mathfrak{a}$
 on $L\big(\Lambda;\widehat{\mathfrak{sl}(3,\mathbf{C})}[\sigma]\big)$ --
$\Lambda\in P_{+}$ being as in \eqref{2.2.4} -- does not centralize the 
action of the Heisenberg
algebra $\tilde{\mathfrak{s}}[\sigma]'$ defined in \eqref{2.2.8}, and 
so the corresponding vacuum space
$\Omega_{L(\Lambda;\widehat{\mathfrak{sl}(3,\mathbf{C})}[\sigma])}$
is not necessarily $\mathfrak{a}$-invariant. Moreover, the principal
realization of $A_{2}^{(2)}$ induces on $\mathfrak{a}\cong
A_{1}^{(1)}$ the gradation of type $(5,1)$ -- and not $(1,2)$ as
needed -- and the level of
$L\big(\Lambda;\widehat{\mathfrak{sl}(3,\mathbf{C})}[\sigma]\big)$
viewed as
an $\mathfrak{a}$-module increases to $\Lambda(c')=4(2k+1)$. 
\end{remark}

\section{Standard modules for $A_{1}^{(1)}$ and $A_{2}^{(2)}$ as 
  representations of VOAs}

In this section we give a first approach to Problem 1. We start by noticing
 that after some appropriate modifications, expressions
\eqref{2.2.2}-\eqref{2.2.3} and \eqref{2.2.5}-\eqref{2.2.6} become the
$q$-traces of certain (twisted) modules 
for VOAs associated to vacuum representations
of $A_{1}^{(1)}$ and $A_{2}^{(1)}$ respectively. As it was previously 
indicated, the ultimate goal would be to interpret both
$\mathcal{O}_{1}(k)$ and $\mathcal{O}_{2}(k)$ 
as the same set of equivalence classes of
irreducible (twisted) modules for some rational VOA or a related structure.  
We study various possibilities for such a structure.

\subsection{GVOAs and modules}

We review here some necessary background on GVOAs and their modules. For the 
formal calculus involved and further results we refer to \cite{B}, 
\cite{FLM}, \cite{FHL}, \cite{DL1}, \cite{Li1}, \cite{Li2}, \cite{Z}.

{\allowdisplaybreaks
\begin{definition}[\cite{DL1}]
Let $S\in \mathbf{Z}_{+}$ and let $G$ be a finite abelian group
endowed with a symmetric nondegenerate
$\tfrac{1}{S}\mathbf{Z}/2\mathbf{Z}$-valued $\mathbf{Z}$-bilinear
form:
\begin{equation*}
  (g,h)\in \tfrac{1}{S}\mathbf{Z}/2\mathbf{Z} \text{ for } g,h\in G.
\end{equation*}
A {\em generalized vertex operator algebra of level $S$} associated with 
the group $G$ and the form $(\cdot\,,\cdot)$ is a vector space $V$ with 
two gradations
\begin{equation*}
  V=\coprod_{n\in \frac{1}{S}\mathbf{Z}}V_{n}=\coprod_{g\in G}V^{g} 
\text{ with } \mbox{wt}(v)=n \text{ for } v\in V_{n},
\end{equation*}
such that
{\allowdisplaybreaks
  \begin{align*}
    & V^{g}=\coprod_{n\in \frac{1}{S}\mathbf{Z}}V_{n}^{g}, \text{
      where } V_{n}^{g}=V_{n}\cap V^{g} \text{ for any } g\in G \text{
      and } n\in \tfrac{1}{S}\mathbf{Z},\\
    & \dim V_{n}<\infty \text{ for all } n\in \tfrac{1}{S}\mathbf{Z},
\,V_{n}=0 \text{ for } n \ll 0,\\
\intertext{which is equipped with a linear map}
& Y(\cdot,z):\, V\longrightarrow (\mbox{End}\, V)\big[\big[z^{1/S},z^{-1/S}
\big]\big] \\
& \phantom{Y(\cdot,z):\,\,\,\, }v\longmapsto Y(v,z)=\sum_{n\in \frac{1}{S}
\mathbf{Z}}v_{n}z^{-n-1}
\intertext{and with two distinguished vectors $\mathbf{1}\in V_{0}^{0}$,
$\omega\in V_{2}^{0}$, satisfying the following conditions: for any $g,h\in
G$, $u,v\in V$ and $m\in \tfrac{1}{S}\mathbf{Z}$}
  & u_{m}V^{h}\subset V^{g+h} \text{ if } u\in V^{g},\\
  & u_{m}v=0 \text{ if } m \gg 0,\\
  & Y(\mathbf{1},z)=\mbox{id}_{V},\\
  & Y(v,z)\mathbf{1}\in V[[z]] \text{ and } \!\lim_{z\rightarrow
0}Y(v,z)\mathbf{1}=v,\\
& Y(v,z)\big|_{V^{h}}=\sum_{n\equiv (g,h)\bmod \mathbf{Z}}v_{n}z^{-n-1}
\text{ if } v\in V^{g}\\
\intertext{(i.e., $n+2\mathbf{Z}\equiv (g,h)\bmod
  \mathbf{Z}/2\mathbf{Z}$); if $u\in V^{g}$ and $v\in V^{h}$ then the
  following generalized Jacobi identity holds:}
& z_{0}^{-1}\!\left(\!\frac{z_{1}-z_{2}}{z_{0}}\!\right)^{(g,h)}\!\delta\!
\left(\!\frac{z_{1}-z_{2}}{z_{0}}\!\right)\!Y(u,z_{1})Y(v,z_{2})\\
   & \!\!-z_{0}^{-1}\!\left(\!\frac{z_{2}-z_{1}}{z_{0}}\!\right)^{(g,h)}\!
\delta\!
\left(\!\frac{z_{2}-z_{1}}{-z_{0}}\!\right)\!Y(v,z_{2})Y(u,z_{1})
 \\
& \!\!=z_{2}^{-1}\delta\!\left(\!\frac{z_{1}-z_{0}}{z_{2}}\!\right)\!
Y(Y(u,z_{0})v,z_{2})\!\left(\!\frac{z_{1}-z_{0}}{z_{2}}\!\right)^{-g},
 \\
\intertext{where $\delta (z)=\sum_{n\in \mathbf{Z}}z^{n}$ and whenever 
$k\in G$ and $w\in V^{k}$,}
& \left(\!\frac{z_{1}-z_{0}}{z_{2}}\!\right)^{-g}\!\delta\!\left(\!
\frac{z_{1}-z_{0}}{z_{2}}\!\right)\cdot
w=\left(\!\frac{z_{1}-z_{0}}{z_{2}}\!\right)^{-(g,k)}\!\delta\!\left(\!
\frac{z_{1}-z_{0}}{z_{2}}\!\right)w;\\
\intertext{furthermore, for any $m,n\in \mathbf{Z}$ one has}
& [L(m),L(n)]=(m-n)L(m+n)+\frac{1}{12}(m^{3}-m)\delta_{m+n,0}(\mbox{rank}\,
V)\mbox{id}_{V},\\
\intertext{where $L(n)=\omega_{n+1}$ for $n\in \mathbf{Z}$, i.e., 
$Y(\omega,z)=\sum_{n\in \mathbf{Z}}L(n)z^{-n-2}$, $\mbox{rank}\,V\in 
\mathbf{C}$, and}
&  L(0)v=nv=\mbox{wt}(v)v \text{ for }n\in \tfrac{1}{S}\mathbf{Z},\, v\in
V_{n},\\
& Y(L(-1)v,z)=\frac{d}{dz}Y(v,z).
\end{align*}
\noindent
This completes the definition of a generalized vertex operator algebra 
(GVOA).} The GVOA defined above is denoted by
$\big(V,Y,\mathbf{1},\omega,S,G,(\cdot\,,\cdot)\big)$ 
or simply by $V$ if no confusion is possible.
\end{definition}}

A $\tfrac{1}{S}\mathbf{Z}$-{\em graded VOA} is a GVOA of level $S$ associated 
with the group $G=\{0\}$. Following \cite{DLiM1}, we call a GVOA of 
this type a $\mathbf{Q}$-{\em graded VOA}. A VOA is then a 
$\mathbf{Q}$-graded VOA of level 1.

{\allowdisplaybreaks
\begin{definition}[\cite{DL1}, \cite{DLiM1}]
Let $\big(V,Y,\mathbf{1},\omega,S\big)$ be a $\mathbf{Q}$-graded VOA. A {\em 
weak $V$-module} is a pair $(M,Y_{M})$, where $M$ is a vector space and 
$Y_{M}(\cdot,z)$ is a linear map $V\rightarrow (\mbox{End}\,M )
\big[\big[z^{1/S},z^{-1/S}\big]\big]$ satisfying the following: 
$Y_{M}(\mathbf{1},z)=\mbox{id}_{M}$, $z^{n}Y_{M}(a,z)u\in M
\big[\big[z^{1/S}
\big]\big]$ for any $a\in V$, $u\in M$, and $n\in \tfrac{1}{S}\mathbf{Z}$ 
sufficiently large, $Y_{M}(L(-1)a,z)=\tfrac{d}{dz}Y_{M}(a,z)$ for $a\in V$, 
and the Jacobi identity
{\allowdisplaybreaks
 \begin{align*}
  & z_{0}^{-1}\delta\!\left(\!\frac{z_{1}-z_{2}}{z_{0}}\!\right)\!Y_{M}
(a,z_{1})
Y_{M}(b,z_{2})-z_{0}^{-1}\delta\!\left(\!\frac{z_{2}-z_{1}}{-z_{0}}\!\right)
\!Y_{M}(b,z_{2})Y_{M}(a,z_{1}) \\
  & \!\!=z_{2}^{-1}\delta\!\left(\!\frac{z_{1}-z_{0}}{z_{2}}\!\right)\!Y_{M}
(Y(a,z_{0})b,z_{2}),
 \end{align*}
\noindent
for $a, b\in V$.} A weak $V$-module $(M,Y_{M})$ is called a (ordinary) 
$V$-{\em module} if $L(0)$ ($=\mbox{Res}_{z}zY_{M}(\omega,z)$)
 acts semisimply on 
$M$ with the decomposition into $L(0)$-eigenspaces $M=\coprod_{h\in 
\mathbf{C}}M_{h}$ such that for any $h\in \mathbf{C}$, $\dim M_{h}< \infty$ 
and $M_{h+n}=0$ for $n\in \tfrac{1}{S}\mathbf{Z}$ sufficiently small. 
A $\mathbf{Q}$-{\em graded weak module} is a weak $V$-module $(M,Y_{M})$ 
together with a $\mathbf{Q}$-gradation $M=\coprod_{n\in \mathbf{Q}}M_{n}$ 
such that $a_{m}M_{n}\subseteq M_{r+n-m-1}$ for $a\in V_{r}$, $r$, $m\in 
\tfrac{1}{S}\mathbf{Z}$, $n\in \mathbf{Q}$. 
An {\em admissible $V$-module} is a $\mathbf{Q}_{\geq 0}$-graded weak 
$V$-module. A $\mathbf{Q}$-graded VOA $V$ is said to be {\em rational} if 
every admissible $V$-module is completely reducible, i.e., a direct sum of 
simple admissible $V$-modules.
\end{definition}}

Let $\big(V,Y,\mathbf{1},\omega\big)$ be a VOA with an automorphism
 $\sigma$ of order $T$. Set $V^{k}=\{a\in V\mid
  \sigma (a)=\exp (2k\pi i/T)a \}$, $0\le k\le T-1$, so that 
$V=\oplus_{k=0}^{T-1}V^{k}$.

{\allowdisplaybreaks
\begin{definition}[\cite{BDM}, \cite{DLiM2}]
A {\em weak} $\sigma$-{\em twisted} $V$-{\em module} is a pair  
$(M,Y_{M})$, where $M$ is a vector space and $Y_{M}(\cdot,z)$ is a 
linear map 
$V\rightarrow (\mbox{End}\,M )\big[\big[z^{1/T},z^{-1/T}\big]\big]$ given by 
$a\mapsto Y_{M}(a,z)=\sum_{n\in \mathbf{Q}}a_{n}z^{-n-1}$, such that for 
$a$, $b\in V$ and $u\in M$ the following hold: $a_{n}u=0$ if $n\gg 0$, 
$Y_{M}(\mathbf{1},z)=\mbox{id}_{M}$, $Y_{M}(a,z)=\sum_{n\in k/T+
\mathbf{Z}}
a_{n}z^{-n-1}$ for $a\in V^{k}$, and the twisted Jacobi identity
{\allowdisplaybreaks
 \begin{align*}
  & z_{0}^{-1}\delta\!\left(\!\frac{z_{1}-z_{2}}{z_{0}}\!\right)\!Y_{M}
(a,z_{1})
Y_{M}(b,z_{2})-z_{0}^{-1}\delta\!\left(\!\frac{z_{2}-z_{1}}{-z_{0}}\!\right)
\!Y_{M}(b,z_{2})Y_{M}(a,z_{1}) \\
  & \!\!=z_{2}^{-1}\left(\!\frac{z_{1}-z_{0}}{z_{2}}\!\right)^{-\frac{k}{T}}
\delta\!\left(\!\frac{z_{1}-z_{0}}{z_{2}}\!\right)\!Y_{M}
(Y(a,z_{0})b,z_{2}),
 \end{align*}
\noindent
for $a\in V^{k}$, $b\in V$.} A weak $\sigma$-twisted $V$-module $(M,Y_{M})$ 
is called a (ordinary) $\sigma$-{\em twisted} $V$-{\em module} if $L(0)$ acts 
semisimply on 
$M$ with the decomposition into $L(0)$-eigenspaces $M=\coprod_{h\in 
\mathbf{C}}M_{h}$ such that for any $h\in \mathbf{C}$, $\dim M_{h}< \infty$ 
and $M_{h+n}=0$ for $n\in \tfrac{1}{T}\mathbf{Z}$ sufficiently small. 
A $\tfrac{1}{T}\mathbf{Z}$-{\em graded weak} $\sigma$-{\em twisted} $V$-{\em 
module} is a weak $\sigma$-twisted $V$-module $(M,Y_{M})$ which carries a 
$\tfrac{1}{T}\mathbf{Z}$-grading $M=\coprod_{n\in \tfrac{1}{T}\mathbf{Z}}
M_{n}$ such that $a_{m}M_{n}\subseteq M_{r+n-m-1}$ for $a\in V_{r}$, $r\in 
\mathbf{Z}$, $m$, $n\in \tfrac{1}{T}\mathbf{Z}$. 
An {\em admissible} $\sigma$-{\em twisted} $V$-{\em module} is a 
$\tfrac{1}{T}\mathbf{N}$-graded weak $\sigma$-twisted $V$-module. 
A  VOA $V$ is said to be $\sigma$-{\em rational} if every admissible 
$\sigma$-twisted $V$-module is completely reducible.
\end{definition}}

It was shown in \cite{DLiM2} that if $(M,Y_{M})$ is a 
weak $\sigma$-twisted $V$-module then the component operators of 
$Y_{M}(\omega,z)$ together with $\mbox{id}_{M}$ induce a 
$Vir$-representation on $M$ with central charge $\mbox{rank\,V}$ and  
$Y_{M}(L(-1)a,z)=\tfrac{d}{dz}Y_{M}(a,z)$ for $a\in V$. It was further 
shown in {\em loc.~cit.} that a $\sigma$-rational VOA has only 
finitely many isomorphism classes of simple admissible $\sigma$-twisted 
modules and that every such module is an ordinary $\sigma$-twisted 
module.

The following are consequences of the above definitions (cf.~\cite{Li2}): if 
$V$ is a VOA, $\sigma \in \text{Aut}(V)$ has order $T$, $0\le k\le T-1$, 
$a\in V^{k}$, $b\in V$, and $(M,Y_{M})$ is a weak $\sigma$-twisted 
$V$-module then
{\allowdisplaybreaks
  \begin{align}
    & [Y_{M}(a,z_{1}),Y_{M}(b,z_{2})]=\sum_{j=0}^{\infty}\frac{1}{j!}\!
\left(\!\left(\!\frac{\partial}{\partial z_{2}}\!\right)^{j}z_{1}^{-1}
\delta\!\left(\!\frac{z_{2}}{z_{1}}\!\right)\left(\!\frac{z_{2}}{z_{1}}\!
\right)^{\frac{k}{T}}\!\right)\!Y_{M}(a_{j}b,z_{2}),\label{3.1.1}  \\
    & Y_{M}(Y(a,z_{0})b,z_{2})=\mbox{Res}_{z_{1}}\left(\!\frac{z_{1}-z_{0}}
{z_{2}}\!\right)^{\frac{k}{T}}\!\Bigg[z_{0}^{-1}\delta\!\left(\!\frac{
z_{1}-z_{2}}{z_{0}}\!\right)\!Y_{M}(a,z_{1})Y_{M}(b,z_{2})\label{3.1.2}\\
    & \phantom{Y_{M}(Y(a,z_{0})b,z_{2})=\mbox{Res}_{z_{1}}\left(\!
\frac{z_{1}-z_{0}}{z_{2}}\!\right)^{\frac{k}{T}}}\quad 
    -z_{0}^{-1}\delta\!\left(\!\frac{-z_{2}+z_{1}}{z_{0}}\!\right)\!
Y_{M}(b,z_{2})Y_{M}(a,z_{1})\Bigg]. \nonumber 
\end{align}
These formulas reduce to the ordinary commutator 
respectively associator relations for VOAs by taking $M=V$ and 
$\sigma=\mbox{id}_{V}$.} As in the untwisted case (cf.~\cite[Lemma 1.2.1]{Z}) 
one can prove:

{\allowdisplaybreaks
\begin{prop+}
  Let $V$ be a VOA with an automorphism $\sigma$ of order $T$. If $M$
  is a simple admissible  
$\sigma$-twisted $V$-module such that $M=\coprod_{n\in
  \frac{1}{T}\mathbf{N}}M_{n}$ and $M_{0}\neq 0$, 
$\dim M_{0}<\infty$, then there
exists $h_{M}\in \mathbf{C}$ such that
$L(0)\big|_{M_{n}}=(n+h_{M})\mbox{{\em id}}_{_{M_{n}}}$ for all $n\in
  \frac{1}{T}\mathbf{N}$.\hfill $\Box$
\end{prop+}

By Proposition 3.1.1, one may write the $q$-trace (conformal character) of an
irreducible $\sigma$-twisted 
$V$-module $M=\coprod_{n\in \frac{1}{T}\mathbf{N}}M_{n}$ as the
following formal power series
\begin{equation}\label{3.1.3}
  \mbox{tr}_{_{M}} q^{L(0)-\frac{rank\, V}{24}}=q^{h_{M}-\frac{rank\,
      V}{24}}\sum_{n=0}^{\infty}(\dim M_{\frac{n}{T}})q^{\frac{n}{T}}.
\end{equation}
Obviously, any admissible $\sigma$-twisted 
$V$-module with finite-dimensional $L(0)$-homogeneous subspaces has 
a well-defined $q$-trace.
\begin{prop+}[\cite{DLiM2}, \cite{Lia}]
  Let $(V,Y,\mathbf{1},\omega)$ be a VOA of rank $r\in \mathbf{C}$ and
  suppose that $h\in V$ is such that $L(n)h=\delta_{n,0}h,\, h_{n}h=
  \delta_{n,1}\lambda\mathbf{1}\text{ for } n\in \mathbf{N}$,
  where $\lambda\in \mathbf{C}$. Then
  $(V,Y,\mathbf{1},\omega+h_{-2}\mathbf{1})$ is a vertex algebra of
  rank $r-12\lambda$.
\end{prop+}
\begin{rema+}
  The vector $\omega+h_{-2}\mathbf{1}$ becomes a new Virasoro element
  of $V$, on which it induces a $\mathbf{C}$-gradation instead of the
  original $\mathbf{Z}$-gradation without altering the other axioms in
  the definition of a VOA. We shall use Proposition 3.1.2 only when the  
  new Virasoro vector induces a $\tfrac{1}{S}\mathbf{Z}$-gradation 
  on $V$ for some $S\in \mathbf{Z}_{+}$ so that the new structure becomes
  a $\mathbf{Q}$-graded VOA (cf.~Definition 1).
\end{rema+}  

We now describe briefly the VOAs associated with affine Lie algebras, which 
are sometimes called affine VOAs (cf.~\cite{K2}). Let $\mathfrak{g}$ be a 
finite-dimensional Lie algebra with the form $\langle \cdot,\!\cdot
\rangle$ normalized as in \S 2.1, and form the associated untwisted
affine Kac-Moody algebra $\tilde{\mathfrak{g}}=\hat{\mathfrak{g}}\rtimes
  \mathbf{C}d$ as in \eqref{2.1.2}-\eqref{2.1.3}. Set 
  $\tilde{\mathfrak{g}}_{\ge 0}=\oplus_{n\ge 0}\mathfrak{g}(n)\oplus
  \mathbf{C}c\oplus \mathbf{C}d$ and let $-h\check{}\neq l\in \mathbf{C}$. 
Define a $\tilde{\mathfrak{g}}_{\ge 0}$-module structure on $\mathbf{C}$ by 
$c\cdot 1=l\cdot 1$, $d\cdot 1=0$, $\mathfrak{g}(n)\cdot 1=0$ for $n\ge 0$, 
and form the induced $\tilde{\mathfrak{g}}$-module $N(l\Lambda_{0})=
U(\tilde{\mathfrak{g}})\otimes_{U(\tilde{\mathfrak{g}}_{\ge 0})}\mathbf{C}$ 
(the so-called vacuum representation of $\tilde{\mathfrak{g}}$ at level $l$). 
Note that $N(l\Lambda_{0})$ is a 
restricted $\tilde{\mathfrak{g}}$-module such that $N(l\Lambda_{0})\cong
U\big(\!\oplus_{n<0}\mathfrak{g}(n)\big)$ as vector spaces, and that we may
 identify $\mathfrak{g}(-1)\otimes 1$ with $\mathfrak{g}$. Set 
  $\mathbf{1}=1\otimes 1\in N(l\Lambda_{0})$ and define the element 
\begin{equation*}
  \omega=\frac{1}{2(l+h\check{})}\sum_{j=1}^{\dim
    \mathfrak{g}}a^{j}(-1)b^{j}(-1)\mathbf{1}\in N(l\Lambda_{0}),
\end{equation*}
where $\{a^{j}\mid j\in \{1,\ldots,\dim \mathfrak{g}\}\}$ and
$\{b^{j}\mid j\in \{1,\ldots,\dim \mathfrak{g}\}\}$ 
are dual bases of $\mathfrak{g}$ with respect to $\langle \cdot,\! \cdot
\rangle$. Recall from \eqref{2.1.5} the series $a(z)$ in
this case and define the map 
\begin{gather}
  Y:\, \mathfrak{g}(-1)\mathbf{1} \longrightarrow (\mbox{End}\,
  N(l\Lambda_{0}))[[z,z^{-1}]]\label{3.1.4}\\
  Y(a(-1)\mathbf{1},z)=a(z),\,a\in \mathfrak{g}.\notag
\end{gather}
One can show (cf., e.g., \cite[Theorem 2.6]{MP2}) that $Y$ extends uniquely 
to $N(l\Lambda_{0})$ in such a way that $N(l\Lambda_{0})$ becomes a VOA with
vacuum vector $\mathbf{1}$ and Virasoro 
element $\omega$ such that $\mathfrak{g}(-1)\mathbf{1}=N(l\Lambda_{0})_{1}$ 
(the weight one subspace of $N(l\Lambda_{0})$). 
Moreover, given any restricted $\hat{\mathfrak{g}}$-module $M$ of level
$l$, there is a canonical extension to $N(l\Lambda_{0})$ of the map
\begin{gather*}
  Y_{M}:\, \mathfrak{g}(-1)\mathbf{1} \longrightarrow (\mbox{End}\,
  M)[[z,z^{-1}]]\\
  Y_{M}(a(-1)\mathbf{1},z)=a(z), \text{ } a\in \mathfrak{g},\notag
\end{gather*}
such that $(M,Y_{M})$ becomes a weak $N(l\Lambda_{0})$-module
(cf. ~\cite{DL1}, \cite{FZ}, \cite{Li1}, \cite{Lia}, \cite{MP2}). Let finally 
$N^{1}(l\Lambda_{0})$ be
the (unique) maximal proper $\tilde{\mathfrak{g}}$-submodule of
$N(l\Lambda_{0})$ and notice that we may identify the irreducible
quotient $N(l\Lambda_{0})/N^{1}(l\Lambda_{0})$ with the standard 
$\tilde{\mathfrak{g}}$-module $L(l\Lambda_{0})$. We summarize
some of the results of the above-mentioned papers in
{\allowdisplaybreaks
\begin{theo+}
  For each $l\neq
  -h\check{}$, $(N(l\Lambda_{0}),Y,\mathbf{1},\omega)$ is a VOA of rank
  $\tfrac{l\dim \mathfrak{g}}{l+h\check{}}$ generated by 
$\mathfrak{g}(-1)\mathbf{1}$ and any restricted 
 $\hat{\mathfrak{g}}$-module of level
  $l$ is a weak $N(l\Lambda_{0})$-module. Every 
$\tilde{\mathfrak{g}}$-submodule of
$N(l\Lambda_{0})$ is an ideal of $N(l\Lambda_{0})$ viewed as a VOA. 
In particular, there
exists an induced structure of simple VOA on
$L(l\Lambda_{0})$. If $l$ is a positive integer, then
$L(l\Lambda_{0})$ is rational and the set of equivalence classes of
simple $L(l\Lambda_{0})$-modules is exactly the set of equivalence
classes of standard $\hat{\mathfrak{g}}$-modules of level $l$.
\end{theo+}  

Any automorphism $\sigma$ of order
$T$ of $\mathfrak{g}$ preserves the form $\langle \cdot,\!\cdot \rangle$ 
and induces a Lie algebra automorphism of $\tilde{\mathfrak{g}}$. Since 
the VOA $N(l\Lambda_{0})$ is generated by $N(l\Lambda_{0})_{1}=
\mathfrak{g}(-1)\mathbf{1}$, it 
follows from the associator formula for VOAs that $\sigma$ 
 induces a VOA automorphism of $N(l\Lambda_{0})$ 
(cf.~\cite[Proposition 6.20]{Lia}). It is easy to see that 
$N^{1}(l\Lambda_{0})$ is invariant under the induced VOA automorphism. 
Hence $\sigma$ induces a VOA 
automorphism of $L(l\Lambda_{0})$ as well. We denote
these induced VOA automorphisms again by $\sigma$. If $M$ is a restricted
$\hat{\mathfrak{g}}[\sigma]$-module of level $l$, the map
\begin{gather}
  Y^{\sigma}_{M}:\, \mathfrak{g}(-1)\mathbf{1} \longrightarrow (\mbox{End}\,
  M)\big[\big[z^{1/T},z^{-1/T}\big]\big]\label{3.1.5}\\
  Y^{\sigma}_{M}(a(-1)\mathbf{1},z)=a(\sigma;z),\notag
\end{gather}
for $a\in \mathfrak{g}_{(j)}$, $j=0,\ldots,T-1$, has a unique extension
to $N(l\Lambda_{0})$ that makes $(M,Y^{\sigma}_{M})$ a weak
$\sigma$-twisted $N(l\Lambda_{0})$-module. This is a consequence of
the theory of local systems of twisted vertex operators developed in
\cite{Li2}, where the following result was obtained.}
{\allowdisplaybreaks
\begin{theo+}
  Let $l\neq -h\check{}$ be a complex number. Then any restricted
  $\hat{\mathfrak{g}}[\sigma]$-module of level $l$ is a weak
  $\sigma$-twisted $N(l\Lambda_{0})$-module. If $l$ is a positive
  integer, then $L(l\Lambda_{0})$ is $\sigma$-rational and the set of
  equivalence classes of simple $\sigma$-twisted
  $L(l\Lambda_{0})$-modules is precisely the set of equivalence
  classes of standard $\hat{\mathfrak{g}}[\sigma]$-modules of level $l$.
\end{theo+}

Given $a\in \mathfrak{g}_{(j)}$, $b\in \mathfrak{g}$, and a restricted
$\hat{\mathfrak{g}}[\sigma]$-module $M$ of level $l$, one gets from the
affine Lie algebra relations \eqref{2.1.3} that}
\begin{equation}\label{3.1.6}
  [a(\sigma;z_{1}),b(\sigma;z_{2})]=z_{1}^{-1}\!\left(\!
\frac{z_{2}}{z_{1}}\!\right)^{\frac{j}{T}}\!\delta\!\left(\!
\frac{z_{2}}{z_{1}}\!\right)\![a,b](\sigma;z_{2})+\langle a,\! b \rangle 
\frac{\partial}{\partial z_{2}}\!\left(\!z_{1}^{-1}\delta\!\left(\!
\frac{z_{2}}{z_{1}}\!\right)\left(\!\frac{z_{2}}{z_{1}}\!\!\right)^
{\frac{j}{T}}\right)\!l\mbox{id}_{M}
\end{equation}
in $(\mbox{End}\,
M)\big[\big[z_{1}^{1/T},z_{1}^{-1/T},z_{2}^{1/T},z_{2}^{-1/T}\big]\big]$. 
Following \cite{DL1}, we define a (noncommutative)
``normal ordering'' operation for $a\in \mathfrak{g}_{(j)}$, $b\in
\mathfrak{g}_{(k)}$, $m\in \frac{j}{T}+\mathbf{Z}$, $n\in
\frac{k}{T}+\mathbf{Z}$, by
\begin{equation*}
  _{\text{{\tiny $\times$}}}^{\text{{\tiny $\times$}}}a_{m}b_{n}
\phantom{}_{\text{{\tiny $\times$}}}^{\text{{\tiny $\times$}}}=
  \begin{cases}
    a_{m}b_{n}, &\text{if $m<0$,} \\
    b_{n}a_{m}, &\text{if $m\ge 0$.}
  \end{cases}  
\end{equation*}
\noindent
Using the computational techniques of \cite{Li2} one can then easily 
prove:} 
{\allowdisplaybreaks 
\begin{lemm+}
  Let $M$ be a restricted $\hat{\mathfrak{g}}[\sigma]$-module of level $l$, 
and let $a\in \mathfrak{g}_{(k)}$, $b\in 
\mathfrak{g}_{(-k)}$ for some fixed $k\in \{0,\ldots,T-1\}$. Then
\begin{align}
  Y_{M}^{\sigma}(a(-1)b(-1)\mathbf{1},z)=\sum_{s\in
    \mathbf{Z}}\Bigg(&\sum_{n\in \mathbf{Z}}\phantom{}_{\text{{\tiny
        $\times$}}}^{\text{{\tiny
        $\times$}}}a\big(n+k/T\big)b\big(s-(n+k/T)\big)\phantom{}_
{\text{{\tiny $\times$}}}^{\text{{\tiny $\times$}}}\label{3.1.7}\\
  &-\tfrac{k}{T}[a,b](s)-\langle a,\!b \rangle
  \delta_{s,0}\binom{\frac{k}{T}}{2} l\mbox{{\em id}}_{M}\!\Bigg)
z^{-s-2}.\nonumber
\end{align}
In particular,
\begin{align}
  \mbox{{\em Res}}_{z}\,zY_{M}^{\sigma}(a(-1)b(-1)\mathbf{1},z)=&\sum_{n\in 
\mathbf{Z}}\phantom{}_{\text{{\tiny
        $\times$}}}^{\text{{\tiny
        $\times$}}}a\big(n+k/T\big)b\big(\!-(n+k/T)\big)\phantom{}_
{\text{{\tiny $\times$}}}^{\text{{\tiny $\times$}}}\label{3.1.8} \\
  &-\tfrac{k}{T}[a,b](0)-\langle a,\!b \rangle
  \binom{\frac{k}{T}}{2} l\mbox{{\em id}}_{M}.\tag*{$\Box$}
\end{align}
\end{lemm+}
\indent
Recall finally that} given two VOAs 
$(V_{i},Y_{_{V_{i}}},\mathbf{1}_{_{V_{i}}},\omega_{_{V_{i}}})$,
$i=1,2$, one may define the tensor product VOA $(V_{1}\otimes
V_{2},Y_{_{\otimes }},\mathbf{1}_{_{\otimes }},\omega_{_{\otimes }})$ by 
setting
{\allowdisplaybreaks
  \begin{equation}\label{3.1.9}
    \begin{split}
      & \mathbf{1}_{_{\otimes }}=\mathbf{1}_{_{V_{1}}}\otimes 
\mathbf{1}_{_{V_{2}}},\,
         \omega_{_{\otimes }}=\omega_{_{V_{1}}}\otimes
         \mathbf{1}_{_{V_{2}}}+\mathbf{1}_{_{V_{1}}}\otimes
         \omega_{_{V_{2}}},\\
      & Y_{_{\otimes }}(a_{1}\otimes a_{2},z)=Y_{_{V_{1}}}(a_{1},z)\otimes
      Y_{_{V_{2}}}(a_{2},z) \text{ for } a_{i}\in V_{i}.
    \end{split}
  \end{equation}
The new VOA is endowed with the tensor product grading and the
central charge for the Virasoro algebra relations 
becomes $\text{rank}\,(V_{1}\otimes
V_{2})=\text{rank}\,V_{1}+\text{rank}\,V_{2}$.} 

\subsection{Structures associated to standard modules for
  $A_{1}^{(1)}$, $A_{2}^{(1)}$, and $A_{2}^{(2)}$}

We start our discussion of Problem 1 with some remarks on the
objects in $\mathcal{O}_{1}(k)$, namely the level $k$
standard $A_{1}^{(1)}$-modules. By Theorem 3.1.4, $V_{k}(A_{1}^{(1)}):=
L(k\Lambda_{0};A_{1}^{(1)})$ is a
rational VOA of rank $c_{1}(k):=\frac{3k}{k+2}$ and the level $k$
standard $A_{1}^{(1)}$-modules are precisely its simple 
modules. These modules have $q$-traces which are essentially the same
as their homogeneously specialized characters when viewed as
$A_{1}^{(1)}$-modules, but we need a setting that makes the former
agree with the $(1,2)$-specialization of the latter. There are two
natural and compatible ways to achieve this, as we shall now
explain.

Recall from \S 2.1 the standard basis $\{e,f,h\}$ of
$\mathfrak{sl}(2,\mathbf{C})$ and the form $\langle
\cdot,\!\cdot \rangle$, and note that $\langle
h,\!h \rangle=2\langle
e,\!f \rangle=2$, $\langle
h,\!e \rangle=\langle
h,\!f \rangle=0$. A Virasoro element of $V_{k}(A_{1}^{(1)})$ is then 
given by 
\begin{equation}\label{3.2.1}
\omega=\frac{1}{2(k+2)}\Big[e(-1)f(-1)\mathbf{1}+f(-1)e(-1)\mathbf{1}+\frac{1}{2}h(-1)^{2}\mathbf{1}\Big]\in
V_{k}(A_{1}^{(1)})_{2}.
\end{equation}
Let $\tilde{\omega}=\omega+\frac{1}{3}h(-2)\mathbf{1}$ and 
$Y(\tilde\omega,z)=\sum_{n\in \mathbf{Z}}\tilde{L}(n)z^{-n-2}$. Then
\begin{align}
& \tilde{L}(0)=L(0)-\frac{1}{3}h(0),\text{ and}\label{3.2.2}\\
& [\tilde{L}(0),e(n)]=-\bigg(n+\frac{2}{3}\bigg)e(n),\,[\tilde{L}(0),f(n)]=
-\bigg(n-\frac{2}{3}\bigg)f(n),\label{3.2.3}\\
& [\tilde{L}(0),h(n)]=-nh(n)\text{ for }n\in 
\mathbf{Z}.\nonumber
\end{align}
\noindent
By Proposition 3.1.2 and Remark 3.1.3, 
$\big(V_{k}(A_{1}^{(1)}),Y,\mathbf{1},\tilde{\omega}\big)$ is a 
$\tfrac{1}{3}\mathbf{Z}$-graded VOA of rank $\tilde{c}_{1}(k):=
c_{1}(k)-\tfrac{8k}{3}$. Since $k\in \mathbf{Z}_{+}$, it follows from 
\cite[Theorem 3.7]{DLiM5} that the set of inequivalent 
simple weak $\big(V_{k}(A_{1}^{(1)}),Y,\mathbf{1},\tilde{\omega}\big)$-modules 
is exactly $\mathcal{O}_{1}(k)$. Moreover, \eqref{3.2.3} implies that for 
$n\in \mathbf{N}$, $e(n)$, $f(n+1)$, and $h(n+1)$ have negative degrees 
with respect to the operator $\tilde{L}(0)$, so that by 
\cite[Theorem 2.20]{DLiM1} any weak (in particular, any admissible) 
$\big(V_{k}(A_{1}^{(1)}),Y,\mathbf{1},\tilde{\omega}\big)$-module is 
completely reducible. Thus 
$\big(V_{k}(A_{1}^{(1)}),Y,\mathbf{1},\tilde{\omega}\big)$ is a rational 
$\tfrac{1}{3}\mathbf{Z}$-graded VOA of rank $\tilde{c}_{1}(k)$ and  
$\mathcal{O}_{1}(k)$ is the complete set of its inequivalent simple modules 
(cf.~\cite[Theorem 3.14]{DLiM1}). By analogy with \eqref{3.1.3},
the $q$-traces of these modules are the modified graded dimensions induced
by $\tilde{L}(0)$:
\begin{equation}\label{3.2.4}
\chi_{_{k_{0},k}}(q):=\text{tr}_{L_{k_{0},k}(A_{1}^{(1)})}q^{\tilde{L}(0)-
\frac{\tilde{c}_{1}(k)}{24}}=q^{h_{k_{0},k}-\frac{\tilde{c}_{1}(k)}{24}}
\sum_{n=0}^{\infty}(\dim
L_{k_{0},k}(A_{1}^{(1)})_{_{\frac{n}{3}}})q^{\frac{n}{3}},
\end{equation}
where $h_{k_{0},k}\in \mathbf{C}$ is the lowest weight of
$L_{k_{0},k}(A_{1}^{(1)})$ with respect to $\tilde{L}(0)$. Using
\eqref{3.1.8}, \eqref{3.2.1} and \eqref{3.2.2} one easily gets that 
\begin{equation}\label{3.2.5}
h_{k_{0},k}=\frac{3k_{0}^{2}-2k_{0}k-k^{2}+2k_{0}-2k}{12(k+2)},
\end{equation}
while \eqref{3.2.3} implies that
\begin{equation*}
\sum_{n=0}^{\infty}(\dim
L_{k_{0},k}(A_{1}^{(1)})_{_{\frac{n}{3}}})q^{\frac{n}{3}}=F_{(1,2)}
\big(e^{-k_{0}\Lambda_{0}-k_{1}\Lambda_{1}}\text{ch}\,L_{k_{0},k}
(A_{1}^{(1)})\big)\Big|_{q\rightarrow q^{\frac{1}{3}}},
\end{equation*}
which substituted in \eqref{3.2.4} gives
\begin{equation}\label{3.2.6}
\chi_{_{k_{0},k}}(q)=F_{(1,2)}\big(e^{-k_{0}\Lambda_{0}-k_{1}\Lambda_{1}}
\text{ch}\,L_{k_{0},k}(A_{1}^{(1)})\big)\Big|_{q\rightarrow
  q^{\frac{1}{3}}}\cdot q^{h_{k_{0},k}-\frac{\tilde{c}_{1}(k)}{24}}.
\end{equation}
\noindent
Up to the indicated modifications, the characters \eqref{2.2.2} and 
\eqref{2.2.3} may therefore be interpreted as the $q$-traces of the 
simple modules for the rational $\tfrac{1}{3}\mathbf{Z}$-graded VOA
$\big(V_{k}(A_{1}^{(1)}),Y,\mathbf{1},\tilde{\omega}\big)$.

Alternatively, one may view these modified characters as 
the $q$-traces of
the simple $\nu$-twisted modules over the $\nu$-rational VOA
$\big(V_{k}(A_{1}^{(1)}),Y,\mathbf{1},\omega\big)$, where $\nu$ is the
$(1,2)$-automorphism of $\mathfrak{sl}(2,\mathbf{C})$ given by 
$\nu(h)=h$, $\nu(e)=\exp(4\pi i/3)e$, $\nu(f)=\exp(2\pi i/3)f$. 
Indeed, let us denote the canonical generators of
$\widehat{\mathfrak{sl}(2,\mathbf{C})}[\nu]$ as in Remark
2.1.1 and notice that
$h=\tfrac{1}{3}(h_{1}^{\nu}-2h_{0}^{\nu})$. Obviously, the degree
operator $d^{\nu}$ of $\widehat{\mathfrak{sl}(2,\mathbf{C})}[\nu]$ induces a
$\tfrac{1}{3}\mathbf{Z}$-gradation on
$L_{k_{0},k}\big(\widehat{\mathfrak{sl}(2,\mathbf{C})}[\nu]\big)$,
and one has 
\begin{equation}\label{3.2.7}
  \sum_{n=0}^{\infty}(\dim
  L_{k_{0},k}\big(\widehat{\mathfrak{sl}(2,\mathbf{C})}[\nu]
\big)_{\frac{n}{3}})q^{\frac{n}{3}}=F_{(1,2)}\big(e^{-k_{0}\Lambda_{0}-k_{1}
\Lambda_{1}}\text{ch}\, L_{k_{0},k}(A_{1}^{(1)})\big)\big|_{q\rightarrow 
q^{\frac{1}{3}}}.
\end{equation}
\noindent
On the other hand, $\big\{
\big(L_{k_{0},k}\big(\widehat{\mathfrak{sl}(2,\mathbf{C})}[\nu]\big),
Y^{\nu}\big)\mid k_{0}\in\{0,1,\ldots,k\} \big\}$ is the complete set of 
simple $\nu$-twisted modules for the
$\nu$-rational VOA $\big(V_{k}(A_{1}^{(1)}), Y,
\mathbf{1},\omega\big)$ by Theorem 3.1.5. From the twisted commutator formula
\eqref{3.1.1} one gets that 
\begin{equation}\label{3.2.8}
  [L(0),e(n+2/3)]=-(n+2/3)e(n+2/3),\, [L(0),f(n+1/3)]=-(n+1/3)f(n+1/3),
\end{equation}
where $L(0)=\text{Res}_{z}zY^{\nu}(\omega,z)$. By \eqref{3.2.7} and
\eqref{3.2.8} the $q$-trace of
$L_{k_{0},k}\big(\widehat{\mathfrak{sl}(2,\mathbf{C})}[\nu]\big)$ is 
{\allowdisplaybreaks 
\begin{align}
  \chi_{_{k_{0},k}}^{\nu}(q): & =\text{tr}_{L_{k_{0},k}
(\widehat{\mathfrak{sl}(2,\mathbf{C})}[\nu])}q^{L(0)-\frac{c_{1}(k)}{24}}
\label{3.2.9} \\
  & =F_{(1,2)}\big(e^{-k_{0}\Lambda_{0}-k_{1}\Lambda_{1}}\text{ch}\,L_{k_{0},k}
(A_{1}^{(1)})\big)\Big|_{q\rightarrow 
  q^{\frac{1}{3}}}\cdot q^{h_{k_{0},k}^{\nu}-\frac{c_{1}(k)}{24}},\notag
\end{align}
where $h_{k_{0},k}^{\nu}\in \mathbf{C}$ is the lowest weight of
$L_{k_{0},k}\big(\widehat{\mathfrak{sl}(2,\mathbf{C})}[\nu]\big)$
with respect to the operator $L(0)$ of \eqref{3.2.8}.} The value of
$h_{k_{0},k}^{\nu}$ is easily obtained from \eqref{3.1.8} and \eqref{3.2.1}:
\begin{equation}\label{3.2.10}
  h_{k_{0},k}^{\nu}=\frac{9k_{0}^{2}-6k_{0}k+k^{2}+6k_{0}+2k}{36(k+2)}.
\end{equation}
Then \eqref{3.2.5} and \eqref{3.2.10} yield
\begin{equation}\label{3.2.11}
  h_{k_{0},k}-\frac{\tilde{c}_{1}(k)}{24}=h_{k_{0},k}^{\nu}-
\frac{c_{1}(k)}{24}=\frac{18k_{0}^{2}-12k_{0}k+2k^{2}+12k_{0}-5k}{72(k+2)}.
\end{equation}
It now follows from \eqref{3.2.6}, \eqref{3.2.9} and \eqref{3.2.11} that
\begin{equation}\label{3.2.12}
  \chi_{_{k_{0},k}}(q)=\chi_{_{k_{0},k}}^{\nu}(q).
\end{equation}  

Summarizing, we have
{\allowdisplaybreaks
\begin{theo+}
  (i) $\big(V_{k}(A_{1}^{(1)}),Y,\mathbf{1},\tilde{\omega}\big)$ is a 
simple rational $\frac{1}{3}\mathbf{Z}$-graded VOA of rank 
$\tilde{c}_{1}(k)$ and 
the complete set of its inequivalent simple modules is $\big\{\big(L_{k_{0},k}
(A_{1}^{(1)}),Y\big)\mid k_{0}\in \{0,1,\ldots,k\}\big\}$.\newline
 (ii) $\big(V_{k}(A_{1}^{(1)}),Y,\mathbf{1},\omega\big)$  is a 
 simple $\nu$-rational VOA  of rank $c_{1}(k)$ and the complete set of its 
inequiva\-lent simple 
 $\nu$-twisted modules is $\big\{\big(L_{k_{0},k}\big(
 \widehat{\mathfrak{sl}(2,\mathbf{C})}[\nu]\big),
 Y^{\nu}\big)$, $k_{0}\in \{0,1,\ldots,k\}\big\}$.\newline
(iii) The $q$-trace of $\big(L_{k_{0},k}(A_{1}^{(1)}),
Y\big)$ and the $q$-trace of  
$\big(L_{k_{0},k}\big(\widehat{\mathfrak{sl}(2,\mathbf{C})}[\nu]\big),
Y^{\nu}\big)$ satisfy
\begin{equation*}
  \chi_{_{k_{0},k}}(q)=\chi_{_{k_{0},k}}^{\nu}(q)=F_{(1,2)}\big(e^{-k_{0}
\Lambda_{0}-k_{1}\Lambda_{1}}\text{{\em ch}}\,L_{k_{0},k}(A_{1}^{(1)})\big)
\Big|_{q\rightarrow
  q^{\frac{1}{3}}}\cdot q^{h_{k_{0},k}^{\nu}-\frac{c_{1}(k)}{24}}
\end{equation*}
for every $k_{0}\in \{0,1,\ldots,k\}$.\hfill $\Box$
\end{theo+}
As we mentioned in the introduction of Section 3, an ideal answer to Problem
1 would be that $\mathcal{O}_{1}(k)$ and $\mathcal{O}_{2}(k)$ are
in fact the same set of simple modules for a rational VOA-like
structure such that appropriate modifications of the characters
\eqref{2.2.2}-\eqref{2.2.3} are the $q$-traces of these modules.} If
true, this would almost certainly require some modular invariance
properties from the functions $\tilde{\chi}_{_{k_{0},k}}(\tau):=
\chi_{_{k_{0},k}}(e^{2\pi i\tau})$, where $q=e^{2\pi
  i\tau}$, $\Im (\tau)>0$. Indeed, it was proved in \cite{Z} that the linear 
span of the $q$-traces of all the simple $V$-modules becomes a 
(finite-dimensional) module for the modular group in case $V$ is rational and 
$C_{2}$-finite, i.e., the space $\{a_{-2}b \mid a, b\in V\}$ has finite 
codimension in $V$. This modular invariance property has been generalized in 
a suitable sense to $\nu$-rational $C_{2}$-finite VOAs in \cite{DLiM4}. Since 
most of the known $\nu$-rational VOAs -- including 
$\big(V_{k}(A_{1}^{(1)}),Y,\mathbf{1},\omega\big)$ -- are $C_{2}$-finite 
(cf.~{\em loc.~cit.}), Theorem 3.2.1 implies that the characters 
$\tilde{\chi}_{_{k_{0},k}}$ 
do in fact satisfy the modular invariance
property in the sense of \cite{DLiM4}. Moreover, it can be shown that 
$\tilde{\chi}_{_{k_{0},k}}$ is a modular function for
some subgroup of finite index of the modular group (\cite{Bo}). 
\begin{rema+}
  The above construction of $\nu$-twisted modules from untwisted
  modules follows from \cite[Proposition 5.4]{Li2}. A different but
  essentially equivalent construction was given in \cite{Li3} by using the
  restricted dual and a certain automorphism, a procedure which in our
  case basically amounts to the change of Virasoro element. 
  Since in the present situation all finite order
  automorphisms are inner, the previous constructions are
  actually a reflection of the fact that under the deformed
  action associated with the semisimple weight one primary vector 
  $\lambda h=\lambda h_{-1}\mathbf{1}$, $\lambda \in \mathbf{Q}$, the
  simple (adjoint) module $V_{k}(A_{1}^{(1)})$ is a 
  so-called $G$-simple current in physical terminology, where $G$ is
  any torsion subgroup of $\text{Aut}\big(V_{k}(A_{1}^{(1)})\big)$
  that contains $\exp(2\pi i\lambda h(0))$ (cf.~\cite{DLiM3}). Although such 
deformed actions also work for general modules, the simple current modules 
are particularly important from a physical point of view. This is mainly 
because they give rise to 
a tensor functor which acts as a permutation on the set of equivalence classes 
of irreducible weak modules (the associated matrix of the left multiplication 
of the equivalence class of a simple current module with respect to the 
standard basis of the Verlinde algebra is a permutation). We refer to 
\cite{DLiM3}, \cite{Li3}, and \cite{Li4} for an in-depth discussion of 
simple current (twisted) modules and various simple current extensions of 
VOAs.
\end{rema+} 

We now concentrate on the objects in $\mathcal{O}_{2}(k)$, namely the 
spaces $\Omega_{k_{0},2k+1}$ defined in
\eqref{2.2.10}. More specifically, we show that these spaces 
are twisted modules for a certain VOA lying inside
$L\big((2k+1)\Lambda_{0};A_{2}^{(1)}\big)$ in such a way that their
$q$-traces are well defined and actually equal to
$\chi_{_{k_{0},k}}^{\nu}(q^{1/2})$.

Recall the setting and notations of \S 2.1-2.2. In particular,
$\mu$ is the minus transpose map of $\mathfrak{sl}(3,\mathbf{C})$, and
$\mathfrak{sl}(3,\mathbf{C})_{[i]}$ is the $(-1)^{i}$-eigenspace of
$\mu$, $i=0,1$. Since $\mathfrak{sl}(3,\mathbf{C})_{[0]}\cong 
\mathfrak{sl}(2,\mathbf{C})$, we may choose a basis $\{E_{1}, F_{1},
H_{1}\}$ of $\mathfrak{sl}(3,\mathbf{C})_{[0]}$ such that
$[H_{1},E_{1}]=2E_{1}$, $[H_{1},F_{1}]=-2F_{1}$, and
$[E_{1},F_{1}]=H_{1}$. Then $H_{0}:=[E_{0},F_{0}]=-H_{1}/2$, where
$E_{0}$ and $F_{0}$ are a lowest respectively highest weight vector of the
$\mathfrak{sl}(3,\mathbf{C})_{[0]}$-module
$\mathfrak{sl}(3,\mathbf{C})_{[1]}$, normalized so that
$[H_{0},E_{0}]=2E_{0}$. The CSA $\mathfrak{s}$ of 
$\mathfrak{sl}(3,\mathbf{C})$ is
chosen to be the principal one, that is 
\begin{equation}\label{3.2.13}
  \mathfrak{s}=\text{Cent}_{\mathfrak{sl}(3,\mathbf{C})}\big(E_{0}+E_{1}\big)
=\mbox{$\mathbf{C}$-span}\,\{E_{0}+E_{1},2F_{0}+F_{1}\}.
\end{equation}
The form $\langle \cdot,\!\cdot \rangle$ being normalized as in \S 2.1,
one has that 
{\allowdisplaybreaks 
\begin{equation}\label{3.2.14}
 \langle E_{0},\!F_{0} \rangle =1,\, \langle E_{1},\!F_{1} \rangle
  =4,\, \langle H_{1},\!H_{1} \rangle =8,\, [H_{1},E_{0}]=-4E_{0},\, 
[H_{1},F_{0}]=4F_{0}.
\end{equation}
Let $\eta=\exp(\pi i/3)$ and recall from \eqref{2.1.1} the principal
automorphism $\sigma$ of $\mathfrak{sl}(3,\mathbf{C})$:}
\begin{equation}\label{3.2.15}
  \sigma(H_{i})=H_{i},\, \sigma(E_{i})=\eta E_{i},\, i=0,1.
\end{equation}
Equivalently, $\sigma=\exp\big(\pi i\text{ad}(H_{1})/6\big)\mu$. Define the
following elements of $\mathfrak{sl}(3,\mathbf{C})$:
{\allowdisplaybreaks 
\begin{equation}\label{3.2.16}
  \begin{split}
  & a_{1}=\frac{1}{3}(E_{0}+E_{1}),\,
  a_{2}=\frac{1}{2}(2F_{0}+F_{1}),\,
  a_{3}=\frac{1}{3}(2E_{0}-E_{1}),\,
  a_{4}=\frac{1}{4}(4F_{0}-F_{1}),\\
  & a_{5}=\frac{1}{2}[E_{1},E_{0}],\,
  a_{6}=\frac{1}{2}[F_{0},F_{1}],\, a_{7}=\frac{1}{2\sqrt{2}}H_{1},\,
  a_{8}=\frac{1}{2\sqrt{6}}\big[E_{1},[E_{1},E_{0}]\big],
\end{split}
\end{equation}
and let $\pi =(12)(34)(56)\in \mathcal{S}_{8}$.} Then $\{a_{i}\mid i\in
\{1,\ldots,8\} \}$ and $\{a_{\pi (i)}\mid i\in
\{1,\ldots,8\} \}$ are $\sigma$-homogeneous dual bases of
$\mathfrak{sl}(3,\mathbf{C})$ such that $\sigma(a_{i})=\eta a_{i}$,
$i=1,3$, $\sigma(a_{i})=\eta^{5} a_{i}$, $i=2,4$,
$\sigma(a_{5})=\eta^{2} a_{5}$, $\sigma(a_{6})=\eta^{4} a_{6}$,
$\sigma(a_{7})=a_{7}$, $\sigma(a_{8})=\eta^{3} a_{8}$, and
$\{a_{1},a_{2}\}$, $\{a_{2},a_{1}\}$ are $\sigma$-homogeneous dual
bases of $\mathfrak{s}$. Consider now the following subalgebras of
$\widetilde{\mathfrak{sl}(3,\mathbf{C})}$:
\begin{equation}\label{3.2.17}
  \hat{\mathfrak{s}}=\mathfrak{s}\otimes \mathbf{C}[t,t^{-1}]\oplus
  \mathbf{C}c,\, \tilde{\mathfrak{s}}=\hat{\mathfrak{s}}\rtimes
  \mathbf{C}d,\, \tilde{\mathfrak{s}}'=[\tilde{\mathfrak{s}},
\tilde{\mathfrak{s}}],\,\tilde{\mathfrak{s}}'_{\pm}=\mathfrak{s}\otimes 
t^{\pm 1}\mathbf{C}[t^{\pm 1}].
\end{equation}
Notice that $\tilde{\mathfrak{s}}'=\tilde{\mathfrak{s}}'_{+}\oplus
\tilde{\mathfrak{s}}'_{-}\oplus \mathbf{C}c$ is a Heisenberg Lie
subalgebra of $\hat{\mathfrak{s}}$, and form the induced irreducible 
$\hat{\mathfrak{s}}$-module (which is irreducible even as an  
$\tilde{\mathfrak{s}}'$-module)
\begin{equation}\label{3.2.18}
  M(2k+1)=U(\hat{\mathfrak{s}})\otimes_{U(\mathfrak{s}\otimes
    \mathbf{C}[t]\oplus \mathbf{C}c)}\mathbf{C},
\end{equation}
where $\mathfrak{s}\otimes
    \mathbf{C}[t]$ acts trivially on $\mathbf{C}$ while $c$ acts as
    multiplication by $2k+1$. Recall that $M(2k+1)\cong
    S(\tilde{\mathfrak{s}}'_{-})$ as vector spaces, and that the action of
    $\hat{\mathfrak{s}}$ extends to $\tilde{\mathfrak{s}}$ by letting $d$
    act as the degree operator, so that $M(2k+1)$ acquires thereby a
    $\mathbf{Z}$-gradation.
\begin{rema+}
  Notice that $\mathfrak{s}_{(1)}=\mathbf{C}a_{1}$ and 
  $\mathfrak{s}_{(5)}=\mathbf{C}a_{2}$, and let 
  \begin{equation}\label{3.2.19}
    \hat{\mathfrak{s}}[\sigma]=\mathfrak{s}_{(1)}\otimes
    t^{1/6}\mathbf{C}[t,t^{-1}]\oplus
    \mathfrak{s}_{(5)}\otimes
    t^{5/6}\mathbf{C}[t,t^{-1}]\oplus \mathbf{C}c.
  \end{equation}
  Then
  $\hat{\mathfrak{s}}[\sigma]=\tilde{\mathfrak{s}}[\sigma]'=\big[
\tilde{\mathfrak{s}}[\sigma],\tilde{\mathfrak{s}}[\sigma]\big]$, where 
$\tilde{\mathfrak{s}}[\sigma]\subset \widetilde{\mathfrak{sl}(3,\mathbf{C})}
[\sigma]$ is defined in \eqref{2.2.8}.
\end{rema+}  

Let now $V_{2k+1}(A_{2}^{(1)})=L\big((2k+1)\Lambda_{0};A_{2}^{(1)}\big)$. When
considered as a $\tilde{\mathfrak{s}}'$-module,
$V_{2k+1}(A_{2}^{(1)})$ decomposes by \cite[Theorem 1.7.3]{FLM} as 
$V_{2k+1}(A_{2}^{(1)})=M(2k+1)\otimes \Omega_{2k+1}\cong
  S(\tilde{\mathfrak{s}}_{-}')\otimes \Omega_{2k+1}$, 
where $\Omega_{2k+1}$ is the vacuum subspace
of $V_{2k+1}(A_{2}^{(1)})$ for the action of
$\tilde{\mathfrak{s}}'$. Set
{\allowdisplaybreaks 
\begin{align}
  & \omega'=\frac{1}{4(k+2)}\sum_{i=1}^{8}a_{i}(-1)a_{\pi
    (i)}(-1)\mathbf{1}\in V_{2k+1}(A_{2}^{(1)}),\label{3.2.20} \\
  & \omega_{1}=\frac{1}{2(2k+1)}\sum_{i=1}^{2}a_{i}(-1)a_{\pi
    (i)}(-1)\mathbf{1}\in M(2k+1).\label{3.2.21}
\end{align}
By Theorems 3.1.4 and 3.1.5,
$\big(V_{2k+1}(A_{2}^{(1)}),Y,\mathbf{1},\omega'\big)$ is a
$\sigma$-rational VOA of rank $c(k):=\frac{4(2k+1)}{k+2}$ and $\big\{
\big(L_{k_{0},2k+1}\big(\widehat{\mathfrak{sl}(3,\mathbf{C})}[\sigma]\big),
Y^{\sigma}\big)\mid
k_{0}\in \{0,1,\ldots,k\} \big\}$ is the complete set of its
simple $\sigma$-twisted modules, the map $Y^{\sigma}$ being as in 
\eqref{3.1.5}.} On the other hand, $M(2k+1)$ is
stable under $Y(a,z)\big|_{M(2k+1)}$ for all $a\in M(2k+1)$, and it is
not difficult to prove that when equipped with the restricted map,
$\big(M(2k+1),Y,\mathbf{1},\omega_{1}\big)$ becomes a VOA of rank
$\dim \mathfrak{s}=2$ (cf., e.g., \cite{DL1}, \cite{Li1}). Moreover,
$\big(M(\sigma;2k+1),Y^{\sigma}\big)$ is a simple $\sigma$-twisted
$M(2k+1)$-module, where $M(\sigma;2k+1)$ is the irreducible
$\tilde{\mathfrak{s}}[\sigma]'$-module defined in \eqref{2.2.9}.

Let $\{\beta_{1},\beta_{2}\}\subset \mathfrak{s}^{*}$ be the basis of
the root system  
of $\mathfrak{sl}(3,\mathbf{C})$ with associated coroots $\{a_{1}+a_{2}, \eta
^{2}a_{1}-\eta a_{2}\}$. Form the root lattice 
$Q=\mathbf{Z}\beta_{1}+\mathbf{Z}\beta_{2}$ and note that 
\begin{equation}\label{3.2.22}
  \Omega_{2k+1}=\coprod_{\beta\in Q}\Omega_{2k+1}^{\beta},\, 
\Omega_{2k+1}^{\beta}=\{v\in \Omega_{2k+1}\mid
  a(0)v=\beta(a)v \text{ for } a\in \mathfrak{s} \}.
\end{equation}
Then
{\allowdisplaybreaks 
\begin{align}
  \Omega_{2k+1}^{0} &  = \{v\in V_{2k+1}(A_{2}^{(1)})\mid a(n)v=0
  \text{ for } a\in \mathfrak{s},\, n\ge 0 \}\label{3.2.23} \\
  & = \{v\in V_{2k+1}(A_{2}^{(1)})\mid a(n)v=0
  \text{ for } a\in M(2k+1),\, n\ge 0 \}, \notag
\end{align}
and it follows from \eqref{3.2.20} and \eqref{3.2.21} that}
\begin{equation}\label{3.2.24}
  \omega_{2}:=\omega'-\omega_{1}\in \Omega_{2k+1}^{0}.
\end{equation}
Let $Y(\omega',z)=\sum_{n\in \mathbf{Z}}L(n)z^{-n-2}$ and $Y(\omega_{i},z)
=\sum_{n\in \mathbf{Z}}L^{i}(n)z^{-n-2}$, $i=1,2$, 
so that $L(n)=L^{1}(n)+L^{2}(n)$. Then $L^{1}(-1)b=0$ for every $b\in
\Omega_{2k+1}^{0}$ by \eqref{3.2.23}, hence
$Y(L^{2}(-1)b,z)=Y(L(-1)b,z)=\frac{d}{dz}Y(b,z)$. Furthermore,
\eqref{3.2.21} implies that $L(1)\omega_{1}=0$ and then by 
\cite[Theorem 5.1]{FZ} one gets that $[L^{1}(m),L^{2}(n)]=0$ and
{\allowdisplaybreaks
\begin{equation*}
  [L^{2}(m),L^{2}(n)]=(m-n)L^{2}(m+n)+\frac{1}{12}(m^{3}-m)
\delta_{m+n,0}c_{2}(k) \text{ for } m, n\in \mathbf{Z},
\end{equation*}
where $c_{2}(k):=c(k)-2=2c_{1}(k)$.} By \eqref{3.1.1} and \eqref{3.2.23}
 $\Omega_{2k+1}^{0}$ is 
closed under vertex operators so that  
when equipped with the
restricted map, $\big(\Omega_{2k+1}^{0},Y,\mathbf{1},\omega_{2}\big)$
becomes a VOA of rank $c_{2}(k)$. 
\begin{rema+}
  (i) Since $\big(\Omega_{2k+1}^{\beta},Y\big)$, 
  $\beta\in Q$, are inequivalent simple $\Omega_{2k+1}^{0}$-modules, the 
VOA $\Omega_{2k+1}^{0}$ is not rational. It is also easy to see that 
$\Omega_{2k+1}^{0}$ is in fact a simple VOA.\newline 
  \noindent
  (ii) It follows from \cite[Theorem 5.2]{FZ} that $\Omega_{2k+1}^{0}$ can
  equivalently be defined as the set $\{v\in V_{2k+1}(A_{2}^{(1)})\mid
  L^{1}(-1)v=0\}$. Using the terminology of
  \cite{FZ}, we may call $\Omega_{2k+1}^{0}$ the commutant of 
  $M(2k+1)$ (the commutant 
  construction corresponds to the coset construction in physics). Notice also 
  that $\Omega_{2k+1}^{0}$ has no 
  weight one vectors, as it can be seen directly
  from the structure of $N((2k+1)\Lambda_{0})$ (cf.~(A.3) in the Appendix
  for $k=1$).
\end{rema+}

Fix $k_{0}\in \{0,1,\ldots,k\}$ and denote the vertex operators
on $L_{k_{0},2k+1}\big(\widehat{\mathfrak{sl}(3,\mathbf{C})}[\sigma]\big)$
associated to $\omega'$ and $\omega_{i}$, $i=1,2$, respectively by
\begin{equation}\label{3.2.25}
  Y^{\sigma}(\omega',z)=\sum_{n\in \mathbf{Z}}L(n)z^{-n-2},\,
  Y^{\sigma}(\omega_{i},z)=\sum_{n\in \mathbf{Z}}L^{i}(n)z^{-n-2},\, i=1,2,
\end{equation}
so that $L(n)=L^{1}(n)+L^{2}(n)$. Then $[L^{1}(m),L^{2}(n)]=0$, 
$m,n\in \mathbf{Z}$, and
 $\Omega_{k_{0},2k+1}$ is stable under $Y^{\sigma}(a,z)$ for
all $a\in \Omega_{2k+1}^{0}$ by \eqref{3.1.1}. Let $v_{k_{0}}\in 
\Omega_{k_{0},2k+1}$
be a highest weight vector of 
$L_{k_{0},2k+1}\big(\widehat{\mathfrak{sl}(3,\mathbf{C})}[\sigma]\big)$. 
Recall that $L(0)$ induces a $\frac{1}{6}\mathbf{Z}$-gradation on
this module and let $(\Omega_{k_{0},2k+1})_{s}=\{v\in \Omega_{k_{0},2k+1}\mid
  L^{2}(0)v=sv\}$, $s\in \mathbf{C}$. Since $L^{1}(0)$ acts as a scalar
operator on $\Omega_{k_{0},2k+1}$ (cf. ~\eqref{3.2.27} below), it follows that
$\big(\Omega_{k_{0},2k+1},Y^{\sigma}\big)$ is a
$\frac{1}{6}\mathbf{Z}$-graded weak $\sigma$-twisted 
$\Omega_{2k+1}^{0}$-module such that 
\begin{equation*}
\Omega_{k_{0},2k+1}=\coprod_{n\in
    \frac{1}{6}\mathbf{N}}(\Omega_{k_{0},2k+1})_{n+\lambda(k_{0})}, 
\,\dim  (\Omega_{k_{0},2k+1})_{\lambda(k_{0})} =1,\, 
\dim  (\Omega_{k_{0},2k+1})_{n+\lambda(k_{0})} < \infty, \, n\in
    \tfrac{1}{6}\mathbf{N}, 
\end{equation*}
where $\lambda(k_{0})\in \mathbf{Q}$ is the lowest weight of
$\Omega_{k_{0},2k+1}$ given by
$L^{2}(0)v_{k_{0}}=\lambda(k_{0})v_{k_{0}}$. Thus 
$\big(\Omega_{k_{0},2k+1},Y^{\sigma}\big)$ is in fact a
$\sigma$-twisted
$\big(\Omega_{2k+1}^{0},Y,\mathbf{1},\omega_{2}\big)$-module whose 
$q$-trace is 
{\allowdisplaybreaks 
\begin{align}
  f_{_{k_{0},2k+1}}(q): & =\text{tr}_{_{\Omega_{k_{0},2k+1}}}q^{L^{2}(0)-
\frac{c_{2}(k)}{24}}=q^{\lambda(k_{0})-\frac{c_{2}(k)}{24}}\sum_{n=0}^{\infty}
\Big(\!\dim \big(\Omega_{k_{0},2k+1}\big)_{\frac{n}{6}+\lambda(k_{0})}\Big)
  q^{\frac{n}{6}}\label{3.2.26} \\
  & =q^{\lambda(k_{0})-\frac{c_{2}(k)}{24}}\dim_{*} \Omega_{k_{0},2k+1},\notag
\end{align}
where $\dim_{*}
  \Omega_{k_{0},2k+1}$ is as in \eqref{2.2.11}.} In order to compute
  $\lambda(k_{0})$, let us again denote the canonical
  generators of $\widehat{\mathfrak{sl}(3,\mathbf{C})}[\sigma]$ as in
  Remark 2.1.1. Note that $H_{1}=\frac{1}{3}(h_{1}^{\sigma}-4h_{0}^{\sigma})$,
 and that from Lemma 3.1.6 and the definitions of $\Omega_{k_{0},2k+1}$ and
  $\omega_{1}$ it follows that 
\begin{equation}\label{3.2.27}
  L^{1}(0)\big|_{\Omega_{k_{0},2k+1}}=\frac{5}{72}
\text{id}_{_{\Omega_{k_{0},2k+1}}}.
  \end{equation}
Applying Lemma 3.1.6 to each term in \eqref{3.2.20} one gets that 
{\allowdisplaybreaks  
\begin{equation}\label{3.2.28}
  L(0)v_{k_{0}}=\frac{18k_{0}^{2}-12k_{0}k+2k^{2}+12k_{0}+41k+20}{144(k+2)}
v_{k_{0}}.
\end{equation}
Then \eqref{3.2.27} and \eqref{3.2.28} give}
\begin{equation*}
  \lambda(k_{0})=\frac{18k_{0}^{2}-12k_{0}k+2k^{2}+12k_{0}+31k}{144(k+2)},
\end{equation*}
so that by \eqref{3.2.11} 
{\allowdisplaybreaks  
\begin{equation}\label{3.2.29}
  \lambda(k_{0})-\frac{c_{2}(k)}{24}=\frac{1}{2}\Big(h_{k_{0},k}-
\frac{\tilde{c}_{1}(k)}{24}\Big)=\frac{1}{2}\Big(h_{k_{0},k}^{\nu}-
\frac{c_{1}(k)}{24}\Big).
\end{equation}
Plugging \eqref{3.2.29} in \eqref{3.2.26} and comparing with
\eqref{3.2.6} and \eqref{3.2.9}, one concludes that}
\begin{equation}\label{3.2.30}
  f_{_{k_{0},2k+1}}(q)=\chi_{_{k_{0},k}}(q^{1/2})=\chi_{_{k_{0},k}}^{\nu}
(q^{1/2}) \text{ for } k_{0}\in \{0,1,\ldots,k\}.
\end{equation}

Summarizing, we proved
{\allowdisplaybreaks
\begin{theo+}
  $\big(\Omega_{2k+1}^{0},Y,\mathbf{1},\omega_{2}\big)$ is a simple VOA of
  rank $c_{2}(k)=2c_{1}(k)$. Moreover, for every $k_{0}\in \{0,1,\ldots,k\}$, 
  $\big(\Omega_{k_{0},2k+1},Y^{\sigma}\big)$ is a $\sigma$-twisted 
  $\Omega_{2k+1}^{0}$-module whose $q$-trace
  $f_{_{k_{0},2k+1}}(q)$ satisfies 
  $f_{_{k_{0},2k+1}}(q)=\chi_{_{k_{0},k}}^{\nu}(q^{1/2})$, 
  where
  $\chi_{_{k_{0},k}}^{\nu}(q)$ is as in Theorem 3.2.1.\hfill 
$\Box$
\end{theo+}  

\begin{remark}
  As in the $A_{1}^{(1)}$-case, one could at first consider 
  $L_{k_{0},2k+1}\big(A_{2}^{(2)}\big)$ just as a simple $\mu$-twisted
  $V_{2k+1}\big(A_{2}^{(1)}\big)$-module, and then use Proposition
  3.1.2 in order to modify $\omega'$ appropriately. Indeed, let
  $\tilde{\omega}'=\omega'+\tfrac{1}{12}H_{1}(-2)\mathbf{1}$. It is easily 
  checked that 
  $\big(V_{2k+1}(A_{2}^{(1)}), Y, \mathbf{1}, \tilde{\omega}'\big)$
  becomes a $\tfrac{1}{6}\mathbf{N}$-graded VOA of rank
  $\tfrac{2(2k+1)(4-k)}{3(k+2)}$ such that its $q$-trace coincides 
  up to a power of $q$  
  with the rescaled $(4,1,1)$-specialization of
  $e^{-(2k+1)\Lambda_{0}}\text{ch}\,L\big((2k+1)\Lambda_{0};A_{2}^{(1)}\big)$.
 Then $\tilde{\omega}'_{1}$ induces the rescaled principal 
gradation on $L_{k_{0},2k+1}\big(A_{2}^{(2)}\big)$, as needed. However, 
in this 
``picture'' one does not get a representation of $\hat{\mathfrak{s}}[\sigma]$ 
on $L_{k_{0},2k+1}\big(A_{2}^{(2)}\big)$ by means of $\mu$-twisted vertex 
operators parametrized by some elements of $V_{2k+1}\big(A_{2}^{(1)}\big)$ 
(see also (A.2) in the Appendix). In order to achieve this, one would have 
to deform the $\mu$-twisted vertex operators into $\sigma$-twisted vertex 
operators as in \cite[Proposition 5.4]{Li2}, which reduces essentially to 
the viewpoint adopted above.
\end{remark}

Let us now see how appropriate the structures of Theorem 3.2.1 and
Theorem 3.2.5 are
in the context of Problems 1 and 2.} As shown in the Appendix for $k=1$,
$\Omega_{2k+1}^{0}$ and $V_{k}(A_{1}^{(1)})$ may in fact be nonisomorphic even
as $\mathbf{Z}$-graded vector spaces. However, it is important
to notice that the $q$-traces $\chi_{_{k_{0},k}}^{\nu}(q)$ and
$f_{_{k_{0},2k+1}}(q)$ are actually equal only up to the transformation
$q\rightarrow q^{\frac{1}{2}}$. Moreover, both the ranks of the VOAs
$V_{k}(A_{1}^{(1)})$ and $\Omega_{2k+1}^{0}$ and the orders of the
automorphisms $\nu$ and $\sigma$ differ by the same factor
$2$. Obviously, it would be preferable not to have these differences,
if possible. This can actually be arranged by using the
recently developed permutation orbifold theory (\cite{BDM}), as we
shall now explain.

Recall from \eqref{3.1.9} that the tensor product VOA
$\big(V_{k}(A_{1}^{(1)})^{\otimes 2},
Y_{_{\otimes}},\mathbf{1}_{_{\otimes}},\omega_{_{\otimes}}\big)$ has
rank $c_{2}(k)=2c_{1}(k)$, and notice that the diagonal action of
$\text{Aut}\big(V_{k}(A_{1}^{(1)})\big)$ on
$V_{k}(A_{1}^{(1)})^{\otimes 2}$ commutes with the action of the
symmetric group $\mathcal{S}_{2}$ on $V_{k}(A_{1}^{(1)})^{\otimes
  2}$. Let $\psi$ be the generator of $\mathcal{S}_{2}$ and form the
6-th order automorphism $\tau=\nu^{2}\psi\in
\text{Aut}\big(V_{k}(A_{1}^{(1)})^{\otimes 2}\big)$. Given a
$\nu$-twisted $V_{k}(A_{1}^{(1)})$-module $M$, one may define a structure
of $\tau$--twisted $V_{k}(A_{1}^{(1)})^{\otimes 2}$-module on $M$ in
the following way. By \cite[Proposition 2.1.1]{Hua}, the identity
\begin{equation}
  \Bigg(\!\exp \!\Bigg(-\sum_{j\in
    \mathbf{Z}_{+}}C_{j}x^{j+1}\frac{d}{dx}\!\Bigg)\!\Bigg)\cdot
  x=x+\frac{x^{2}}{2}\label{3.2.31}
\end{equation}
uniquely determines a sequence $\{C_{j}\}_{_{j\in
    \mathbf{Z}_{+}}}\in \mathbf{Q}^{\infty}$, and one defines 
  the operator
\begin{equation}
  \Delta_{2}(z)=\exp\!\Bigg(\sum_{j\in
    \mathbf{Z}_{+}}C_{j}z^{-\frac{j}{2}}L(j)\!\Bigg)2^{-L(0)}
z^{-\frac{1}{2}L(0)}\in \big(\text{End}\,V_{k}(A_{1}^{(1)})\big)
\big[\big[z^{1/2},z^{-1/2}\big]\big],\label{3.2.32}
\end{equation}
where $L(j)=\text{Res}_{z}z^{j+1}Y(\omega,z)$, $j\in \mathbf{Z}_{+}$,
and $2^{-L(0)}z^{-\frac{1}{2}L(0)}$ acts on each homogeneous subspace
$V_{k}(A_{1}^{(1)})_{n}$ as $2^{-n}z^{-n/2}$, $n\in
\mathbf{Z}$. Recall that $\eta=\exp(\pi i/3)$ and let $r\in \{0,1,2\}$
and $a\in V_{k}(A_{1}^{(1)})$ be such that $\nu (a)=\eta^{-2r}a$. 
Define the following operators:
{\allowdisplaybreaks
  \begin{equation}\label{3.2.33}
    \begin{split}
      & Y^{\tau}(a\otimes
      \mathbf{1},z)=Y^{\nu}(\Delta_{2}(z)a,z^{1/2})\in (\text{End}\,
      M)\big[\big[z^{1/6},z^{-1/6}\big]\big],\\
      & Y^{\tau}(\mathbf{1}\otimes
      a,z)=\eta^{2}\lim_{z^{1/6}\rightarrow \eta z^{1/6}}Y^{\tau}(a\otimes
      \mathbf{1},z)\in (\text{End}\,
      M)\big[\big[z^{1/6},z^{-1/6}\big]\big].
    \end{split}
  \end{equation}
It was proved in \cite{BDM} that $\big\{Y^{\tau}(a\otimes
      \mathbf{1},z),\,Y^{\tau}(\mathbf{1}\otimes
      b,z)\mid a,b\in V_{k}(A_{1}^{(1)})\big\}$ generates 
a local system of $\tau$-twisted vertex operators in the sense
of \cite{Li2}. By using the theory of {\em loc. ~cit.}, it was further shown in
\cite[\S 3 and \S 6]{BDM} that the vertex map $Y^{\tau}$ extends to
the whole space $V_{k}(A_{1}^{(1)})^{\otimes 2}$ in such a way that
 $(M,Y^{\tau})$ becomes a $\tau$-twisted
$V_{k}(A_{1}^{(1)})^{\otimes 2}$-module.} Then one can prove:
{\allowdisplaybreaks
\begin{theorem}
  The VOA $\big(V_{k}(A_{1}^{(1)})^{\otimes 2},
Y_{_{\otimes}},\mathbf{1}_{_{\otimes}},\omega_{_{\otimes}}\big)$ of 
rank $c_{2}(k)$ is $\tau$-rational and its irreducible $\tau$-twisted
modules are exactly
$\big(L_{k_{0},k}\big(\widehat{\mathfrak{sl}(2,\mathbf{C})}[\nu]\big),
Y^{\tau}\big)$,
$k_{0}\in \{0,1,\ldots,k\}$. Moreover, the $q$-trace
$\chi_{_{k_{0},k}}^{\tau}(q)$ of
$\big(L_{k_{0},k}\big(\widehat{\mathfrak{sl}(2,\mathbf{C})}[\nu]\big),
Y^{\tau}\big)$
coincides with the $q$-trace $f_{_{k_{0},2k+1}}(q)$ of the
$\sigma$-twisted $\Omega_{2k+1}^{0}$-module
$\big(\Omega_{k_{0},2k+1},Y^{\sigma}\big)$, and one has
\begin{equation}
  \chi_{_{k_{0},k}}^{\tau}(q)=f_{_{k_{0},2k+1}}(q)=\chi_{_{k_{0},k}}^{\nu}
(q^{1/2}),\, k_{0}\in \{0,1,\ldots,k\}, \label{3.2.34}
\end{equation}
where $\chi_{_{k_{0},k}}^{\nu}(q)$ is as in Theorem 3.2.1.
\end{theorem}}
{\allowdisplaybreaks
\begin{proof}
  The first statement follows from Theorems 6.8, 6.4, and 3.8 in
  \cite{BDM}, so that it remains to prove \eqref{3.2.34}. Recall that
  $\omega_{_{\otimes }}=\omega \otimes \mathbf{1}+\mathbf{1}\otimes
  \omega$ and write $Y^{\tau}(\omega_{_{\otimes }},z)=\sum_{n\in
    \mathbf{Z}}L^{\tau}(n)z^{-n-2}$. By \eqref{3.2.33} one has
  \begin{equation}
    Y^{\tau}(\omega_{_{\otimes }},z)=Y^{\nu}(\Delta_{2}(z)\omega
    ,z^{1/2})+\eta^{2}\lim_{z^{1/6}\rightarrow \eta
      z^{1/6}}Y^{\nu}(\Delta_{2}(z)\omega ,z^{1/2}).\label{3.2.35} 
  \end{equation}
  The constants $C_{j}$ in \eqref{3.2.31} may be computed reccurrently  
for small values of $j$. For instance,
  $C_{1}=-\frac{1}{2}$ and $C_{2}=\frac{1}{4}$. Moreover, the Virasoro 
algebra commutation relations for the VOA $V_{k}(A_{1}^{(1)})$ 
imply that $L(2)\omega=\tfrac{c_{1}(k)}{2}\mathbf{1}$ and  
$L(1)\omega=L(j)\omega=0$ for $j\ge 3$ (cf.~\cite{FHL}). 
Then \eqref{3.2.32} yields
  $$\Delta_{2}(z)\omega = \frac{z^{-1}}{4}\Big(\omega
  +\frac{c_{1}(k)C_{2}}{2}z^{-1}\Big)=\frac{z^{-1}}{4}\Big(\omega
  +\frac{c_{1}(k)}{8}z^{-1}\Big),$$
  so that
  \begin{equation}
    Y^{\nu}(\Delta_{2}(z)\omega
    ,z^{1/2})=\frac{z^{-1}}{4}Y^{\nu}(\omega,z^{1/2})+\frac{c_{1}(k)}{32}
\text{id}_{_{M_{k_{0},k}}},\label{3.2.36}
  \end{equation}
  where 
  $M_{k_{0},k}:=L_{k_{0},k}\big(\widehat{\mathfrak{sl}(2,\mathbf{C})}
[\nu]\big)$. Plugging \eqref{3.2.36} in
  \eqref{3.2.35} and then extracting the coefficient of $z^{-2}$ from
  \eqref{3.2.35} one gets that $L^{\tau}(0)=\frac{1}{2}L(0)+\frac{c_{1}(k)}
{16}\text{id}_{_{M_{k_{0},k}}}$, 
  and thus
  $$L^{\tau}(0)-\frac{c_{2}(k)}{24}\text{id}_{_{M_{k_{0},k}}}=\frac{1}{2}
\Big(L(0)-\frac{c_{1}(k)}{24}\text{id}_{_{M_{k_{0},k}}}\!\Big)$$
  since $c_{2}(k)=2c_{1}(k)$. It then follows from \eqref{3.2.9} that
  $$\text{tr}_{_{M_{k_{0},k}}}q^{L^{\tau}(0)-\frac{c_{2}(k)}{24}}=
\text{tr}_{_{M_{k_{0},k}}}q^{\frac{1}{2}\big(\!L(0)-\frac{c_{1}(k)}{24}\!
\big)}=\chi_{_{k_{0},k}}^{\nu}(q^{1/2}),$$
  which is the same as \eqref{3.2.34}.
\end{proof}}
\begin{remark}
  (i) As used here, a $\nu$-twisted
  $V_{k}(A_{1}^{(1)})$-module is the
  same as a $\nu^{-1}$-twisted $V_{k}(A_{1}^{(1)})$-module in the
  sense of \cite{BDM}. Note that $\nu^{-1}=\nu^{2}$ in our
  case.\newline
  (ii) An interesting similarity between the automorphisms $\tau$
  and $\sigma$ is that both of them are
  compositions of two commuting automorphisms of orders 3 and 2 in 
  $\text{Aut}\big(V_{k}(A_{1}^{(1)})^{\otimes 2}\big)$ and
  $\text{Aut}\big(\Omega_{2k+1}^{0}\big)$ 
  respectively (cf.~\eqref{3.2.15}). Note also that in both 
  cases the order 3 automorphism is inner and the involution is outer.
\end{remark}

In view of Theorems 3.2.5 and 3.2.7, one may ask whether the VOAs
$V_{k}(A_{1}^{(1)})^{\otimes 2}$ and $\Omega_{2k+1}^{0}$ are
isomorphic. As shown in
the Appendix for $k=1$, this is not necessarily true. On the other
hand, nothing was mentioned so far about 
the irreducibility of the $\sigma$-twisted $\Omega_{2k+1}^{0}$-modules
$\Omega_{k_{0},2k+1}$, $k_{0}\in \{0,1,\ldots ,k\}$. Note that in
the level one case -- which is excluded from the present discussion   
-- the $\Omega_{1}^{0}$-module $\Omega_{0,1}$ is trivially irreducible
since $\dim \Omega_{0,1}=1$, as it can be seen from \eqref{2.2.5} and 
\eqref{2.2.11}. The general results on coset constructions and dual pairs 
for VOAs recently obtained in \cite{DM} seem to indicate that these 
$\sigma$-twisted $\Omega_{2k+1}^{0}$-modules are indeed simple, but this 
is highly non-trivial for $k\ge 1$. We
shall discuss such irreducibility questions in a different setting in \S 4.4.
{\allowdisplaybreaks
\begin{remark}
  (i) As a step in the direction suggested in Problem 2, one could
  investigate whether
  $L_{k_{0},k}\big(\widehat{\mathfrak{sl}(2,\mathbf{C})}[\nu]\big)$
  and $\Omega_{k_{0},2k+1}$ are isomorphic $Vir$-modules with
  central charge $c_{2}(k)$ (cf.~Theorems 3.2.5 and 3.2.7). Note
  though that even if true, this 
  would not provide a fully satisfactory answer to Problem 2 
  since 
  $L_{k_{0},k}\big(\widehat{\mathfrak{sl}(2,\mathbf{C})}[\nu]\big)$ 
  and $\Omega_{k_{0},2k+1}$ are not simple 
  $Vir$-modules.\newline
\noindent
(ii) Let $V$ be a VOA with a finite order automorphism $\sigma$ and denote by 
$A_{\sigma}(V)$ the twisted analogue of Zhu's algebra $A(V)$ 
(cf.~\cite{DLiM2}, \cite{Z}). Using the results of \cite{DLiM2} one 
can actually show that 
$A_{\tau}\big(V_{k}(A_{1}^{(1)})^{\otimes 2}\big)\cong
A_{\nu}\big(V_{k}(A_{1}^{(1)})\big)\cong 
A_{\sigma}\big(V_{2k+1}(A_{2}^{(1)})\big)$. 
However, this is more a reformulation of the fact that
$|\mathcal{O}_{1}(k)|=|\mathcal{O}_{2}(k)|=k+1$ rather than a conceptual 
explanation to Theorem 2.2.1.
\end{remark}}

\section{The GVOA $\Omega_{2k+1}^{A}$ and its action on the spaces in 
$\mathcal{O}_{2}(k)$}

Continuing our study of Problem 1, we shall now embed
$\Omega_{2k+1}^{0}$ into a larger structure, namely a simple GVOA which 
acts irreducibly on each of the spaces 
$\Omega_{k_{0},2k+1}$ without altering
the $q$-traces $f_{k_{0},2k+1}(q)$ of \eqref{3.2.26} (Theorems 4.4.3 and 
4.4.8). This 
new structure is defined as the quotient $\Omega_{2k+1}^{A}$ of the
vacuum space $\Omega_{2k+1}$ by the action of a certain infinite
abelian group $A$, and it acts on the spaces $\Omega_{k_{0},2k+1}$ by
means of quotient relative $\sigma$-twisted vertex operators 
(introduced in \cite{DL2}) whose component operators generate an algebra that 
includes the $\sigma$-twisted $\mathcal{Z}$-algebra of
\cite{LW}. This is done by using the diagonal action of
$\widehat{\mathfrak{sl}(3,\mathbf{C})}[\sigma]$ on the tensor product
of $2k+1$ copies of the basic standard
$\widehat{\mathfrak{sl}(3,\mathbf{C})}[\sigma]$-module, together with 
the fact that each standard
$\widehat{\mathfrak{sl}(3,\mathbf{C})}[\sigma]$-module of level $2k+1$
can be isomorphically embedded into this tensor product (cf.~\S 2.1). As we 
explain at the end of \S 4.4, the GVOA $\Omega_{2k+1}^{A}$ appears to be a more
appropriate tool than the VOA $\Omega_{2k+1}^{0}$ for further
investigating Problem 1. The main reason for that is the ``equivalence
theorem'' \cite[Theorem 5.5]{LW}, which essentially establishes an 
equivalence between appropriately defined module categories for 
$\widetilde{\mathfrak{sl}(3,\mathbf{C})}[\sigma]$ and 
(an algebra that includes) the $\sigma$-twisted $\mathcal{Z}$-algebra.

The methods used below are based on \cite{DL1}, \cite{DL2}, \cite{L2}, 
\cite{LW}, where a great many technical ingredients were used. Since the 
proofs of our main results in \S 4.4 require 
a careful analysis of most of these ingredients, we first recall 
in \S 4.1-\S 4.3 some of the constructions of {\em loc.~cit.} and adapt 
them to our particular case. We refer to the above-mentioned papers for  
additional details.

Throughout this section we let $l=2k+1$, where $k\in \mathbf{Z}_{+}$ is 
fixed, and we assume that 
$\mathfrak{g}=\mathfrak{sl}(3,\mathbf{C})$ with
the bilinear form  
$\langle \cdot\,,\cdot \rangle$ normalized as in \S 2.1. The notations
introduced below will partially overlap 
with some of the previous sections, as we need to work in the setting
of \cite{DL1} and \cite{DL2}.

\subsection{The new setting}

Recall from \S 3.2 the root basis $\{\beta_{1},\beta_{2}\}$ and
the root lattice $Q$ of
$\mathfrak{g}$, and let $\Phi$ be its root system. Let
$\sigma_{1}$ be the reflection with respect to $\beta_{1}$ and $\mu$ be the
automorphism of $Q$ induced by the Dynkin diagram automorphism
determined by $\mu\beta_{1}=\beta_{2}$,
$\mu\beta_{2}=\beta_{1}$. Recall \eqref{3.2.13}-\eqref{3.2.16} and note that
the linear extension of $\sigma_{1}\mu$ to
$Q\otimes_{\mathbf{Z}}\mathbf{C}$ coincides with the restriction to
$\mathfrak{s}$ of the principal automorphism $\sigma$ of
$\mathfrak{g}$, where $\mathfrak{s}$ and
$Q\otimes_{\mathbf{Z}}\mathbf{C}$ are identified via the form $\langle
\cdot\,,\cdot \rangle$. Since $\sigma_{1}\mu\equiv \sigma$ is a
``twisted Coxeter element'' of order 6, one has
\begin{equation}\label{4.1.1}
  \sum_{p=0}^{5}\sigma^{p}\alpha=0,\quad
  \sum_{p=0}^{5}p\langle \sigma^{p}\alpha\,,\beta \rangle \equiv 0
  \bmod 6 \text{ for $\alpha, \beta\in Q$.}
\end{equation}
Moreover, $Q$ is a positive definite even lattice, $\langle
\sigma^{3}\alpha\,,\alpha \rangle=-\langle \alpha\,,\alpha \rangle \in
2\mathbf{Z}\text{ for }\alpha \in Q$, and $\langle \cdot\,,\cdot
\rangle$ is $\sigma$-invariant. Denote by $\langle \eta \rangle$ the 
cyclic group of order 6 generated by $\eta=\exp(\pi i/3)$. Set
\begin{equation}\label{4.1.2}
  L=Q_{1}\oplus \ldots \oplus Q_{l},
\end{equation}
where each $Q_{i}$ is a copy of $Q$. We write $\alpha_{i} \in Q_{i}$
for the element corresponding to $\alpha \in Q$, and we extend the
form $\langle \cdot\,,\cdot
\rangle$ to $L$ so that $Q_{i}$ and $Q_{j}$ are orthogonal if $i\neq j$.

Consider the $\sigma$-invariant alternating
$\mathbf{Z}$-bilinear maps 
$c_{0}$, $c_{0}^{\sigma}: \, L\times L \longrightarrow 
  \mathbf{Z}/6\mathbf{Z}$ 
defined by $c_{0}(\alpha,\beta)=3\langle\alpha \,, \beta \rangle
  +6\mathbf{Z}$ respectively
  $c_{0}^{\sigma}(\alpha,\beta)=\sum_{p=0}^{5}(3+p)\langle\sigma^{p}\alpha\,, 
\beta \rangle +6\mathbf{Z}$ 
for $\alpha, \beta\in L$. Notice that
$c_{0}^{\sigma}(\alpha,\beta)\equiv 0 \bmod 6$ for every
$\alpha, \beta\in L$, and set
\begin{equation}\label{4.1.3}
  c(\alpha,\beta)=\eta^{c_{0}(\alpha,\beta)}=(-1)^{\langle\alpha,
    \beta \rangle},\quad
  c_{\sigma}(\alpha,\beta)=\eta^{c_{0}^{\sigma}(\alpha,\beta)}=1,\quad
  \alpha, \beta\in L.
\end{equation}
Define also
\begin{equation}\label{4.1.4}
  \epsilon_{0}(\alpha,\beta)=\langle\alpha+\sigma \alpha\,, \beta
  \rangle + 6\mathbf{Z},\quad
  \epsilon(\alpha,\beta)=\eta^{\epsilon_{0}(\alpha,\beta)},\quad 
  \alpha, \beta\in L.
\end{equation}
Then $\epsilon_{0}(\cdot\,,\cdot)$ is a $\sigma$-invariant
2-cocycle on $L$ satisfying
\begin{equation}\label{4.1.5}
  \epsilon_{0}(\alpha,\beta)-\epsilon_{0}(\beta,\alpha)=c_{0}(\alpha,\beta)-
c_{0}^{\sigma}(\alpha,\beta)=c_{0}(\alpha,\beta),
\end{equation}
and it follows from \eqref{4.1.1} together with $\sigma^{3}=-1$ that
$\epsilon_{0}(\alpha,\beta)\equiv 0 \bmod 3$. There\-fore,
$\epsilon(\cdot\,,\cdot)$ is $\sigma$-invariant and
$\epsilon(\alpha,\beta)\in \{\pm 1\}$ for $\alpha,\beta \in L$. Up to 
equivalence, the commutator maps $c$ and
$c_{\sigma}$ uniquely determine two central extensions of $L$ 
\begin{equation}\label{4.1.6}
  1\longrightarrow \langle \eta \rangle \longrightarrow \hat{L}
  \bar{\longrightarrow} L \longrightarrow 1 \,\text{ and }\, 
  1\longrightarrow \langle \eta \rangle \longrightarrow \hat{L}_{\sigma}
  \bar{\longrightarrow} L \longrightarrow 1
\end{equation}  
by the commutator relations 
\label{3.3.27}
\begin{equation}\label{4.1.7}
  aba^{-1}b^{-1}=c(\bar{a}\,,\bar{b}) \text{ for } a,b\in
    \hat{L} \text{ respectively }
  aba^{-1}b^{-1}=c_{\sigma}(\bar{a}\,,\bar{b})=1 \text{ for } a,b\in
    \hat{L}_{\sigma}.
\end{equation}
Notice that 
$\hat{L}_{\sigma}$ is a split extension of $L$, so that
$\hat{L}_{\sigma}=\langle \eta \rangle \times L$. One
has a set-theoretic identification between the groups $\hat{L}$ and 
$\hat{L}_{\sigma}$ such that
the respective group multiplications $\times$ and $\times_{\sigma}$
satisfy
\begin{equation}\label{4.1.8}
  a\times b=\epsilon(\bar{a}\,,\bar{b})\,a\times_{\sigma}b.
\end{equation}
Moreover, $\sigma$ lifts to an automorphism $\hat{\sigma}$ of
$\hat{L}$ such that
\begin{equation}\label{4.1.9}
  \hat{\sigma}\eta =\eta,\quad \overline{\hat{\sigma}a}=\sigma \bar{a}
  \text{ for $a\in \hat{L}$},
\end{equation}
and such a lifting is unique up to multiplication by a lifting of the 
identity automor\-phism of $L$ (cf. ~\cite[Proposition 5.4.1]{FLM}). Since
$\epsilon(\cdot\,,\cdot)$ is $\sigma$-invariant, one gets from
\eqref{4.1.8} that $\hat{\sigma}$ is also an automorphism of
$\hat{L}_{\sigma}$ covering $\sigma$. Let
\begin{align}
  e\, : \, L & \longrightarrow \hat{L}\label{4.1.10}\\
  \alpha & \longmapsto e_{\alpha}=(1,\alpha)\nonumber
\end{align}
be the section corresponding to the cocycle
$\epsilon_{0}(\cdot\,,\cdot)$, i.e., $\bar{e}_{\alpha}=\alpha$, 
  $e_{\alpha}e_{\beta}=\epsilon(\alpha\,,\beta)\,e_{\alpha+\beta}$, 
  $\alpha,\beta \in L$. 
Note that the bilinearity of $\epsilon_{0}(\cdot\,,\cdot)$ implies
that $e_{0}=1 \, (=(1,0))$.  
\begin{remark}
(i) The lattice in
  \eqref{4.1.2} is denoted by $L_{0}$ in \cite[Ch.~13 \&
  14]{DL1}, where one defines $L$ to be the direct sum of $l$ copies
  of the weight lattice of $\mathfrak{g}$ instead (so that $L_{0}$ is
  the dual lattice of $L$ in $L\otimes_{\mathbf{Z}}\mathbf{C}$). For
  our purposes though, it will suffice to use only the root lattice 
  $Q$.\newline
(ii) The 2-cocycle $\epsilon_{0}(\cdot\,,\cdot)$ is cohomologous to the one 
defined in \cite[(2.13)]{DL2}, from which it differs by the 2-coboundary
  $L\times L \ni (\alpha,\beta)\mapsto
  2\langle\alpha,\beta \rangle +6\mathbf{Z}\in \mathbf{Z}/6\mathbf{Z}$.
\end{remark}

Applying the above discussion to a single
copy of the root lattice $Q$, one gets from \cite[Ch.~6]{FLM} and 
\cite[\S 8]{LW} that there exist $x_{\alpha}\in \mathfrak{g}$ ($\alpha \in
\Phi$) such that $\mathfrak{g}=\mathfrak{s}\oplus \coprod_{\alpha \in
    \Phi}\mathbf{C}x_{\alpha}$ and 
{\allowdisplaybreaks 
  $$\sigma x_{\alpha}=x_{\sigma \alpha},\,
  \mathfrak{s}'=0,\, [h,x_{\alpha}]=\langle h\,,\alpha
  \rangle \,x_{\alpha},\, [x_{\alpha},x_{\beta}]=
  \begin{cases}
    \epsilon(\alpha\,,-\alpha)\, \alpha & \text{ if } \alpha+\beta=0,\\
    \epsilon(\alpha\,,\beta)\, x_{\alpha+\beta} & \text{ if } \alpha+\beta \in \Phi,\\
    0 & \text{ if } \alpha+\beta \notin \Phi \cup \{0\},
  \end{cases}$$
for $h\in \mathfrak{s}$ and $\alpha, \beta \in \Phi$.} It is
well-known that the bilinear
form $\langle\cdot\,,\cdot \rangle$ of $\mathfrak{g}$ satisfies
  $$\langle h\,,x_{\alpha}\rangle =0,\quad \langle
  x_{\alpha}\,,x_{\beta}\rangle =
  \begin{cases}
    \epsilon(\alpha\,,-\alpha) & \text{ if }
    \alpha+\beta=0,\\
    0 & \text{ if }
    \alpha+\beta \neq 0
  \end{cases}$$
for $\alpha, \beta \in \Phi$ and $h\in \mathfrak{s}$. Recall the
$\sigma$-decomposition $\mathfrak{g}=\coprod_{j\in \mathbf{Z}/6\mathbf{Z}}
\mathfrak{g}_{(j)}$ and let $\mathfrak{g}_{(n)}=\mathfrak{g}_{(n \bmod 6)}$ 
for $n\in \mathbf{Z}$. Note also that
$\mathfrak{s}=\mathfrak{s}_{(1)}\oplus \mathfrak{s}_{(5)}$ (cf.~
\eqref{3.2.19}). Let further $x_{(n)}$ denote the projection of $x\in
\mathfrak{g}$ on $\mathfrak{g}_{(n)}$, $n\in
\mathbf{Z}$, and set
\begin{equation}\label{4.1.11}
  x(z)=\sum_{n\in
\mathbf{Z}}(x\otimes t^{n})z^{-n-1},\quad x(\sigma;z)=\sum_{n\in
\tfrac{1}{6}\mathbf{Z}}(x_{(6n)}\otimes t^{n})z^{-n-1}.
\end{equation}
We denote by $\langle\cdot\,,\cdot \rangle$ as well the bilinear 
form induced on the space
\begin{equation}\label{4.1.12}
  \mathfrak{h}=L\otimes_{\mathbf{Z}}\mathbf{C}.
\end{equation}

Let now $\mathfrak{h}_{*}$
be a subspace of $\mathfrak{h}$ on which $\langle\cdot\,,\cdot
\rangle$ remains nonsingular, so that
\begin{equation}\label{4.1.13}
  \mathfrak{h}=\mathfrak{h}_{*}\oplus \mathfrak{h}_{*}^{\perp},
\end{equation}
where $^{\perp}$ denotes the orthogonal complement. We also assume
that
\begin{equation}\label{4.1.14}
  \sigma \mathfrak{h}_{*}=\mathfrak{h}_{*}
\end{equation}
and we write
\begin{equation}\label{4.1.15}
 \begin{aligned}
\mathfrak{h} & \longrightarrow \mathfrak{h}_{*}^{\perp}\\
h & \longmapsto h'
\end{aligned}
\text{\quad , \quad}
\begin{aligned}
  \mathfrak{h} & \longrightarrow \mathfrak{h}_{*}\\
  h & \longmapsto h''
\end{aligned}
\end{equation}
for the projection maps to $\mathfrak{h}_{*}^{\perp}$ and
$\mathfrak{h}_{*}$ respectively. By \eqref{4.1.14}, these maps commute with
the action of $\sigma$. In \S 4.2-\S 4.3 we shall 
work with an arbitrary space $\mathfrak{h}_{*}$ satisfying  
\eqref{4.1.13}-\eqref{4.1.14}. This flexibility for the choice of
$\mathfrak{h}_{*}$ will allow us to recover the
usual (unrelativized) untwisted and twisted vertex operators by
taking $\mathfrak{h}_{*}=0$, in which case
the index $_{*}$ will be removed from the notation. In \S 4.4  
though we shall specify $\mathfrak{h}_{*}$ to 
be the image of $\mathfrak{s}$ under its diagonal embedding in
$\mathfrak{h}$.

\subsection{Relative untwisted vertex operators}

Form the induced $\hat{L}$-module and $\mathbf{C}$-algebra
  $\mathbf{C}\{L\}=\mathbf{C}[\hat{L}]\otimes_{\mathbf{C}[\langle \eta
    \rangle ]}\mathbf{C}\cong \mathbf{C}[L]$ (linearly), 
where $\mathbf{C}[\cdot]$ denotes the group algebra and $\eta$ acts on
$\mathbf{C}$ as multiplication by $\eta$ (here $\langle \eta
    \rangle$ is understood as an abstract group disjoint from
    $\mathbf{C}^{\times}$). For $a\in \hat{L}$, let $\iota (a)=a\otimes 1$ be 
    the image of $a$ in $\mathbf{C}\{L\}$. Then the action of
    $\hat{L}$ on $\mathbf{C}\{L\}$ and the product in
    $\mathbf{C}\{L\}$ are given by:
\begin{equation}\label{4.2.1}
  a\cdot \iota (b)=\iota (a) \iota (b)=\iota (ab),\quad \eta \cdot
  \iota (b)=\eta \iota(b)
\end{equation}  
for $a,b\in \hat{L}$. We endow $\mathbf{C}\{L\}$ with the 
$\mathbf{C}$-gradation determined by
\begin{equation}\label{4.2.2}
  \mbox{wt}(\iota (a))=\frac{1}{2}\langle \bar{a}'\,,\bar{a}'\rangle
  \text{ for } a\in \hat{L},
\end{equation}
and we define a grading-preserving action of $\mathfrak{h}$ on
$\mathbf{C}\{L\}$ by
\begin{equation}\label{4.2.3}
  h\cdot \iota (a)=\langle h'\,,\bar{a}\rangle\, \iota (a)
\end{equation}
for $h\in \mathfrak{h}$. The automorphism $\hat{\sigma}$ of $\hat{L}$ 
acts canonically and in a grading-preserving fashion on
$\mathbf{C}\{L\}$ such that
  $\hat{\sigma}\iota (a)=\iota (\hat{\sigma}a)$, $a\in \hat{L}$, 
and $\hat{\sigma}(\iota (a)\iota (b))=\hat{\sigma}(a\cdot \iota
  (b))=\hat{\sigma}(a)\cdot \hat{\sigma}\iota (b)=\hat{\sigma}\iota
  (a)\,\hat{\sigma}\iota (b)$. Then $\mathfrak{h}$ acts as algebra 
derivations and $\hat{\sigma}(h\cdot \iota (a))=\sigma(h)\cdot \hat{\sigma}
\iota (a)$, $h\in \mathfrak{h}$, $a\in \hat{L}$. Define also an action
  $z^{h}\cdot \iota (a)=z^{\langle h'\,,\bar{a}\rangle}\iota (a)$ 
    for $h\in \mathfrak{h}$, $a\in \hat{L}$, so that 
  $\hat{\sigma}(z^{h}\cdot \iota (a))=z^{\sigma(h)}\cdot \hat{\sigma}\iota
  (a)$.

Form the algebras
\begin{equation}\label{4.2.4}
  \hat{\mathfrak{h}}=\mathfrak{h}\otimes \mathbf{C}[t,t^{-1}]\oplus
  \mathbf{C}c,\, \tilde{\mathfrak{h}}=\hat{\mathfrak{h}}\rtimes
  \mathbf{C}d, \, \tilde{\mathfrak{h}}'=[\tilde{\mathfrak{h}},
\tilde{\mathfrak{h}}], \, 
  \tilde{\mathfrak{h}}'_{\pm}=\mathfrak{h}\otimes
  t^{\pm 1}\mathbf{C}[t^{\pm 1}], 
\end{equation}
and notice that $\tilde{\mathfrak{h}}'=\tilde{\mathfrak{h}}'_{+}\oplus 
\tilde{\mathfrak{h}}'_{-}\oplus
  \mathbf{C}c$ is a Heisenberg Lie subalgebra of the affine Lie
  algebra $\tilde{\mathfrak{h}}$. On $\hat{\mathfrak{h}}$ we define a 
weight gradation associated with $\mathfrak{h}_{*}$ by
\begin{equation}\label{4.2.5}
  \mbox{wt}(x\otimes t^{m})=0,\, \mbox{wt}(y\otimes t^{n})=-n,\,
  \mbox{wt}(c)=0
\end{equation}
for $x\in \mathfrak{h}_{*}$, $y\in \mathfrak{h}_{*}^{\perp}$, and
$m,n\in \mathbf{Z}$. By analogy with \eqref{3.2.18}, we
consider the induced irreducible $\tilde{\mathfrak{h}}'$- and 
$\hat{\mathfrak{h}}$-module 
  $M(1,\mathfrak{h})=U(\hat{\mathfrak{h}})\otimes_{U(\mathfrak{h}\otimes 
\mathbf{C}[t]\oplus \mathbf{C}c)}\mathbf{C}\cong S(\tilde{\mathfrak{h}}'_{-})$ 
(linearly), where $\mathfrak{h}\otimes \mathbf{C}[t]$ acts trivially on
$\mathbf{C}$ and $c$ acts as 1. It follows from \eqref{4.2.5} that
$M(1,\mathfrak{h})$ is $\mathbf{Z}$-graded so that
$\mbox{wt}(1)=0$. Moreover, $\sigma$ acts in a natural 
grading-preserving way on
$\hat{\mathfrak{h}}$ (fixing $c$) and on $M(1,\mathfrak{h})$. Set
\begin{equation}\label{4.2.6}
  V_{L}=M(1,\mathfrak{h})\otimes_{\mathbf{C}}\mathbf{C}\{L\}\cong
  S(\tilde{\mathfrak{h}}'_{-})\otimes_{\mathbf{C}}\mathbf{C}[L]\,
  \text{ (linearly)},
\end{equation}
and endow $V_{L}$ with the tensor product $\mathbf{C}$-gradation. In
particular, $\mbox{wt}(\iota (1))=0$, $\mathbf{C}\{L\}$ being naturally
identified with $1\otimes \mathbf{C}\{L\}$. Then $\hat{L}$,
$\tilde{\mathfrak{h}}'$, $\mathfrak{h}$, $z^{h}$ ($h\in \mathfrak{h}$)
act naturally on $V_{L}$ by acting either on $M(1,\mathfrak{h})$ or
$\mathbf{C}\{L\}$ as indicated above. The
automorphism $\hat{\sigma}$ acts in a grading-preserving way
on $V_{L}$ via $\sigma\otimes \hat{\sigma}$, and this action is
compatible with the other actions (cf.~\cite[(3.23)-(3.25)]{DL2}). 
The relations
\begin{equation}\label{4.2.7}
   i(h)=\big(h_{1}(-1)+\ldots +h_{l}(-1)\big)\cdot \iota (1),\, h\in
  \mathfrak{s},\text{ and } i(x_{\alpha})=\iota (e_{\alpha_{1}})+\ldots +\iota
  (e_{\alpha_{l}}),\, \alpha \in \Phi,
\end{equation}
define a linear injection $i$ from $\mathfrak{g}$ to $V_{L}$. For
$\alpha \in \mathfrak{h}$, $n\in \mathbf{Z}$, we denote by $\alpha (n)$
the operator on $V_{L}$ determined by $\alpha \otimes t^{n}$, and we
let $\alpha (z)=\sum_{n\in \mathbf{Z}}\alpha (n)z^{-n-1}$.

We can now introduce the relative untwisted vertex operators of
\cite{DL1}. These operators are parametrized by elements of the space
$V_{L}$ on which they also act, and they are defined 
  relative to $\mathfrak{h}_{*}$. We use a ``normal ordering'' procedure
$_{\text{{\tiny $\circ$}}}^{\text{{\tiny
      $\circ$}}}\cdot\phantom{}_{\text{{\tiny $\circ$}}}^{\text{{\tiny
      $\circ$}}}$ to signify that the enclosed expression is to be
reordered if necessary so that all the operators $\alpha (n)$ ($\alpha
\in \mathfrak{h}$, $n<0$) and $a\in \hat{L}$ are to be placed to the
left of all the operators $\alpha (n)$ and $z^{\alpha}$ ($\alpha
\in \mathfrak{h}$, $n\ge 0$) before the expression is evaluated. Using
an obvious formal integration notation, we set
\begin{equation}\label{4.2.8}
  Y_{*}^{(l)}(\iota (a),z)=Y_{*}^{(l)}(a,z)=_{\text{{\tiny $\circ$}}}^{
\text{{\tiny
      $\circ$}}}\!\exp\bigg[\int\!\!\big(\bar{a}'(z)-\bar{a}'(0)z^{-1}\big)
\bigg]\,a\,z^{\bar{a}'}\phantom{}_{\text{{\tiny $\circ$}}}^{\text{{\tiny
      $\circ$}}}, \, a\in \hat{L}.
\end{equation}
If $v\in V_{L}$ is of the form
  $v=\alpha_{1}(-n_{1})\cdots \alpha_{j}(-n_{j})\cdot \iota (a)$, 
 $a\in \hat{L}$, $\alpha_{k}\in \mathfrak{h}$,
$n_{k}\in \mathbf{Z}_{+}$, $1\le k\le j$, we define
\begin{equation}\label{4.2.9}
  Y_{*}^{(l)}(v,z)=_{\text{{\tiny $\circ$}}}^{\text{{\tiny
      $\circ$}}}\!\!\left[\frac{1}{(n_{1}-1)!}\left(\!\frac{d}{dz}\!\right)^
{n_{1}-1}\!\alpha_{1}'(z)\right]\cdots \left[\frac{1}{(n_{j}-1)!}\left(\!
\frac{d}{dz}\!\right)^{n_{j}-1}\!\alpha_{j}'(z)\right]Y_{*}^{(l)}(a,z)\,_{
\text{{\tiny $\circ$}}}^{\text{{\tiny
      $\circ$}}},
\end{equation}
{\allowdisplaybreaks
and we extend this definition to $V_{L}$ by linearity. One gets a
well-defined linear map
\begin{align}
    V_{L} & \longrightarrow (\text{End}\,V_{L})\{z\}\label{4.2.10}\\
    v  & \longmapsto Y_{*}^{(l)}(v,z)=\sum_{n\in
      \mathbf{C}}v_{n}z^{-n-1},\, v_{n}\in \text{End}\,V_{L},\nonumber
\end{align}
where for any vector space $W$, $W\{z\}$ denotes the linear space of
$W$-valued formal series in $z$.} The case $\mathfrak{h}_{*}=0$ leads to 
the (ordinary)
untwisted vertex operators as defined in \cite[Ch.~8]{FLM}. In this
case, the operators $Y_{*}^{(l)}(v,z)$ are denoted just by
$Y^{(l)}(v,z)$. We
refer to \cite{DL1} and \cite{DL2} for a detailed discussion of the properties
satisfied by the operators $Y_{*}^{(l)}(v,z)$. One has in particular
\begin{equation}\label{4.2.11}
  [\tilde{\mathfrak{h}}_{*}',Y_{*}^{(l)}(v,z)]=0 \text{ for } v\in
  V_{L},
\end{equation}
{\allowdisplaybreaks
where
$\tilde{\mathfrak{h}}_{*}'=[\tilde{\mathfrak{h}}_{*},\tilde{\mathfrak{h}}_{*}]$
is the Heisenberg algebra associated with the abelian Lie algebra
$\mathfrak{h}_{*}$ (see \eqref{4.2.4}). 

Notice that when applied to a
single copy of $Q$ in the case $\mathfrak{h}_{*}=0$, the above
construction yields a well-defined linear map
\begin{align*}
    V_{Q} & \longrightarrow (\text{End}\,V_{Q})[[z,z^{-1}]]\\
    v  & \longmapsto Y^{(1)}(v,z)=\sum_{n\in
      \mathbf{Z}}v_{n}z^{-n-1},\nonumber
\end{align*}
the lattice $Q$ being even. Furthermore, using the notation}
\begin{equation}\label{4.2.12}
  \mathfrak{s}_{i}=Q_{i}\otimes_{\mathbf{Z}}\mathbf{C},\,
  \hat{\mathfrak{s}}_{i}=\mathfrak{s}_{i}\otimes
  \mathbf{C}[t,t^{-1}]\oplus \mathbf{C}c,\,
  \tilde{\mathfrak{s}}_{i_{\pm}}'=\mathfrak{s}_{i}\otimes
  t^{\pm1}\mathbf{C}[[t^{\pm1}]],\, 1\le i\le l,
\end{equation}
one gets from \eqref{4.1.2} that
  $\mathbf{C}\{L\}\cong \otimes_{i=1}^{l}\mathbf{C}\{Q_{i}\}$ 
    and $S(\tilde{\mathfrak{h}}_{-}')\cong
  \otimes_{i=1}^{l}S(\tilde{\mathfrak{s}}_{i_{-}}')$ 
as vector spaces, and thus
  $V_{L}\cong \otimes_{i=1}^{l}V_{Q_{i}}$ linearly 
(cf.~\eqref{4.2.6}). For $1\le j\le l$ define a linear injection $i_{j}$ 
by 
{\allowdisplaybreaks
\begin{align}
  i_{j}:\, \mathfrak{g} & \longrightarrow V_{L} \label{4.2.13}\\
  h & \longmapsto h_{j}(-1)\cdot \iota (1),\, h\in
  \mathfrak{s},\notag \\
  x_{\alpha} & \longmapsto \iota (e_{\alpha_{j}}),\, \alpha \in
  \Phi,\notag 
\end{align}
so that $i(x)=\sum_{j=1}^{l}i_{j}(x)$ for $x\in
\mathfrak{g}$ by \eqref{4.2.7}.} Then from \cite[Ch.~13]{DL1} and 
\cite[\S 7.2]{FLM} one gets:
{\allowdisplaybreaks 
\begin{theorem}
  The linear map $\rho : \hat{\mathfrak{g}}\rightarrow 
  \text{{\em End}}\,V_{Q}$ given by
\begin{equation}\label{4.2.14}
  \rho (c)=1,\quad \rho (x(z))=Y^{(1)}(i(x),z),\, x\in
\mathfrak{g},
\end{equation}
defines a level 1 $\hat{\mathfrak{g}}$-module structure on $V_{Q}$
such that $V_{Q}$ is thereby isomorphic to the basic module $L(\Lambda_{0}
;\hat{\mathfrak{g}})$. Moreover, the linear map $\pi :
\hat{\mathfrak{g}}\rightarrow 
  \text{{\em End}}\,V_{L}$ determined by
\begin{equation}\label{4.2.15}
  \pi (c)=l,\quad \pi (x(z))=Y^{(l)}(i(x),z),\, x\in
\mathfrak{g},
\end{equation}
defines a level $l$ $\hat{\mathfrak{g}}$-module structure on $V_{L}$
such that the representation \eqref{4.2.15} is isomorphic to the tensor
product of the representations
\begin{equation}\label{4.2.16}
  \pi_{j} : \hat{\mathfrak{g}}\rightarrow 
  \text{{\em End}}\,V_{Q_{j}},\, 1\le j\le l,
\end{equation}
given by
  $\pi_{j}(c)=1,\, \pi_{j}(x(z))=Y^{(1)}(i_{j}(x),z),\, x\in \mathfrak{g}$.
\end{theorem}}
\begin{remark}
  (i) The
representations \eqref{4.2.14}-\eqref{4.2.16} extend to 
$\tilde{\mathfrak{g}}$ by
letting $d$ act as the degree operator (cf.~\eqref{4.2.2} and 
\eqref{4.2.5} for
$\mathfrak{h}_{*}=0$), in which case \eqref{4.2.14} is the well-known
Frenkel-Kac-Segal construction of the affine Lie algebra
$\tilde{\mathfrak{g}}$. \\
\noindent
(ii) Under the isomorphism in the second part of
Theorem 4.2.1, the action of  
  $Y^{(l)}(i_{j}(x),z)$ on $V_{L}$ coincides with the action of
  $1\otimes \ldots \otimes 1\otimes Y^{(1)}(i_{j}(x),z)\otimes 1\ldots
  \otimes 1$, where the nontrivial factor is in the j-th
  position.
\end{remark}
By \cite[Ch.~13]{DL1}, $V_{L}$ is a completely
reducible $\hat{\mathfrak{g}}$-module and its 
$\hat{\mathfrak{g}}$-submodule $U(\hat{\mathfrak{g}})\cdot \iota (1)$
is isomorphic to the standard $\hat{\mathfrak{g}}$-module
$L(l\Lambda_{0};\hat{\mathfrak{g}})$, so that we may view  
$L(l\Lambda_{0};\hat{\mathfrak{g}})$ as a subspace of $V_{L}$.
\begin{corollary}
On the subspace $L(l\Lambda_{0};\hat{\mathfrak{g}})$ of $V_{L}$ the map
$Y^{(l)}(\cdot \,,z)$ of \eqref{4.2.15} coincides with the map 
$Y(\cdot \,,z)$ of Theorem 3.1.4 and Theorem 3.2.5.
\end{corollary}
\begin{proof}
Let $x\in \mathfrak{g}$, $n\in \mathbf{Z}$, and denote by $x(n)$ the operator
$\pi (x\otimes t^{n})$ acting on $V_{L}$. Then the map $i$ can be
rewritten as $i(x)=x(-1)\cdot \iota (1)$, and by comparing the
expression for $x(z)$ in \eqref{4.1.11} with the generic expansion
$Y^{(l)}(i(x),z)=\sum_{n\in \mathbf{Z}}x_{n}z^{-n-1}$ (cf.~\eqref{4.2.10})
we see that in fact $x_{n}=x(n)$. For $x\in
\mathfrak{g}$ one has therefore that
$Y^{(l)}(i(x),z)=Y(x(-1)\mathbf{1},z)$, where $\mathbf{1}=1\otimes
\iota (1)$ and $Y$ is the map defined in \eqref{3.1.4}. Furthermore,
from $L(l\Lambda_{0};\hat{\mathfrak{g}})\cong U(\hat{\mathfrak{g}})\cdot
\iota (1)$ one gets that $L(l\Lambda_{0};\hat{\mathfrak{g}})$ is
generated by the set $\{i(x)\mid x\in
\mathfrak{g}\}\cup \{\mathbf{1}\}$ through the map $Y^{(l)}(\cdot
\,,z)$, which by construction satisfies the same iterate formula 
as the map $Y(\cdot \,,z)$ (cf.~\eqref{3.1.2}). This proves the corollary.
\end{proof}

\subsection{Relative twisted vertex operators}

We now concentrate on the twisted counterpart of the above
constructions. Let again $\mathfrak{h}_{*}\subset \mathfrak{h}$ be as
in \eqref{4.1.13}-\eqref{4.1.15} and recall from \eqref{4.1.6}-\eqref{4.1.8} 
the central extension $\hat{L}_{\sigma}$ of $L$ and 
the set-theoretic identification
between the groups $\hat{L}$ and $\hat{L}_{\sigma}$.
 Notice that $1-\sigma=\sigma ^{-1}$. It follows from 
\cite[Proposition 6.1]{L2} (see
also \cite{Ca}) that there exists a unique homomorphism $\psi :
\hat{L}_{\sigma}\rightarrow \mathbf{C}^{\times}$ such that $\psi
(\eta)=\eta$ and $L$ maps to the identity. Let $T=\mathbf{C}_{\psi }$
be the one-dimensional $\hat{L}_{\sigma}$-module affording $\psi $,
and give $T$ the trivial $\mathbf{C}$-gradation. We let $\hat{\sigma}$
act on $T$ as a grading-preserving linear automorphism 
(cf.~\eqref{4.1.8}-\eqref{4.1.9}). By Remark 3.2.3, the 
$\sigma$-decomposition of $\mathfrak{h}$ is given by 
  $\mathfrak{h}=\coprod_{n\in
    \mathbf{Z}/6\mathbf{Z}}\mathfrak{h}_{(n)}=\mathfrak{h}_{(1)}\oplus
  \mathfrak{h}_{(5)}$. Let $\mathfrak{h}_{(n)}=\mathfrak{h}_{(n\bmod 6)}$, 
$n\in \mathbf{Z}$, and form the $\sigma$-twisted affine Lie algebra 
$\tilde{\mathfrak{h}}[\sigma]=\hat{\mathfrak{h}}[\sigma]\rtimes
  \mathbf{C}d$ and its subalgebras
\begin{equation}\label{4.3.1}
  \hat{\mathfrak{h}}[\sigma]=\oplus_{n\in
    \tfrac{1}{6}\mathbf{Z}}\mathfrak{h}_{(6n)}\otimes t^{n}\oplus
  \mathbf{C}c,\quad \tilde{\mathfrak{h}}[\sigma]'_{\pm}=\oplus_{n\in
    \tfrac{1}{6}\mathbf{Z},\, \pm n>0}\!\mathfrak{h}_{(6n)}\otimes
  t^{n},
\end{equation}
so that $\hat{\mathfrak{h}}[\sigma]=\tilde{\mathfrak{h}}[\sigma]'_{+}\oplus
\tilde{\mathfrak{h}}[\sigma]'_{-}\oplus \mathbf{C}c$. 
Then 
$\tilde{\mathfrak{h}}[\sigma]':=\big[\tilde{\mathfrak{h}}[\sigma],
\tilde{\mathfrak{h}}[\sigma]\big]$ is a
Heisenberg Lie algebra and $\tilde{\mathfrak{h}}[\sigma]'=
\hat{\mathfrak{h}}[\sigma]$. On $\hat{\mathfrak{h}}[\sigma]$ we
define the weight gradation associated with $\mathfrak{h}_{*}$
by $\mbox{wt}(x\otimes t^{m})=0$, $\mbox{wt}(y\otimes t^{n})=-n$, 
$\mbox{wt}(c)=0$, $m,n\in \tfrac{1}{6}\mathbf{Z}$, 
$x\in {\mathfrak{h}_{*}}_{(6m)}$,
$y\in {\mathfrak{h}_{*}^{\perp}}_{\!(6n)}$ (cf.~\eqref{4.2.5}). Form the
induced irreducible $\hat{\mathfrak{h}}[\sigma]$-module
  $M(\sigma;1,\mathfrak{h})=U\big(\hat{\mathfrak{h}}[\sigma]\big)\otimes_{
U(\oplus_{n\ge 0}\mathfrak{h}_{(6n)}\oplus \mathbf{C}c)}\mathbf{C}\cong S\big(
\tilde{\mathfrak{h}}[\sigma]'_{-}\big)$ (linearly), 
where $\oplus_{n\ge 0}\mathfrak{h}_{(6n)}$ acts trivially on
$\mathbf{C}$ and $c$ acts as 1. Then $M(\sigma;1,\mathfrak{h})$ has a
natural $\tfrac{1}{6}\mathbf{Z}$-grading which is compatible with the
action of $\hat{\mathfrak{h}}[\sigma]$ and such that $\mbox{wt}(1)=0$. 
Moreover, $\sigma$ acts in a grading-preserving way on
$\hat{\mathfrak{h}}[\sigma]$ (fixing $c$) and on
$M(\sigma;1,\mathfrak{h})$ (as an algebra isomorphism). Let 
$\hat{\mathfrak{h}}_{*}[\sigma]_{-}$ and
$\hat{\mathfrak{h}}_{*}^{\perp}[\sigma]_{-}$ be defined as in
\eqref{4.3.1} and set
\begin{equation}\label{4.3.2}
  V_{L}^{T}=M(\sigma;1,\mathfrak{h})\otimes_{\mathbf{C}}T\cong
  S\big(\tilde{\mathfrak{h}}[\sigma]'_{-}\big)\cong
  S\big(\hat{\mathfrak{h}}_{*}[\sigma]_{-}\big)\otimes
  S\big(\hat{\mathfrak{h}}_{*}^{\perp}[\sigma]_{-}\big)\, \text{
    (linearly)},
\end{equation}
(cf.~\eqref{4.1.13}). Using the gradations of $M(\sigma;1,\mathfrak{h})$ 
and $T$,
we see that $V_{L}^{T}$ is naturally $\mathbf{Q}$-graded. Clearly,
$\hat{\mathfrak{h}}[\sigma]$ and $\hat{L}_{\sigma}$ act on $V_{L}^{T}$
by acting either on $M(\sigma;1,\mathfrak{h})$ or $T$, and
$\hat{\sigma}$ extends to a linear automorphism of $V_{L}^{T}$ such that
$\hat{\sigma}(u\otimes t)=\sigma (u)\otimes \hat{\sigma}(t)=t\,\sigma
(u)\otimes 1$ for $u\in M(\sigma;1,\mathfrak{h})$, $t\in T$.

For $\alpha \in \mathfrak{h}$ let $\alpha_{(6n)}(n)$ be the operator on 
$M(\sigma;1,\mathfrak{h})$
corresponding to $\alpha_{(6n)}\otimes t^{n}$, $n\in
    \tfrac{1}{6}\mathbf{Z}$. Set
\begin{equation}\label{4.3.3}
  \alpha (\sigma ;z)=\sum_{n\in
    \tfrac{1}{6}\mathbf{Z}}\alpha_{(6n)}(n)\,z^{-n-1}.
\end{equation}
We shall also use the scalar function
\begin{equation}\label{4.3.4}
  \tau (\alpha)=2^{-\langle
    \alpha, \alpha\rangle}(1-\eta^{-1})^{\langle \sigma
    \alpha, \alpha \rangle}(1-\eta^{-2})^{\langle \sigma^{2}
    \alpha, \alpha \rangle}\,\text{ for } \alpha \in L,
\end{equation}
which obviously satisfies $\tau (\sigma \alpha)=\tau (\alpha)$. Given
$a\in \hat{L}$, one defines the relative $\sigma$-twisted vertex operator 
$Y_{*}^{\sigma,(l)}(\iota (a),z)$ acting on
$V_{L}^{T}$ as follows: for $\alpha \in \mathfrak{h}$ introduce first the
operators
\begin{equation}\label{4.3.5}
  E_{*(1)}^{\pm}(\alpha,z)=\exp\left[\sum_{n\in
    \tfrac{1}{6}\mathbf{Z},\, \pm
    n>0}\!\!\frac{\alpha'_{(6n)}(n)}{n}\,z^{-n}\right]\in
  \big(\mbox{End}\,S\big(\hat{\mathfrak{h}}_{*}^{\perp}[\sigma]_{-}\big)\big)
\big[\big[z^{\pm 1/6}\big]\big]
\end{equation}
(which by \eqref{4.3.2} is a subspace of 
$\big(\mbox{End}\,V_{L}^{T}\big)\{z\}$), and then set
{\allowdisplaybreaks
  \begin{align}
        Y_{*}^{\sigma,(l)}(\iota
    (a),z)&=Y_{*}^{\sigma,(l)}(a,z)=6^{-\langle
        \bar{a}', \bar{a}'\rangle/2}\,\tau (\bar{a}')\,_{
\text{{\tiny $\circ$}}}^{\text{{\tiny
      $\circ$}}}\exp\bigg[\int \bar{a}'(\sigma;z)\bigg]\,a\,z^{-\langle
        \bar{a}', \bar{a}'\rangle/2}\,\phantom{}_{\text{{\tiny $\circ$}}}^{
\text{{\tiny
      $\circ$}}}\label{4.3.6} \\
&=6^{-\langle
        \bar{a}', \bar{a}'\rangle/2}\,\tau (\bar{a}')\,E_{*(1)}^{-}
(-\bar{a},z)E_{*(1)}^{+}(-\bar{a},z)\,a\,z^{-\langle
        \bar{a}', \bar{a}'\rangle/2},\notag
  \end{align}
where $a$ is viewed as an element of $\hat{L}_{\sigma}$ by means of
the set-theoretic identification between $\hat{L}$ and
$\hat{L}_{\sigma}$.} For $\alpha_{i} \in
\mathfrak{h}$, $n_{i}\in \mathbf{Z}_{+}$, $1\le i\le j$, and 
$v=\alpha_{1}(-n_{1})\cdots\alpha_{j}(-n_{j})\cdot \iota (a)\in
V_{L}$, set
\begin{equation}\label{4.3.7}
  W_{*}(v,z)=_{\text{{\tiny $\circ$}}}^{\text{{\tiny
      $\circ$}}}\!\!\left[\frac{1}{(n_{1}-1)!}\left(\!\frac{d}{dz}\!\right)^{
n_{1}-1}\!\alpha_{1}'(\sigma ;z)\right]\cdots \left[\frac{1}{(n_{j}-1)!}\left(
\!\frac{d}{dz}\!\right)^{n_{j}-1}\!\alpha_{j}'(\sigma ;z)\right]Y_{*}^{\sigma,
(l)}(a,z)\,_{\text{{\tiny $\circ$}}}^{\text{{\tiny
      $\circ$}}}
\end{equation}
(the right-hand side being an operator on $V_{L}^{T}$), and then
extend this definition to all $v\in V_{L}$ by linearity. Let
$\{\gamma_{1},\ldots ,\gamma_{d}\}$ be an orthonormal basis of
$\mathfrak{h}_{*}^{\perp}$ and define the constants $c_{_{mni}}\in 
\mathbf{C}$, $m,n\in \mathbf{N}$, $i\in \{0,1,\ldots,5\}$, as in 
\cite[(4.41)]{DL2}. In particular, $c_{_{00i}}=0$ for all $i$. The 
following operator is independent of the choice of 
orthonormal basis in $\mathfrak{h}_{*}^{\perp}$ (cf.~\cite[(4.42)]{DL2}, 
\cite[\S 9.2]{FLM}):
\begin{equation}\label{4.3.8}
  \Delta_{z*}=\sum_{m,n\ge
  0}\sum_{i=0}^{5}\sum_{j=1}^{d}c_{_{mni}}(\sigma^{-i}\gamma_{j})(m)
\gamma_{j}(n)
 z^{-m-n}\in (\mbox{End}\,V_{L})[[z^{-1}]].
\end{equation} 
Since $c_{_{00i}}=0$ for all $i$, $\exp (\Delta_{z*})$ is
well defined on $V_{L}$ and $\exp (\Delta_{z*})v\in V_{L}[z^{-1}]$ for
$v\in V_{L}$. The 
  relative $\sigma$-twisted vertex operator $Y_{*}^{\sigma,(l)}(v,z)$
is then defined by
\begin{equation}\label{4.3.9}
  Y_{*}^{\sigma,(l)}(v,z)=W_{*}(\exp (\Delta_{z*})v,z)\,\text{ for } \,
  v\in V_{L}.
\end{equation}

Summarizing, we get a well-defined linear map (the relative 
$\sigma$-twisted vertex operator map)
{\allowdisplaybreaks
 \begin{align}
  V_{L} & \longrightarrow \big(\mbox{End}\,V_{L}^{T}\big)\{z\}\label{4.3.10} \\
    v & \longmapsto Y_{*}^{\sigma,(l)}(v,z)=\sum_{n\in
      \mathbf{C}}v_{n}z^{-n-1},\, v_{n}\in
    \mbox{End}\,V_{L}^{T},\nonumber
  \end{align}
which satisfies $Y_{*}^{\sigma,(l)}(v,z)\in
\big(\mbox{End}\,V_{L}^{T}\big)\big[\big[z^{1/6},z^{-1/6}\big]\big]$
if $\mathfrak{h}_{*}=0$.} We refer to \cite{DL2} for a
discussion of the properties of these 
operators. Of particular importance for our purposes is the
fact that
\begin{equation}\label{4.3.11}
  \big[\hat{\mathfrak{h}}_{*}[\sigma],Y_{*}^{\sigma,(l)}(v,z)\big]=0
  \text{ for } v\in V_{L}
\end{equation}
(cf.~\cite[Proposition 4.6]{DL2}), where
$\hat{\mathfrak{h}}_{*}[\sigma]$ is defined as in
\eqref{4.3.1}-\eqref{4.3.2}.

By using a single copy of $Q$ and the space $\mathfrak{h}_{*}=0$ we get a 
well-defined linear map
{\allowdisplaybreaks
  \begin{align*}
    V_{Q} & \longrightarrow \big(\mbox{End}\,V_{Q}^{T}\big)\big[\big[z^{1/6},
z^{-1/6}\big]\big] \\
    v & \longmapsto Y^{\sigma,(1)}(v,z)=\sum_{n\in
      \tfrac{1}{6}\mathbf{Z}}v_{n}z^{-n-1}.
  \end{align*}
Let $\hat{\mathfrak{s}}_{i}[\sigma]=\oplus_{n\in
      \tfrac{1}{6}\mathbf{Z}}\mathfrak{s}_{i_{(6n)}}\otimes
    t^{n}\oplus \mathbf{C}c$ and $\hat{\mathfrak{s}}_{i}[\sigma]_{\pm}=
\oplus_{n\in
      \tfrac{1}{6}\mathbf{Z},\, \pm n>0}\!\mathfrak{s}_{i_{(6n)}}\otimes
    t^{n}$, $1\le i\le l$. Since  
  $S\big(\hat{\mathfrak{h}}[\sigma]_{-}\big)=S\big(\tilde{\mathfrak{h}}
[\sigma]_{-}'\big)=\otimes_{i=1}^{l}S\big(\hat{\mathfrak{s}}_{i}[\sigma]_{-}
\big)$, \eqref{4.3.2} implies that}
\begin{equation}\label{4.3.12}
  V_{L}^{T}\cong \otimes_{i=1}^{l}V_{Q_{i}}^{T} \,\text{ (linearly)}.
\end{equation}
Then from \cite[\S 9]{L2} and \cite[\S 8-9]{LW} (see also \cite[\S
7.4]{FLM}, \cite{Fig}) one
gets the following 
{\allowdisplaybreaks
\begin{theorem}
  The linear map $\rho^{\sigma}:\hat{\mathfrak{g}}[\sigma]\rightarrow
  \text{{\em End}}\,V_{Q}^{T}$ given by
 \begin{equation}\label{4.3.13}
   \rho^{\sigma}(c)=1,\quad
   \rho^{\sigma}(x(\sigma;z))=Y^{\sigma,(1)}(i(x),z),\, x\in
   \mathfrak{g},
 \end{equation}
 defines a level 1 $\hat{\mathfrak{g}}[\sigma]$-module structure on
 $V_{Q}^{T}$ such that $V_{Q}^{T}\cong
 L(\Lambda_{1};\hat{\mathfrak{g}}[\sigma])$. Moreover, the linear map
 $\pi^{\sigma}: \hat{\mathfrak{g}}[\sigma]\rightarrow
  \text{{\em End}}\,V_{L}^{T}$ determined by
  \begin{equation}\label{4.3.14}
    \pi^{\sigma}(c)=l,\quad \pi^{\sigma}(x(\sigma;z))=Y^{\sigma,(l)}(i(x),z),
\, x\in
   \mathfrak{g},
 \end{equation}
 defines a level $l$ $\hat{\mathfrak{g}}[\sigma]$-module structure on
 $V_{L}^{T}$, and this representation is isomorphic to the tensor
 product of the representations
 \begin{equation}\label{4.3.15}
   \pi^{\sigma}_{j}: \hat{\mathfrak{g}}[\sigma]\longrightarrow
  \text{{\em End}}\,V_{Q_{j}}^{T},\quad 1\le j\le l,
 \end{equation}
 given by
   $\pi^{\sigma}_{j}(c)=1,\, \pi^{\sigma}_{j}(x(\sigma;z))=Y^{\sigma,(1)}
(i_{j}(x),z),\, x\in \mathfrak{g}$.
\end{theorem}}
\begin{remark}
  Under the isomorphism in the second part of Theorem 4.3.1 the action of the 
operator 
$Y^{\sigma,(l)}(i_{j}(x),z)$ on
  $V_{L}^{T}$ coincides with the action of $1\otimes \ldots \otimes
  1\otimes Y^{\sigma,(1)}(i_{j}(x),z)\otimes \ldots \otimes
  1$. Moreover, the representations 
   \eqref{4.3.13}-\eqref{4.3.15} extend to $\tilde{\mathfrak{g}}[\sigma]$ by
   letting $d$ act as the degree operator.
\end{remark}
\begin{corollary}
On the subspace $L(l\Lambda_{0};\hat{\mathfrak{g}})$ of $V_{L}$ the map
   $Y^{\sigma,(l)}(\cdot \,,z)$ of \eqref{4.3.14} coincides with
   the map $Y^{\sigma}(\cdot
   \,,z)$ of Theorem 3.1.5 and Theorem 3.2.5.
\end{corollary}
\begin{proof}
By \cite[Theorem
   7.1]{DL2}, $\big(V_{L},Y^{(l)},\mathbf{1},\omega
   \big)$ is a VOA and $\big(V_{L}^{T},Y^{\sigma,(l)}\big)$ is a
   $\sigma$-twisted $V_{L}$-module (see \cite{DL1} and \cite{DL2} for the
   construction of the Virasoro 
   element $\omega \in V_{L}$). In particular, the map 
   $Y^{\sigma,(l)}$ satisfies the twisted associator formula
   \eqref{3.1.2} with $Y_{M}=Y^{\sigma,(l)}$ and $Y=Y^{(l)}$. Let
   $M$ be any standard level $l$ 
   $\hat{\mathfrak{g}}[\sigma]$-module. According to \S 2.1,
   \eqref{4.3.12}, and Theorem 4.3.1, $M$ can be isomorphically
   embedded into $L\big(\Lambda_{1}
;\hat{\mathfrak{g}}[\sigma]\big)^{\otimes \,l}\cong
V_{L}^{T}$. Denoting
  by $x_{(6n)}(n)$ the operator $\pi^{\sigma}(x_{(6n)}\otimes t^{n})$
  on $V_{L}^{T}$, $x\in
   \mathfrak{g}$, $n\in \tfrac{1}{6}\mathbf{Z}$, and comparing the
   expression for $x(\sigma;z)$ from \eqref{4.1.11} with the generic
   expansion $Y^{\sigma,(l)}(i(x),z)=\sum_{n\in
     \tfrac{1}{6}\mathbf{Z}}x_{n}z^{-n-1}$ (cf.~\eqref{4.3.10}), we see
   that $x_{n}=x_{(6n)}(n)$. For $x\in
   \mathfrak{g}$ one has therefore that 
   $Y^{\sigma,(l)}(i(x),z)=Y^{\sigma}(x(-1)\mathbf{1},z)$ when acting on $M$,
   where (as before) $\mathbf{1}=1\otimes \iota (1)\in
   L(l\Lambda_{0};\hat{\mathfrak{g}})$ and $Y^{\sigma}(\cdot
   \,,z)$ is the map defined in \eqref{3.1.5}. Then the corollary follows from 
   the twisted iterate formula \eqref{3.1.2} and Corollary 4.2.3.
\end{proof}

\subsection{The GVOA $\Omega_{2k+1}^{A}$ and its action on 
the vacuum spaces $\Omega_{k_{0},2k+1}$}

Set now 
{\allowdisplaybreaks
\begin{equation}\label{4.4.1}
  \begin{split}
  \Omega_{*} & =\{v\in V_{L} \mid h(n)v=0 \text{ for } h\in
  \mathfrak{h}_{*},\, n>0 \}\subset V_{L},\\
  \Omega_{*}^{\sigma} & =\{v\in V_{L}^{T} \mid h(n)v=0 \text{ for } h\in
  \mathfrak{h}_{*},\, n\in \tfrac{1}{6}\mathbf{Z}_{+} \}\subset
  V_{L}^{T}.
  \end{split}
\end{equation}
Clearly, $\Omega_{*}$ is the vacuum space for the action of the
Heisenberg algebra $\tilde{\mathfrak{h}}_{*}'$ on $V_{L}$, while
$\Omega_{*}^{\sigma}$ is the vacuum space for the action of the
Heisenberg algebra $\hat{\mathfrak{h}}_{*}[\sigma]$ on $V_{L}^{T}$.} By
\eqref{4.2.11} and \eqref{4.3.11}, these spaces are preserved by 
the operators $Y_{*}^{(l)}(v,z)$, $v\in V_{L}$, and $Y_{*}^{\sigma
  ,(l)}(v,z)$, $v\in V_{L}$, respectively. Since
$S\big(\tilde{\mathfrak{h}}_{-}'\big)=S\big((\tilde{\mathfrak{h}}_{*})'_{-}
\big)\otimes
S\big((\tilde{\mathfrak{h}}^{\perp}_{*})'_{-}\big)$, one gets from 
\eqref{4.2.6} that
\begin{equation}\label{4.4.2}
  \Omega_{*}=S\big((\tilde{\mathfrak{h}}^{\perp}_{*})'_{-}\big)\otimes
  \mathbf{C}\{L\},
\end{equation}
and thus $V_{L}$ may be decomposed as
{\allowdisplaybreaks
\begin{equation}\label{4.4.3}
  V_{L}=\Omega_{*}\oplus V_{*},\,
  V_{*}=\text{$\mathbf{C}$-span}\,\{h(n)V_{L} \mid h\in
  \mathfrak{h}_{*},
  n<0\}=(\tilde{\mathfrak{h}}_{*})'_{-}S\big((\tilde{\mathfrak{h}}_{*})'_{-}
\big)\otimes \Omega_{*}.
\end{equation}
Note that by construction the
operators $Y_{*}^{\sigma,(l)}$ satisfy
\begin{equation}\label{4.4.4}
  Y_{*}^{\sigma,(l)}(v,z)=0 \, \text{ if }\, v\in V_{*}.
\end{equation}
In similar fashion, \eqref{4.3.2} yields
  $\Omega_{*}^{\sigma}=S(\hat{\mathfrak{h}}_{*}^{\perp}[\sigma]_{-})\otimes T$.

Recall \eqref{4.1.2} and \eqref{4.1.12} and consider the 
diagonal embedding}
{\allowdisplaybreaks
  \begin{align}
    _{*}\,:\,\mathfrak{s} & \longrightarrow \mathfrak{h}\label{4.4.5}\\
    \alpha & \longmapsto \alpha_{*}=\alpha_{1}+\ldots
    +\alpha_{l},\nonumber
  \end{align}
where ${\alpha}_{i}\in Q_{i}\otimes_{\mathbf{Z}}\mathbf{C}$ corresponds to 
$\alpha \in \mathfrak{s}$. Obviously, the map $_{*}$ defines a group 
isomorphism between $Q$ and $Q_{*}$.} We assume henceforth that the subspace
$\mathfrak{h}_{*}$ of $\mathfrak{h}$ (see \eqref{4.1.13}) is
\begin{equation}\label{4.4.6}
  \mathfrak{h}_{*}=\{\alpha_{*} \mid \alpha \in
\mathfrak{s}\}=\mathfrak{s}_{*}.
\end{equation}
Then $\mathfrak{h}_{*}^{\perp}=\text{$\mathbf{C}$-span}\,\{\alpha_{i}-
\alpha_{j} \mid \alpha \in \mathfrak{s},\, i,j=1,\ldots l \}$ 
(cf.~\eqref{4.1.13}). Define
\begin{equation}\label{4.4.7}
  A=\{e_{\alpha_{*}} \mid \alpha \in 2Q \}\subset \hat{L}.
\end{equation}
By the choice of the cocycle $\epsilon_{0}$ 
and of the section $e$ (cf.~\eqref{4.1.4}-\eqref{4.1.5} and 
\eqref{4.1.10}), $A$ is a central
subgroup of $\hat{L}$ isomorphic to the subgroup $2Q_{*}$ of $Q_{*}$
and such that $A\cap \langle \eta \rangle =1$ ($=e_{0}$). Moreover, 
$\langle L\,, \bar{A} \rangle \in 2\mathbf{Z}$, $\mathfrak{h}_{*}\cap
L$ $\mathbf{C}$-spans $\mathfrak{h}_{*}$, and $\text{rank}\, A = 
\text{rank}\,\mathfrak{g}=\dim
  \mathfrak{s}=2$, 
so that $\bar{A}$ $\mathbf{C}$-spans $\mathfrak{h}_{*}$. Recall the
projections \eqref{4.1.15} and notice that
\begin{equation}\label{4.4.8}
  \alpha_{i}''=l^{-1}\alpha_{*} \, \text{ for }\, \alpha \in
  \Phi,\, i=1,\ldots ,l,
\end{equation}
so that
\begin{equation}\label{4.4.9}
  L''=\{ \beta ''\mid \beta \in L \}=l^{-1}Q_{*}.
\end{equation}

Let $G=L''/\bar{A}=l^{-1}Q_{*}/2Q_{*}\cong
 l^{-1}Q/2Q$  
and notice that the smallest positive integer $S$ such that $\langle
\alpha '\,,\beta ' \rangle \in S^{-1}\mathbf{Z}$, $\langle
\alpha \,,\beta \rangle \in S^{-1}\mathbf{Z}$, and $\langle
\alpha '\,,\alpha ' \rangle \in 2S^{-1}\mathbf{Z}$ for $\alpha$,
$\beta \in L$ is precisely 
\begin{equation}\label{4.4.10}
  S=l.
\end{equation}
A well-defined symmetric nondegenerate
$l^{-1}\mathbf{Z}/2\mathbf{Z}$-valued
$\mathbf{Z}$-bilinear form on $G$ is then given by
{\allowdisplaybreaks
\begin{align*}
  G\times G & \longrightarrow
  l^{-1}\mathbf{Z}/2\mathbf{Z}\\
  \big(l^{-1}\alpha_{*}+\bar{A},l^{-1}\beta_{*}+\bar{A}\big)
  & \longmapsto l^{-1}\langle \alpha\,,\beta \rangle
    +2\mathbf{Z} \, \text{ for }\, \alpha ,\beta \in Q.
\end{align*}

{\em In the remainder of this section we shall denote the
standard $\hat{\mathfrak{g}}$-module
$L\big(l\Lambda_{0};\hat{\mathfrak{g}}\big)$ simply by $L(l,0)$}. Recall that 
$L(l,0)\cong U(\hat{\mathfrak{g}})\cdot \iota (1)\subset V_{L}$ 
(cf.~Corollary 4.2.3).} The action of $\hat{\mathfrak{g}}$ on
$V_{L}$ being given by Theorem 4.2.1, one gets from 
\cite[Proposition 14.5]{DL1} and the complete reducibility of $V_{L}$ as a 
$\hat{\mathfrak{g}}$-module that any $\hat{\mathfrak{g}}$-submodule $M$ of 
$V_{L}$ is $A$-stable and has an $\hat{\mathfrak{h}}_{*}$-stable complement in
$V_{L}$. In particular, this holds for $M=L(l,0)$. Define the following spaces:
{\allowdisplaybreaks
\begin{align}
  & \Omega_{L(l,0)}=\Omega_{*}\cap L(l,0),\label{4.4.11}\\
  & W_{L(l,0)}=\text{$\mathbf{C}$-span}\, \{v-a\cdot v \mid v\in
  L(l,0),\, a\in A\}\subset L(l,0),\label{4.4.12}\\
  & \Omega_{L(l,0)}^{A}=\Omega_{L(l,0)}/\big(\Omega_{*}\cap
  W_{L(l,0)}\big)=\big(\Omega_{L(l,0)}+W_{L(l,0)}\big)/W_{L(l,0)}.
\label{4.4.13}
\end{align}
Note that by \eqref{3.2.22}, \eqref{4.4.1}, and
\eqref{4.4.11}, one has in fact that 
$\Omega_{L(l,0)}=\Omega_{2k+1}$. Since $L(l,0)$ has an $A$-stable and 
$\hat{\mathfrak{h}}_{*}$-stable complement in
$V_{L}$, it follows that $\tilde{\mathfrak{h}}_{*}'$ acts naturally} on 
$L(l,0)/W_{L(l,0)}$, so that $\Omega_{L(l,0)}^{A}$ may be described 
alternatively as
{\allowdisplaybreaks
\begin{align}
  \Omega_{L(l,0)}^{A} & = \{v\in L(l,0)/W_{L(l,0)} \mid h(n)v=0 \text{
    for } h\in \mathfrak{h}_{*},\, n>0 \}\nonumber\\
  & = \{v\in L(l,0)/W_{L(l,0)} \mid h(n)v=0 \text{
    for } h\in \mathfrak{s},\, n>0 \}\subset \Omega_{*}^{A},\nonumber
\end{align}
where for $h\in \mathfrak{s}$, $h(n)$ denotes the operator $\pi
(h\otimes t^{n})=h_{*}(n)$ of \eqref{4.2.15} and}
\begin{equation*}
  \Omega_{*}^{A}=\{v\in V_{L}/W_{V_{L}} \mid h(n)v=0 \,\text{ for }\, h\in 
\mathfrak{h}_{*},\, n>0 \},
\end{equation*}
with $W_{V_{L}}=\text{$\mathbf{C}$-span}\, \{v-a\cdot v \mid v\in
  V_{L},\, a\in A \}\subset V_{L}$. One may then define a $G$-gradation on 
$\Omega_{L(l,0)}^{A}$ in the following way: note first that 
by \eqref{4.4.2} and \eqref{4.4.8}-\eqref{4.4.9} $\Omega_{*}$ 
decomposes as
\begin{equation}\label{4.4.14}
  \Omega_{*}=\coprod_{\lambda \in
    l^{-1}Q}\Omega_{*}^{\lambda_{*}},\, 
  \Omega_{*}^{\lambda_{*}}=\{v\in \Omega_{*} \mid h(0)v=\langle
  h_{*}\,,\lambda_{*} \rangle v,\,  h\in \mathfrak{s}
  \}= \!\sum_{a\in \hat{L},\,
    \bar{a}''=\lambda_{*}}\!S\big((\tilde{\mathfrak{h}}_{*}^{\perp})_{-}'\big)
\otimes \iota (a).
\end{equation}
Since $\Omega_{L(l,0)}$ is
stable under the operators $h(0)$, $h\in \mathfrak{s}$, we see that
  $\Omega_{L(l,0)}=\coprod_{\lambda \in
    l^{-1}Q}\Omega_{L(l,0)}^{\lambda_{*}}$, where 
$\Omega_{L(l,0)}^{\lambda_{*}}=\Omega_{L(l,0)}\cap
  \Omega_{*}^{\lambda_{*}}=L(l,0)\cap \Omega_{*}^{\lambda_{*}}$. 
It then follows from \eqref{4.4.13} that $\Omega_{L(l,0)}^{A}$ is $G$-graded,  
with $$\big(\Omega_{L(l,0)}^{A}\big)^{g}=\!\!\sum_{\lambda \in
    l^{-1}Q,\,\lambda_{*}+2Q_{*}=g}\!\!
    \big(\Omega_{L(l,0)}^{\lambda_{*}}+W_{L(l,0)}\big)/W_{L(l,0)}
    \, \text{ for } \, g\in G.$$
Note that by \eqref{4.4.7} and \eqref{4.4.11}-\eqref{4.4.12} one
has that if $\lambda ,\mu \in l^{-1}Q$ are such that
$\lambda_{*}-\mu_{*}\in 2Q_{*}$ then 
$\big(\Omega_{L(l,0)}^{\lambda_{*}}+W_{L(l,0)}\big)/W_{L(l,0)}=
\big(\Omega_{L(l,0)}^{\mu_{*}}+W_{L(l,0)}\big)/W_{L(l,0)}$, 
so that in particular
\begin{equation}\label{4.4.15}
  \big(\Omega_{L(l,0)}^{A}\big)^{0}=
\big(\Omega_{L(l,0)}^{0}+W_{L(l,0)}\big)/W_{L(l,0)}.
\end{equation}  

Recall now the space $\Omega_{l}^{0}$ from \eqref{3.2.23}
(with $l=2k+1$) and the element $\omega_{2}=\omega'-\omega_{1}\in
\Omega_{l}^{0}$ (cf.~\eqref{3.2.20}-\eqref{3.2.21} and \eqref{3.2.24}). Then 
clearly
\begin{equation}\label{4.4.16}
  \Omega_{l}^{0}= \Omega_{L(l,0)}^{0}.
\end{equation}
Moreover, by using the actions of Theorem 4.2.1, one can rewrite
$\omega'$, $\omega_{1}$, and $\omega_{2}$ respectively as 
  \begin{equation}\label{4.4.17}
      \omega' =\frac{1}{2(l+3)}\sum_{i=1}^{8}a_{i}(-1)a_{\pi
        (i)}(-1)\cdot \iota (1),\, 
      \omega_{1} =\frac{1}{2}\sum_{i=1}^{2}\gamma^{i}(-1)\gamma^{\pi
        (i)}(-1)\cdot \iota (1),\, \omega_{2}=\omega'-\omega_{1},
   \end{equation}  
where $\gamma^{i}=\big(\tfrac{a_{i}}{\sqrt{l}}\big)_{*}\in
\mathfrak{h}_{*}$ for $i=1,2$ (so that $\{\gamma^{1},\gamma^{2}\}$ and
$\{\gamma^{2},\gamma^{1}\}$ are dual bases of $\mathfrak{h}_{*}$) and
$\mathbf{1}$ is identified with $1\otimes \iota (1)$. Then
$\omega_{2}$ becomes exactly the element $\omega_{\mathcal{G}_{l},\,
  \mathcal{H}_{l}}$ of \cite[(14.53)]{DL1}.

We shall need the following
\begin{lemma}
  With the above notations, one has
      $\Omega_{L(l,0)}^{0}\cap W_{L(l,0)}=\{0\}$.
\end{lemma}  
\begin{proof}
  Recall the
  group $A$ defined in \eqref{4.4.7} and the action of $\hat{L}$ on
  $\mathbf{C}\{L\}$ given by \eqref{4.2.1}, and let
  $U=\Omega_{L(l,0)}^{0}\cap W_{L(l,0)}$ and $v\in U$. Since $U\subset
  \Omega_{L(l,0)}\cap W_{L(l,0)}$, it follows from
  \eqref{4.4.2}, \eqref{4.4.11}, and \eqref{4.4.12} that there exist positive 
integers
  $n$, $m_{i}$ ($1\le i\le n$), such that
   \begin{equation}\label{4.4.18}
      v=\sum_{i=1}^{n}\sum_{j=1}^{m_{i}}\big[u_{ij}\otimes \iota
      (b_{ij})-a_{i}\!\cdot \!(u_{ij}\otimes \iota
      (b_{ij}))\big],
     \end{equation}
  where $u_{ij}\in
      S\big((\tilde{\mathfrak{h}}^{\perp}_{*})'_{-}\big)$, $b_{ij}\in
      \hat{L}$, and $a_{i}\in A$. Recall also the linear maps of 
\eqref{4.2.7}, \eqref{4.2.13}, and \eqref{4.4.5} and
  notice that for $h\in \mathfrak{s}$, $n\in \mathbf{Z}$, the operator
  $h(n)$ acts on $v$ as the ``unrelativized'' operator $\pi (h\otimes t^{n})
=h_{*}(n)$
  of \eqref{4.2.15}. Since $a_{i}\!\cdot \!(u_{ij}\otimes \iota
      (b_{ij}))=u_{ij}\otimes \iota
      (a_{i}b_{ij})$, one gets in particular 
{\allowdisplaybreaks
    \begin{equation}\label{4.4.19}
      \begin{split}
      & h(0)(u_{ij}\otimes \iota
      (b_{ij}))=\langle h_{*}\,,\bar{b}_{ij}'' \rangle \,
      u_{ij}\otimes \iota (b_{ij}),\\
      & h(0)\big[a_{i}\!\cdot \!(u_{ij}\otimes \iota
      (b_{ij}))\big]=\langle h_{*}\,,\bar{b}_{ij}''+\bar{a}_{i}'' \rangle \, 
a_{i}\!\cdot \!(u_{ij}\otimes \iota
      (b_{ij}))
      \end{split}
     \end{equation}
   for $h\in \mathfrak{s}$, $1\le i\le n$, $1\le
        j\le m_{i}$.} In terms of the gradation \eqref{4.4.14},
      \eqref{4.4.19} means that the corresponding vectors have degrees
      $\bar{b}_{ij}''$ and $\bar{b}_{ij}''+\alpha^{i}_{*}$
      respectively, where $\alpha^{i}\in 2Q$ are such that
      $a_{i}=e_{\alpha^{i}_{*}}$, $1\le i\le n$. Denote by
      $\lambda_{1},\ldots ,\lambda_{m}\in \frac{1}{l}Q_{*}$ the
      distinct degrees occurring among the vectors in the right-hand
      side of \eqref{4.4.18}. Then one may rewrite \eqref{4.4.18} as
      \begin{equation}\label{4.4.20}
        v=v_{1}+\ldots +v_{m}
      \end{equation}
      with $v_{i}\in \Omega_{*}$ of degree $\lambda_{i}$. If $m=1$ one gets 
that
      $\bar{b}_{ij}''=\bar{b}_{ij}''+\alpha^{i}_{*}=\lambda_{1}$ for
      all $i$ and $j$, hence $\alpha^{i}_{*}=0$ and 
      $a_{i}=e_{0}=1$, $1\le i\le n$, by the choice of the section $e$
      (cf.~\eqref{4.1.10}). Then \eqref{4.4.18} implies that $v=0$ in this 
case. Assume now
      that $m\ge 2$ and choose $h\in \mathfrak{s}$ such that the
      scalars $\langle h_{*}\,,\lambda_{i} \rangle$ are distinct (such
      $h$ exists since the $\lambda_{i}$'s are distinct linear forms on
      $\mathfrak{h}_{*}$). Then using \eqref{4.4.20} and the fact that $v\in
      \Omega_{L(l,0)}^{0}$ one obtains 
      $v_{1}+\ldots +v_{m}=v$ and $\sum_{k=1}^{m}\langle
      h_{*}\,,\lambda_{k} \rangle
      ^{j}v_{k}=h(0)^{j}v=0$ for $1\le j\le m-1$. 
   This is by construction a system of $m\times m$ linear equations with 
nonsingular matrix,
   so that $v_{k}\in \mathbf{C}v\subset U$ for $1\le k\le m$. In
   particular $v_{i}\in \Omega_{L(l,0)}^{0}$ and thus 
   $\lambda_{i}=0$ for all $i$. Therefore,
   $\bar{b}_{ij}''=\bar{b}_{ij}''+\alpha^{i}_{*}=0$ for $1\le i\le n$,
   $1\le j\le m_{i}$, so that $\alpha^{i}_{*}=0$, $1\le i\le n$, and
   consequently $a_{i}=e_{\alpha^{i}_{*}}=1$, $1\le i\le n$. Then 
\eqref{4.4.18} again implies that $v=0$, which completes the proof. 
\end{proof}

Lemma 4.4.1 and \eqref{4.4.15}-\eqref{4.4.16} yield a linear isomorphism
\begin{equation}\label{4.4.21}
  \Omega_{l}^{0} \stackrel{\cong}{\longrightarrow}
  \big(\Omega_{L(l,0)}^{0}+W_{L(l,0)}\big)/W_{L(l,0)}=
  \big(\Omega_{L(l,0)}^{A}\big)^{0}
\end{equation}
by means of which we identify $\Omega_{l}^{0}$ and 
$\big(\Omega_{L(l,0)}^{A}\big)^{0}$. It follows from \eqref{4.2.9}
that
  $Y_{*}^{(l)}(a\cdot v,z)=a\,Y_{*}^{(l)}(v,z)$ for $v\in
  V_{L}$, $a\in A$. (This is true even for $a\in \hat{L}$ such that $\bar{a}\in L\cap \mathfrak{h}_{*}$.) Therefore, the vertex map $Y_{*}^{(l)}$ of
\eqref{4.2.10} 
induces a well-defined linear map
{\allowdisplaybreaks
  \begin{align}
    \overline{Y}_{*}^{(l)}(\cdot\,,z):\, \Omega_{L(l,0)}^{A} & \longrightarrow
    \big(\mbox{End}\,\Omega_{L(l,0)}^{A}\big)\big[\big[z^{1/l},z^{-1/l}\big]
\big]\label{4.4.22}\\
    v+W_{L(l,0)} & \longmapsto
    \overline{Y}_{*}^{(l)}(v+W_{L(l,0)},z)=Y_{*}^{(l)}(v,z)=\sum_{n\in
      \tfrac{1}{l}\mathbf{Z}}\big(v+W_{L(l,0)}\big)_{n}\,z^{-n-1},\nonumber
\end{align}
which is called the quotient
  relative untwisted vertex operator map defined with respect to the
group $A$ (cf.~\cite[Ch.~4 \& 14]{DL1}).} 
\begin{lemma}
  For $v\in
  \Omega_{l}^{0}$ one has
  $\overline{Y}_{*}^{(l)}(v,z)=Y^{(l)}(v,z)=Y(v,z)$ as (series of)
  operators acting on $\Omega_{L(l,0)}^{A}$.
\end{lemma}
\begin{proof}
  Let $v\in \Omega_{l}^{0}$. By
  \eqref{4.4.14} one may write $v=\sum_{i=1}^{k}v_{i}\otimes \iota
  (a_{i})$ for some $k\in \mathbf{Z}_{+}$, $0\neq v_{i}\in
  S\big((\tilde{\mathfrak{h}}^{\perp}_{*})'_{-}\big)$, and $a_{i}\in
  \hat{L}$, $1\le i\le k$, such that the $\iota (a_{i})$'s are linearly independent. Then the relations
 $0=h(0)v=\sum_{i=1}^{k}\langle h_{*}\,,\bar{a}_{i} \rangle \,
  v_{i}\otimes \iota
  (a_{i})$, $h\in \mathfrak{s}$, imply that $\langle h_{*}\,,
\bar{a}_{i} \rangle =0$, $h\in
  \mathfrak{s},\, 1\le i\le k$. 
  This means that $\bar{a}_{i}''=0$, or equivalently $\bar{a}_{i}'=
\bar{a}_{i}$, $1\le i\le k$. It then follows from \eqref{4.2.8} and
  \eqref{4.2.9} that $Y_{*}^{(l)}(v,z)=Y^{(l)}(v,z)$, which combined with
  \eqref{4.4.21}, \eqref{4.4.22}, and Corollary 4.2.3 proves the lemma.
\end{proof}

Recall that $l=2k+1$ and that
$\text{rank}\,\Omega_{l}^{0}=c_{2}(k)$. 
By using the results of \cite[Ch.~14]{DL1} together with
Theorem 4.2.1, Corollary 4.2.3, \eqref{4.4.21}, and Lemma
4.4.2 we get:
{\allowdisplaybreaks
\begin{theorem}
  Let $\overline{Y}_{*}^{(l)}(\omega_{2}+W_{L(l,0)},z)=\sum_{n\in
    \mathbf{Z}}\bar{L}_{*}(n)z^{-n-2}$. Then the structure
    \begin{equation*}
    \Big(\Omega_{L(l,0)}^{A},\,\overline{Y}_{*}^{(l)},\,\iota
    (1)+W_{L(l,0)},\,\omega_{2}+W_{L(l,0)},\,l,\,G,\,
    (\cdot\,,\cdot)\Big)
  \end{equation*}
  is a simple GVOA of rank $c_{2}(k)$ 
  generated by the set $\{i(x_{\alpha})+W_{L(l,0)} \mid \alpha \in
  \Phi 
  \}$ that consists of $\bar{L}_{*}(0)$-eigenvectors with eigenvalue
  $1-l^{-1}$. Moreover, by restricting the map
  $\overline{Y}_{*}^{(l)}$ to the subspace $\big(\Omega_{L(l,0)}^{A}\big)^{0}
\equiv \Omega_{l}^{0}$, the structure 
      $\big(\Omega_{l}^{0},\,\overline{Y}_{*}^{(l)},\,\iota
    (1)+W_{L(l,0)},\,\omega_{2}+W_{L(l,0)}\big)$
      becomes a simple VOA of rank $c_{2}(k)$ such
   that $\overline{Y}_{*}^{(l)}(v,z)=Y(v,z)$ for every $v\in
   \Omega_{l}^{0}$.\hfill $\Box$ 
\end{theorem}}
Note that by \eqref{4.4.12}, one can actually rewrite $\iota
    (1)+W_{L(l,0)}$ and $\omega_{2}+W_{L(l,0)}$ as $\iota (A)$ and
    $\omega_{2}\otimes \iota (A)$ respectively, where $\omega_{2}$ is
    identified with $\omega_{2}\otimes \iota (1)$.

We now focus on the operators $Y_{*}^{\sigma
  ,(l)}(\cdot\,,z)$ introduced in \S 4.3. We show first:
\begin{lemma}
  Let $a\in A$ and $v\in V_{L}$. Then $Y_{*}^{\sigma
  ,(l)}(a\cdot v,z)=Y_{*}^{\sigma
  ,(l)}(v,z)$ as (series of)
  operators acting on $V_{L}^{T}$.
\end{lemma}  
\begin{proof}
  Since $Y_{*}^{\sigma
  ,(l)}(\cdot\,,z)$ is linear, it suffices to prove the lemma for 
$v\in V_{L}$ of
the form $v=v^{*}\otimes \iota (b)$, where $v^{*}\in
S\big(\tilde{\mathfrak{h}}_{-}'\big)$ and $b\in \hat{L}$. Recall
\eqref{4.1.4} and \eqref{4.1.10}. It follows from 
\cite[Proposition 4.6]{DL2} that if 
$a\in A$ (or more generally, if $a\in \hat{L}$ is such that $\bar{a}\in L\cap
\mathfrak{h}_{*}$) and $v\in V_{L}$, then 
$Y_{*}^{\sigma
  ,(l)}(a\cdot v,z)=\epsilon (\bar{a}\,,\bar{b})\,a\,Y_{*}^{\sigma
  ,(l)}(v,z)$, 
where in the right-hand side $a$ is understood to act on $T$ via the
set-theoretic identification between $\hat{L}$ and
$\hat{L}_{\sigma}$ given by \eqref{4.1.8} (see also \eqref{4.3.6}). Thus 
$Y_{*}^{\sigma
  ,(l)}(a\cdot v,z)=a\,Y_{*}^{\sigma
  ,(l)}(v,z)$ 
by \eqref{4.4.7} and \eqref{4.1.4}, and it is
therefore enough to show that $a$ acts as 1 on $T$ in order to
complete the proof. Using \eqref{4.4.7} and the fact that 
$\epsilon (2\alpha \,,\beta)=1$ for
$\alpha$, $\beta \in L$, we see that the
set-theoretic identification \eqref{4.1.8} restricts to an actual group
isomorphism on $A$, so that $T$ becomes an $A$-module. But
$\hat{L}_{\sigma}=\langle \eta \rangle \times L$ and $L$ acts as 1 on
$T$, hence $(1\,,\alpha )$ acts as 1 on $T$ for every $\alpha \in
L$. In particular, $T$ must be a trivial $A$-module, as needed. 
\end{proof}

It follows from Lemma 4.4.4, \eqref{4.4.10}, and the definition of relative
$\sigma$-twisted vertex operators in \S 4.3 that $Y_{*}^{\sigma
  ,(l)}(\cdot\,,z)$ induces a well-defined quotient linear map 
{\allowdisplaybreaks
  \begin{align}
    \overline{Y}_{*}^{\sigma ,(l)}(\cdot\,,z):\, \Omega_{L(l,0)}^{A} & 
\longrightarrow
    \big(\mbox{End}\,V_{L}^{T}\big)\big[\big[z^{1/6l},z^{-1/6l}\big]\big]
\label{4.4.23}\\
    v+W_{L(l,0)} & \longmapsto
    \overline{Y}_{*}^{\sigma ,(l)}(v+W_{L(l,0)},z)=Y_{*}^{\sigma
      ,(l)}(v,z),\nonumber
  \end{align}
where $Y_{*}^{\sigma
      ,(l)}(\cdot\,,z)$ denotes the restriction to $\Omega_{L(l,0)}$ of the 
map $Y_{*}^{\sigma
      ,(l)}(\cdot\,,z)$ of \eqref{4.3.10}. We may therefore define the
    following subalgebra of $\mbox{End}\,V_{L}^{T}$:
 \begin{align}
    \overline{\mathcal{Y}}_{V_{L}^{T}} &
    =\text{$\mathbf{C}$-span}\,\big\{\text{Res}_{z}\,z^{n}
\overline{Y}_{*}^{\sigma
      ,(l)}(v+W_{L(l,0)},z)\mid v\in \Omega_{L(l,0)}^{A},\,n\in
    \tfrac{1}{6l}\mathbf{Z}\big\}\label{4.4.24} \\
    & = \text{$\mathbf{C}$-span}\,\big\{\text{Res}_{z}\,z^{n}Y_{*}^{\sigma
      ,(l)}(v,z)\mid v\in \Omega_{L(l,0)},\,n\in
    \tfrac{1}{6l}\mathbf{Z}\big\}.\nonumber
  \end{align}

Let now $L\big(\Lambda ;\hat{\mathfrak{g}}[\sigma]\big)$ be any level
$l$ standard $\hat{\mathfrak{g}}[\sigma]$-module.} Then 
$\Lambda=k_{0}\Lambda_{0}+(l-2k_{0})\Lambda_{1}$
for some $k_{0}\in \{0,1,\ldots ,k\}$ by \eqref{2.2.4}, so that $L\big(\Lambda
;\hat{\mathfrak{g}}[\sigma]\big)=L_{k_{0},l}
\big(\hat{\mathfrak{g}}[\sigma]\big)$
in the notation of \S 2.1. According to the proof of Corollary 4.3.3, 
$L_{k_{0},l}\big(\hat{\mathfrak{g}}[\sigma]\big)$ can be isomorphically 
embedded into $L\big(\Lambda_{1}
;\hat{\mathfrak{g}}[\sigma]\big)^{\otimes \,l}\cong
V_{L}^{T}$. Furthermore, Remark 3.2.3 and \eqref{4.4.6} imply that the
principal Heisenberg subalgebra $\tilde{\mathfrak{s}}[\sigma]'$
($=\hat{\mathfrak{s}}[\sigma]$) of $\tilde{\mathfrak{g}}[\sigma]$ may
be identified with $\hat{\mathfrak{h}}_{*}[\sigma]$ ($\cong
\hat{\mathfrak{s}}_{*}[\sigma]$). Consequently, \eqref{2.2.10} may be
rewritten as 
\begin{equation}\label{4.4.25}
  L_{k_{0},l}\big(\hat{\mathfrak{g}}[\sigma]\big)=S\big(\hat{\mathfrak{h}}_{*}
  [\sigma]_{-}\big)\otimes
\Omega_{k_{0},l},
\end{equation}
where $\Omega_{k_{0},l}\in \mathcal{O}_{2}(k)$ is the corresponding
vacuum space (see \eqref{2.2.10} and \eqref{O}). Note also that
\eqref{4.3.14} and \eqref{4.4.1} actually imply that 
\begin{equation}\label{4.4.26}
  \Omega_{k_{0},l}=
  L_{k_{0},l}\big(\hat{\mathfrak{g}}[\sigma]\big)\cap \Omega_{*}^{\sigma}.
\end{equation}

For $\phi \in \mathfrak{h}$ let us define the following operators:
\begin{equation}\label{4.4.27}
    E_{(l)}^{\pm}(\phi ,z)=\exp\left[\sum_{n\in \tfrac{1}{6}\mathbf{Z},
    \pm n>0}\!\!\frac{\phi_{(6n)}\otimes
    t^{n}}{ln}z^{-n}\right]=E_{(1)}^{\pm}(\phi /l ,z),
\end{equation}
  where $E_{(1)}^{\pm}(\cdot \,,z)$ is the unrelativized version of
  the operator $E_{*(1)}^{\pm}(\cdot \,,z)$ defined in
  \eqref{4.3.5}.
\begin{lemma}
  The component operators of $Y_{*}^{\sigma ,(l)}(v,z)$, $v\in
  L(l,0)$, preserve each of the standard
  $\hat{\mathfrak{g}}[\sigma]$-modules
  $L_{k_{0},l}\big(\hat{\mathfrak{g}}[\sigma]\big)$, $k_{0}\in
  \{0,1,\ldots ,k\}$.
\end{lemma}
\begin{proof}
  Let $k_{0}\in
  \{0,1,\ldots ,k\}$ and set
  $M=L_{k_{0},l}\big(\hat{\mathfrak{g}}[\sigma]\big)$. By
  \eqref{4.4.3}-\eqref{4.4.4} and \eqref{4.4.11} it suffices 
  to show that $Y_{*}^{\sigma ,(l)}(v,z)$ preserves $M$ for any $v\in
  \Omega_{L(l,0)}$. Recall from Theorem 4.4.3 that
  $$\Omega_{L(l,0)}=\text{$\mathbf{C}$-span}\,\big\{v^{(1)}_{\,n_{1}}\cdots
  v^{(j)}_{\,n_{j}}\cdot \iota (1) \mid v^{(i)}\in \{i(x_{\alpha}) \mid
  \alpha \in \Phi\}\cup \{\iota (1)\}, n_{i}\in \tfrac{1}{l}\mathbf{Z}, 
i=1,\ldots ,j\big\},$$
  where $v^{(i)}_{\,n_{i}}$ is understood as a component operator of
  $Y_{*}^{(l)}(v^{(i)},z)$. Thus it is enough to
  prove that $M$ is invariant under $Y_{*}^{\sigma ,(l)}(v,z)$ for
  $v\in \Omega_{L(l,0)}$ 
  of the form
  \begin{equation}\label{4.4.28}
    v=v^{(1)}_{\,n_{1}}\cdots
  v^{(j)}_{\,n_{j}}\cdot \iota (1),\, \text{ with } \,v^{(i)}\in 
\{i(x_{\alpha}) \mid
  \alpha \in \Phi\},\,n_{i}\in \tfrac{1}{l}\mathbf{Z},\, i=1,\ldots
  ,j.
  \end{equation}
  Let therefore $\alpha^{m}\in \Phi$ be such that
  $v^{(m)}=i(x_{_{\alpha^{m}}})$, $m=1,\ldots ,j$, and set $\alpha
  =\sum_{m=1}^{j}\alpha^{m}\in Q$. It follows from \eqref{4.2.1} and
  \eqref{4.2.8}-\eqref{4.2.9} (see also \eqref{4.4.2}) that one may
  write $v$ as
  \begin{equation}\label{4.4.29}
    v=\sum_{r_{1},\ldots ,r_{j}=1}^{l}u_{r_{1},\ldots ,r_{j}},
  \end{equation}
  where for any $r_{1},\ldots
    ,r_{j}\in \{1,2,\ldots ,l\}$ one has 
  $$u_{r_{1},\ldots ,r_{j}}=v_{r_{1},\ldots ,r_{j}}\!\otimes
    \iota (e_{r_{1},\ldots ,r_{j}}), \text{ with } v_{r_{1},\ldots ,r_{j}}\in
    S\big((\tilde{\mathfrak{h}}^{\perp}_{*})'_{-}\big), \,
    e_{r_{1},\ldots ,r_{j}}=e_{\alpha^{1}_{\,r_{1}}+\ldots +
        \alpha^{j}_{\,r_{j}}},$$ 
  and $\alpha^{m}_{\,r_{m}}$ denotes the $r_{m}$-th copy of
  $\alpha^{m}\in \Phi$ in $L$. 
  Then \eqref{4.4.29} and \eqref{4.3.8} imply that 
  $$\exp (\Delta_{z*})u_{r_{1},\ldots ,r_{j}}=\exp
    (\Delta_{z})u_{r_{1},\ldots ,r_{j}},$$ 
  which combined with \eqref{4.3.9} and \eqref{4.3.7} gives
  \begin{equation}\label{4.4.30}
    Y_{*}^{\sigma,(l)}(u_{r_{1},\ldots ,r_{j}},z)=W_{*}\big(\exp
    (\Delta_{z*})u_{r_{1},\ldots ,r_{j}},z\big)=W_{*}\big(\exp
    (\Delta_{z})u_{r_{1},\ldots ,r_{j}},z\big).
  \end{equation}
  Note next that for any $r_{1},\ldots ,r_{j}\in \{1,2,\ldots ,l\}$ and
  $p\in \{0,1,\ldots ,5\}$ one
  has by \eqref{4.4.8} that
  $$\sigma^{p}(\alpha^{1}_{\,r_{1}}+\ldots +
        \alpha^{j}_{\,r_{j}})'=\sigma^{p}(\alpha^{1}_{\,r_{1}}+\ldots +
        \alpha^{j}_{\,r_{j}})-l^{-1}\sigma^{p}(\alpha_{*}),$$
  and thus
  \begin{equation}\label{4.4.31}
    \begin{split}
      \langle \sigma^{p}(\alpha^{1}_{\,r_{1}}+\ldots +
        &\alpha^{j}_{\,r_{j}})', (\alpha^{1}_{\,s_{1}}+\ldots +
        \alpha^{j}_{\,s_{j}})'\rangle \\
        &=\langle
        \sigma^{p}(\alpha^{1}_{\,r_{1}}+\ldots +
        \alpha^{j}_{\,r_{j}}), (\alpha^{1}_{\,s_{1}}+\ldots +
        \alpha^{j}_{\,s_{j}})\rangle -l^{-1}\langle
        \sigma^{p}(\alpha ), \alpha \rangle
    \end{split}  
  \end{equation}
  for all $r_{1},\ldots ,r_{j},s_{1},\ldots ,s_{j}\in \{1,2,\ldots
  ,l\}$ and
  $p\in \{0,1,\ldots ,5\}$. Recall \eqref{4.3.4} and let 
  $$ C_{1}(\alpha)=\tau
  \big(\tfrac{\alpha_{*}}{l}\big)^{-1},\,C_{2}(\alpha)=\frac{\langle
    \alpha\,,\alpha
    \rangle}{2l},\,C_{3}(\alpha)=6^{C_{2}(\alpha)}C_{1}(\alpha).$$
  It is then clear from \eqref{4.3.4} and \eqref{4.4.31} that for any 
$r_{1},\ldots
        ,r_{j}\in \{1,2,\ldots ,l\},$
  \begin{equation}\label{4.4.32}
    \tau \big((\alpha^{1}_{\,r_{1}}+\ldots +
        \alpha^{j}_{\,r_{j}})'\big)=C_{1}(\alpha)\tau
        \big(\alpha^{1}_{\,r_{1}}+\ldots + 
        \alpha^{j}_{\,r_{j}}\big).
  \end{equation}
  Recall now \eqref{4.3.3} and notice that by \eqref{4.4.27} one has that 
  $\big[E_{(l)}^{\pm}(\alpha_{*},z),\beta'(\sigma ;z)\big]=0$ for any
  $\beta \in \mathfrak{h}$. It then follows from
  \eqref{4.3.4}-\eqref{4.3.6} and \eqref{4.4.31}-\eqref{4.4.32} that
  for any $r_{1},\ldots 
        ,r_{j}\in \{1,2,\ldots ,l\}$,
  \begin{equation}\label{4.4.33}
  Y_{*}^{\sigma,(l)}\big(\iota (e_{r_{1},\ldots
    ,r_{j}}),z\big)=C_{3}(\alpha)z^{C_{2}(\alpha)} 
      E_{(l)}^{-}(\alpha_{*},z)Y^{\sigma,(l)}\big(\iota 
      (e_{r_{1},\ldots
    ,r_{j}}),z\big)E_{(l)}^{+}(\alpha_{*},z).
  \end{equation}
  Furthermore, since $v_{r_{1},\ldots ,r_{j}}\in
  S\big((\tilde{\mathfrak{h}}^{\perp}_{*})'_{-}\big)$ we get from
  \eqref{4.4.33} and \eqref{4.3.7} that
  $$W_{*}(u_{r_{1},\ldots ,r_{j}},z)=C_{3}(\alpha)z^{C_{2}(\alpha)}
      E_{(l)}^{-}(\alpha_{*},z)W(u_{r_{1},\ldots
        ,r_{j}},z)E_{(l)}^{+}(\alpha_{*},z),$$
  where $W(\cdot \,,z)$ stands for the unrelativized version of the
  operator $W_{*}(\cdot \,,z)$ defined in \eqref{4.3.7}. Then \eqref{4.4.30}
 yields 
  \begin{equation}\label{4.4.34}
    Y_{*}^{\sigma,(l)}(u_{r_{1},\ldots
      ,r_{j}},z)=C_{3}(\alpha)z^{C_{2}(\alpha)}E_{(l)}^{-}(\alpha_{*},z)
    Y^{\sigma,(l)}(u_{r_{1},\ldots
      ,r_{j}},z)E_{(l)}^{+}(\alpha_{*},z).
  \end{equation}
  From \eqref{4.4.29} and \eqref{4.4.34} we finally get
  \begin{equation}\label{4.4.35}
    Y_{*}^{\sigma,(l)}(v,z)=C_{3}(\alpha)z^{C_{2}(\alpha)}
    E_{(l)}^{-}(\alpha_{*},z)Y^{\sigma,(l)}(v,z)E_{(l)}^{+}(\alpha_{*},z).
  \end{equation}
  The lemma now follows from \eqref{4.4.35}, Theorem 4.3.1, and
  Corollary 4.3.3. 
\end{proof}

Let $\alpha \in \Phi$ and consider the operator on $V_{L}^{T}$ determined by
\begin{equation}\label{4.4.36}
  Z(\alpha ,z)=E_{(l)}^{-}(\alpha ,z)\,x_{\alpha}(\sigma
  ;z)\,z\,E_{(l)}^{+}(\alpha ,z)\in
  \overline{U\big(\hat{\mathfrak{g}}[\sigma]\big)}\big[\big[z^{1/6},
z^{-1/6}\big]\big],
\end{equation}
where $E_{(l)}^{\pm}(\cdot \,,z)$ are defined as in \eqref{4.4.27} 
 and $\overline{U\big(\hat{\mathfrak{g}}[\sigma]\big)}$ is a certain
completion of $U\big(\hat{\mathfrak{g}}[\sigma]\big)$ (see
\cite{MP1}). Recall the subalgebra  
$\overline{\mathcal{Y}}_{V_{L}^{T}}$ of $\mbox{End}\,V_{L}^{T}$ 
from \eqref{4.4.24} and notice that $Z(\alpha ,z)$ becomes exactly the
$Z$-operator $Z(\alpha ,\zeta)$ defined
in \cite[(3.18)]{LW} if one substitutes $z$ by $\zeta ^{-6}$ in
\eqref{4.4.36}. Then we can prove the following:
\begin{proposition}
  $\overline{\mathcal{Y}}_{V_{L}^{T}}$ acts
  irreducibly on each of the spaces $\Omega_{k_{0},l}$, $k_{0}\in
  \{0,1,\ldots ,k\}$.
\end{proposition}
\begin{proof}
Let $k_{0}\in
  \{0,1,\ldots ,k\}$. Note first that \eqref{4.3.11} and
  \eqref{4.4.23} imply in particular 
that the component operators of 
$\overline{Y}_{*}^{\sigma ,(l)}(v,z)$, $v\in \Omega_{L(l,0)}^{A}$,
preserve the vacuum space 
$\Omega_{*}^{\sigma}$ defined in \eqref{4.4.1}. By Lemma 4.4.5 and
\eqref{4.4.23}, these operators also preserve the
$\hat{\mathfrak{g}}[\sigma]$-module 
  $L_{k_{0},l}\big(\hat{\mathfrak{g}}[\sigma]\big)$, and then it follows from
  \eqref{4.4.26} that $\overline{\mathcal{Y}}_{V_{L}^{T}}$ preserves
  $\Omega_{k_{0},l}$. 

In order to prove that the action of
$\overline{\mathcal{Y}}_{V_{L}^{T}}$ on $\Omega_{k_{0},l}$ is
irreducible, we note that the operators
$\overline{Y}_{*}^{\sigma ,(l)}(i(x_{\alpha}),z)$, $\alpha \in \Phi$,
are in fact closely related to the operators $Z(\alpha ,z)$ defined in 
\eqref{4.4.36}. Indeed, a simplified version of the computations in Lemma
4.4.5 (in this case just a straightforward consequence of
\eqref{4.3.4}-\eqref{4.3.6}, 
\eqref{4.4.8}, Theorem 4.3.1, Remark 4.3.2, and the fact that
$\tau(\alpha_{i})\tau(\alpha_{i}')^{-1}=\tau \big(\tfrac{\alpha_{*}}{l}\big)$
for $i=1,\ldots,l$) yields 
{\allowdisplaybreaks
  \begin{align}
   \pi ^{\sigma}Z(\alpha ,z)&=E_{(l)}^{-}(\alpha_{*} ,z)\,Y^{\sigma
      ,(l)}(i(x_{\alpha}),z)\,z\,E_{(l)}^{+}(\alpha_{*} ,z)\label{4.4.37} \\
    &=\sum_{i=1}^{l}E_{(l)}^{-}
(\alpha_{*} ,z)\,Y^{\sigma , (1)}(\iota
    (e_{\alpha_{i}}),z)\,z\,E_{(l)}^{+}(\alpha_{*} ,z)\notag \\
    &= (6z)^{-1}z\sum_{i=1}^{l}\tau
    (\alpha_{i})E_{(1)}^{-}(-\alpha_{i}+\alpha_{*}/l
    ,z)E_{(1)}^{+}(-\alpha_{i}+\alpha_{*}/l ,z)\notag  \\
    &= (6z)^{\frac{1}{l}}z
\sum_{i=1}^{l}\tau
    (\alpha_{i})E_{*,(1)}^{-}(-\alpha_{i},z)E_{*,(1)}^{+}(-\alpha_{i},z)
(6z)^{-1+\frac{1}{l}}\notag \\
    &=
    \tau \big(\tfrac{\alpha_{*}}{l}\big)6^{-\frac{1}{l}}z^{1-\frac{1}{l}}
\sum_{i=1}^{l}Y_{*}^{\sigma ,(1)}(\iota
    (e_{\alpha_{i}}),z)=\tau \big(\tfrac{\alpha_{*}}{l}\big)6^{-\frac{1}{l}}
z^{1-\frac{1}{l}}Y_{*}^{\sigma
      ,(l)}(i(x_{\alpha}),z),\notag 
   \end{align}
so that the operator $\pi ^{\sigma}Z(\alpha ,z)$ is essentially a
relative $\sigma$-twisted vertex operator.} Denote by 
$\mathcal{Z}_{V_{L}^{T}}'$ the subalgebra of
$\text{End}\,V_{L}^{T}$ generated by the component operators of the
fields $\pi ^{\sigma}Z(\alpha ,z)$, $\alpha \in 
\Phi$.  Clearly, $\mathcal{Z}_{V_{L}^{T}}'\subset 
\overline{\mathcal{Y}}_{V_{L}^{T}}$. It then follows from 
\cite[Proposition 3.1]{LW} and the ``equivalence
theorem'' \cite[Theorems 5.5 and 5.6]{LW} that $\mathcal{Z}_{V_{L}^{T}}'$
 acts irreducibly on $\Omega_{k_{0},l}$, so that the same must be
 true for $\overline{\mathcal{Y}}_{V_{L}^{T}}$. 
\end{proof}

Recall from \eqref{3.2.24}-\eqref{3.2.25} and \eqref{4.4.17} that the
Virasoro 
element $\omega_{2}\in
\Omega_{l}^{0}$ induces a $\frac{1}{6}\mathbf{Z}$-gradation on
$\Omega_{k_{0},l}$ via the operator
$L^{2}(0)=\text{Res}_{z}\,zY^{\sigma}(\omega_{2},z)$, and that the
corresponding $q$-trace 
is given by $f_{_{k_{0},l}}(q)$ (cf.~\eqref{3.2.26} and \eqref{3.2.30}).
{\allowdisplaybreaks
\begin{lemma}
  With the above notations, one has
  \begin{equation}\label{4.4.38}
    \overline{Y}_{*}^{\sigma ,
  (l)}(\omega_{2}+W_{L(l,0)},z)=Y^{\sigma}(\omega_{2},z).
  \end{equation}
  In particular, the $q$-trace
  $$\chi_{_{k_{0},k}}^{\sigma}(q):=\text{{\em
      tr}}_{_{\Omega_{k_{0},l}}}\,q^{\bar{L}_{*}(0)-\tfrac{c_{2}(k)}{24}}\,\,
    \text{ {\em ({\em where} $\bar{L}_{*}(0)=\text{Res}_{z}\,z
    \overline{Y}_{*}^{\sigma ,(l)}(\omega_{2}+W_{L(l,0)},z)$)}}$$
is well defined and coincides with $f_{_{k_{0},l}}(q)$.
\end{lemma}}
\begin{proof}
  One can use the definition of relative $\sigma$-twisted vertex
operators in \S 4.3 and imitate the proof of Lemma
4.4.2 in order to obtain from Corollary 4.3.3 that when acting on 
$\Omega_{*}^{\sigma}$,  
\begin{equation}\label{4.4.39}
  Y_{*}^{\sigma , (l)}(v,z)=Y^{\sigma,(l)}(v,z)=Y^{\sigma}(v,z) \, 
\text{ for }\, v\in \Omega_{l}^{0}.
\end{equation}  
By Proposition 4.4.6, \eqref{4.4.39} holds even as an identity between (series
of) operators on $\Omega_{k_{0},l}$. Then \eqref{4.4.38} follows from
Lemma 4.4.1, \eqref{4.4.21}, \eqref{4.4.23}
and \eqref{4.4.39}.
\end{proof}

Summarizing, we have
{\allowdisplaybreaks
\begin{theorem}
  Let $k_{0}\in \{0,1,\ldots ,k\}$ and $\Omega_{k_{0},l}\in
  \mathcal{O}_{2}(k)$. Then the structure 
  $$\Big(\Omega_{L(l,0)}^{A},\,\overline{Y}_{*}^{(l)},\,\iota
    (A),\,\omega_{2}\otimes \iota
    (A),\,l,\,G,\,
    (\cdot\,,\cdot)\Big)$$
  is a simple GVOA of rank $c_{2}(k)$ which acts irreducibly on
  $\Omega_{k_{0},l}$ by means of the operators
  $\overline{Y}_{*}^{\sigma ,(l)}(\cdot \,,z)$, in such a way that
  $\Omega_{k_{0},l}$ has a well-defined $q$-trace
  $\chi_{_{k_{0},k}}^{\sigma}(q)$ given by
  $$\chi_{_{k_{0},k}}^{\sigma}(q)=f_{_{k_{0},l}}(q)
  =\chi_{_{k_{0},k}}^{\tau}(q),$$
  where $\chi_{_{k_{0},k}}^{\tau}(q)$ is as in Theorem 3.2.7.\hfill $\Box$
\end{theorem}}
\begin{remark}
  It follows from Lemma 4.4.5 and Proposition 4.4.6 that the whole
  VOA $L(l,0)$ acts in fact irreducibly on each of the spaces
  $\Omega_{k_{0},l}$, $k_{0}\in \{0,1,\ldots ,k\}$, by means of the
  operators $\overline{Y}_{*}^{\sigma ,(l)}(\cdot\,,z)$. Furthermore,
  under this action the Virasoro element $\omega' \in L(l,0)$ (see
 \eqref{4.4.17} and \eqref{3.2.20}) induces the same
$q$-trace $\chi_{_{k_{0},k}}^{\sigma}(q)$ for the space
$\Omega_{k_{0},l}$. The techniques of
Lemma 4.4.5 and Proposition 4.4.6 
seem therefore flexible enough to allow the construction of
structures of (G)VOA type which are ``larger'' than $\Omega_{L(l,0)}^{A}$
and still satisfy the 
properties described in Theorem 4.4.8 (except possibly for the rank).
\end{remark}  

We conclude our partial discussion of Problems 1 and 2 with some remarks
about potential further developments. A natural context for Theorem
4.4.8 would obviously be the ``twisted
representation theory'' of the GVOA $\Omega_{L(l,0)}^{A}$. This fact
alone motivates therefore an axiomatic study of the notion of 
``twisted module for a GVOA'', of which the spaces $\Omega_{k_{0},l}$
should be natural examples. It is quite clear that the techniques
used for proving Lemma 4.4.5 and Proposition 4.4.6 are easily
adapted to the more 
general setting of an arbitrary positive
definite even lattice, thus extending to a reasonably large class of
examples. We therefore believe that formulas such as 
\eqref{4.4.35} and \eqref{4.4.37} together with the duality
pro\-perties of ordinary twisted modules could lead to a relatively 
simple-looking Jacobi identity for a twisted GVOA-module, which is
usually the main axiom for structures of this type. (The very
  technical nature of the 
  generalized twisted Jacobi identity 
  established in \cite{DL2} makes this identity somewhat 
  inappropriate for 
  such purposes.) The
appropriate axiomatic setting once developed, one expects the
GVOA $\Omega_{L(l,0)}^{A}$ to be ``$\sigma$-rational'' and the spaces
$\Omega_{k_{0},l}$, $k_{0}\in \{0,1,\ldots ,k\}$, to be in fact its only
simple $\sigma$-twisted modules, in view of the
``equivalence theorem'' \cite[Theorems 5.5 and 5.6]{LW}.
Theorem 3.2.7, Theorem 4.4.8, and Remark 3.2.8 (ii) would then suggest a strong
connection between the 
$\tau$-twisted representation theory of the $\tau$-rational VOA
$V_{k}\big(A_{1}^{(1)}\big)^{\otimes \,2}$ and the ``$\sigma$-twisted
representation theory'' of the ``$\sigma$-rational'' GVOA
$\Omega_{L(l,0)}^{A}$, or between appropriate substructures 
of these algebras. We believe that such  
connections would eventually answer Problems 1 and 2 in the affirmative.

\section*{Appendix: Some character formulas}

We give here the formula for the unspecialized character of
$V=L(3\Lambda_{0};A_{2}^{(1)})$. Recall the notations from \S 2.1
in this case. In particular, $\alpha_{0}$, $\alpha_{1}$, $\alpha_{2}$
are the simple roots and
$\delta=\alpha_{0}+\alpha_{1}+\alpha_{2}$. Set $u=e^{-\alpha_{0}}$,
$v=e^{-\alpha_{1}}$, $w=e^{-\alpha_{2}}$. The computations were made
by using the techniques developed in \cite{KP} and \cite{K1}, that is, by
expressing $\mbox{ch}\,V$ in terms of the string functions for the
maximal dominant weights and explicitly calculating their orbits
under the affine Weyl group. One gets
{\allowdisplaybreaks 
\begin{align*}
  e^{-3\Lambda_{0}}\mbox{ch}V=
 & \bigg(\sum_{k=0}^{\infty}\mbox{dim}V_{3\Lambda_{0}-k\delta}(uvw)^{k}\bigg)
\Bigg[\sum_{m,n\in \mbox{{\bf Z}}}(uvw)^{3(m^{2}+n^{2}-mn)}v^{-3m}w^{-3n}
\Bigg]\tag{A.1} \\
  & +\bigg(\sum_{k=0}^{\infty}\mbox{dim}V_{3\Lambda_{0}-\alpha_{0}-k\delta}
(uvw)^{k+1}\bigg)\Bigg[\sum_{m,n\in \mbox{{\bf Z}}}(uvw)^{3(m^{2}+n^{2}-mn)-
(2m-n)}v^{1-3m}w^{-3n} \\
  & \quad \times \big(1+(uvw)^{2(2m-n)}v^{-2}\big)+ \sum_{m,n\in
    \mbox{{\bf 
        Z}}}(uvw)^{3(m^{2}+n^{2}-mn)-(2n-m)}v^{-3m}w^{1-3n} \\
  & \quad \times \big(1+(uvw)^{2(2n-m)}w^{-2}\big)+\sum_{m,n\in
    \mbox{{\bf 
        Z}}}(uvw)^{3(m^{2}+n^{2}-mn)-(m+n)}v^{1-3m}w^{1-3n} \\
  & \quad \times\big(1+(uvw)^{2(m+n)}v^{-2}w^{-2}\big)\Bigg]+
\bigg(\sum_{k=0}^{\infty}\mbox{dim}V_{3\Lambda_{0}-2\alpha_{0}-\alpha_{2}-
k\delta}(uvw)^{k}\bigg) \\
  & \quad \times\Bigg[\sum_{m,n\in \mbox{{\bf Z}}}(uvw)^{3(m^{2}+
n^{2}-mn-m-n)+2}
  v^{1+3(2n-m)}w^{-1+3n}\Bigg] \\
  & +\bigg(\sum_{k=0}^{\infty}\mbox{dim}V_{3\Lambda_{0}-2\alpha_{0}-
\alpha_{1}-k\delta}(uvw)^{k}\bigg) \\
  & \quad \times\Bigg[\sum_{m,n\in \mbox{{\bf Z}}}
(uvw)^{3(m^{2}+n^{2}-mn+m-n)+2}v^{-1-3m}w^{1-3n}\Bigg].
\end{align*}}

Using the table Affine $A_{2}$-Level 3-Class 0-Rp 4 from 
\cite[p.~408]{KMPS}, one can obtain concrete expressions for
every specialization
$F_{(s_{0},s_{1},s_{2})}(e^{-3\Lambda_{0}}\mbox{ch}V)$ ``near the
top'' of the module $V$. For instance,
\begin{equation*}
  F_{(4,1,1)}(e^{-3\Lambda_{0}}\mbox{ch}V)\big|_{q\rightarrow 
    q^{\frac{1}{6}}}=1+q^{\frac{4}{6}}+2q^{\frac{5}{6}}+2q+2q^{\frac{7}{6}}+
\ldots+46q^{\frac{20}{6}}+\mathcal{O}(q^{\frac{21}{6}}).\tag{A.2}
\end{equation*}
Unfortunately, as shown by computer experiments, it is not likely that
(A.2) could be rewritten as a product formula in which
the graded dimension of a symmetric algebra is easily identifiable 
(cf. Remark 3.2.6). The same can actually be said about the character
\begin{equation*}
  F_{(1,0,0)}(e^{-3\Lambda_{0}}\mbox{ch}V)=1+8q+44q^{2}+192q^{3}+\ldots+
\mathcal{O}(q^{23}),
\end{equation*}  
which up to the factor $q^{-\frac{1}{6}}$ is just the $q$-trace of the
VOA $V_{3}(A_{2}^{(1)})$. The latter formula may of course be used for
computing the graded dimension of the vacuum space
$\Omega_{3}$. Indeed, after dividing out the graded dimension of the
Fock space $M(3)$ (which by \eqref{3.2.18}  
is $\prod_{n=1}^{\infty}(1-q^{n})^{-2})$ one gets
\begin{equation*}
  \dim_{*} \Omega_{3}=1+6q+27q^{2}+98q^{3}+\ldots+\mathcal{O}(q^{23}).
\end{equation*}

For the purposes described in \S 3.2, we are in fact more interested in
comparing the $q$-traces of the VOAs $\Omega_{3}^{0}$ and
$V_{1}(A_{1}^{(1)})^{\otimes 2}$. If we substitute 
$u$ by $qv^{-1}w^{-1}$ in formula (A.1), then  
multiply it with $\prod_{n=1}^{\infty}(1-q^{n})^{2}$ and collect the
terms containing only zero powers of $v$ and $w$, we get 
\begin{equation*}
  \mbox{tr}_{_{\Omega_{3}^{0}}}\,q^{L^{2}(0)-\frac{c_{2}(1)}{24}}=
q^{-\frac{1}{12}}(1+3q^{2}+8q^{3}+16q^{4}+\ldots+\mathcal{O}(q^{22})),\tag{A.3}
\end{equation*}
since $c_{2}(1)=2$ (cf.~\S 3.2). On the other hand, using the
formulas of \cite[\S 21.8]{KMPS} for the unspecialized character of
$L(\Lambda_{0};A_{1}^{(1)})=V_{1}(A_{1}^{(1)})$, one gets (recall from \S 3.2 
that $c_{1}(1)=1$)
\begin{equation*}
  \mbox{tr}_{V_{1}(A_{1}^{(1)})}\,\,q^{L(0)-\frac{c_{1}(1)}{24}}=
q^{-\frac{1}{24}}(1+3q+4q^{2}+7q^{3}+\ldots+\mathcal{O}(q^{23})),\tag{A.4}
\end{equation*}
and consequently
\begin{equation*}
  \mbox{tr}_{V_{1}(A_{1}^{(1)})^{\otimes
      2}}\,\,q^{L_{_{\otimes}}(0)-\frac{c_{2}(1)}{24}}=q^{-\frac{1}{12}}
(1+6q+17q^{2}+38q^{3}+\ldots+\mathcal{O}(q^{23})).\tag{A.5}
\end{equation*}  
  
As we already noted in \S 3.2, the conformal vectors of $\Omega_{3}^{0}$ and
$V_{1}(A_{1}^{(1)})$ are such that
$\mbox{rank}\,\Omega_{3}^{0}=2\,\mbox{rank}\,V_{1}(A_{1}^{(1)})$,
which obviously prevents
these spaces from being isomorphic as VOAs. Formulas (A.3) and (A.4)
show actually that $\Omega_{3}^{0}$ and
$V_{1}(A_{1}^{(1)})$ are not even isomorphic as graded vector
spaces. Moreover, (A.3) and (A.5) imply that although they have the
same rank, the VOAs $\Omega_{3}^{0}$ and
$V_{1}(A_{1}^{(1)})^{\otimes 2}$ cannot be isomorphic either.

\end{document}